\documentclass[11pt]{article}

\usepackage{amsmath,amssymb,amsthm}
\usepackage{graphicx,verbatim}
\usepackage{algorithm,algorithmic}
\usepackage[small,bf]{caption}
\usepackage{subcaption} 
\usepackage{cases}
\usepackage[table,dvipsnames]{xcolor}
\newif\ifshowchanges
\showchangestrue
\ifshowchanges
\else
\fi
\usepackage[shortlabels]{enumitem}
\usepackage{hyperref,url}
\usepackage{cite}
\usepackage{booktabs,adjustbox,multirow}  
\usepackage{makecell}
\usepackage{import} 
\usepackage{aligned-overset} 
\usepackage[toc,page]{appendix}
\usepackage[normalem]{ulem}
\hypersetup{
	colorlinks=true,
	linkcolor=red,     
	urlcolor=magenta,
	citecolor={blue},
}
\usepackage{tikz,psfrag,pgfplots}
\pgfplotsset{compat=newest}
\usetikzlibrary{plotmarks,arrows.meta}
\usepgfplotslibrary{patchplots,external}
\usepackage{lmodern}
\usepackage{dsfont} 
\usepackage{bm,bbm} 
\usepackage{upgreek} 
\usepackage{pifont}

\newtheorem{lemma}{{Lemma}}
\newtheorem{theorem}{{Theorem}}
\newtheorem{corollary}{{Corollary}}

\newtheorem{assumption}{{Assumption}}
\newtheorem{remark}{{Remark}}

\newtheorem{example}{{Example}}

\newcommand{\grad}{\nabla}              
\newcommand{\zero}{\mathbf{0}}          
\newcommand{\one}{\mathbf{1}}           

\newcommand{\real}{\mathbb{R}}          

\newcommand{\diag}{\operatorname{diag}}         
\newcommand{\col}{\operatorname{col}}           

\newcommand{\Ex}{\mathds{E}}                    

\newcommand{\minimize}{\mathop{\text{minimize}}}       

\newcommand{\tran}{^{\textit{\footnotesize\texttt{T}}}} 

\newcommand{\define}{\triangleq}            

\def\bxi{\boldsymbol{\xi}}                  

\def\B{\mathbf{B}}                          
\def\C{\mathbf{C}}                          
\def\D{\mathbf{D}}                          
\def\F{\mathbf{F}}                          
\def\I{\mathbf{I}}                          
\def\T{\mathbf{T}}                          
\def\W{\mathbf{W}}                          

\def\e{\mathbf{e}}                          
\def\f{\mathbf{f}}                          
\def\g{\mathbf{g}}                          
\def\t{\mathbf{t}}                          
\def\v{\mathbf{v}}                          
\def\x{\mathbf{x}}                          
\def\y{\mathbf{y}}                          

\newcommand{\cF}{\mathcal{F}}               
\newcommand{\cO}{\mathcal{O}}               

\def\evdots{%
	\vbox{%
		\baselineskip=2pt
		\lineskiplimit=0pt
		\kern6pt
		\hbox{$.$}\hbox{$.$}\hbox{$.$}%
	}%
}
\usepackage[margin=1in,
includefoot,
headsep=20pt 
]{geometry}
\usepackage{sectsty}
\allsectionsfont{\sffamily}
\usepackage{fancyhdr} 
\AtBeginDocument{\ifshowchanges\else\hypersetup{linkcolor=black,citecolor=black,urlcolor=black}\fi}
\begin{document}
	
	\title{\bfseries   On the Spectral Dependence of Decentralized Stochastic Gradient Tracking}
	\author{%
		\small
		\begin{tabular}{@{}cc@{}}
			\begin{minipage}[t]{0.47\textwidth}
				\centering
				{\sffamily	Sulaiman A. Alghunaim}\\
				{\footnotesize Department of Electrical Engineering, Kuwait University, Kuwait}\\
				{\tt\footnotesize sulaiman.alghunaim@ku.edu.kw}
			\end{minipage}
			&
			\begin{minipage}[t]{0.47\textwidth}
				\centering
				{\sffamily	Haoyuan Cai}\\
				{\footnotesize   School of Engineering, \'Ecole Polytechnique F\'ed\'erale de Lausanne (EPFL), Switzerland}\\
				{\tt\footnotesize haoyuan.cai@epfl.ch}
			\end{minipage}
			\\
			~& \\
			\multicolumn{2}{c}{
				\begin{minipage}[t]{0.47\textwidth}
					\centering
					{\sffamily	Ali H. Sayed}\\
					{\footnotesize   School of Engineering, \'Ecole Polytechnique F\'ed\'erale de Lausanne (EPFL), Switzerland}\\
					{\tt\footnotesize ali.sayed@epfl.ch}
				\end{minipage}
			}
		\end{tabular}
	}

	\date{}
	\maketitle
	\thispagestyle{fancy} 

	\begin{abstract}
		We study how the spectral geometry of the mixing matrix enters the
		convergence of decentralized stochastic gradient tracking (GT) for smooth
		strongly convex problems.  General bounds
		contain the factor $\Delta_\lambda^{-3}$ in the higher-order stochastic
		term, where $\Delta_\lambda$ is the absolute spectral gap. We show that, even for two
		agents with identical quadratic objectives, the exact steady-state
		mean-square disagreement scales as $\alpha^2\sigma^2\Delta_\lambda^{-3}$
		as the nontrivial eigenvalue approaches $-1$, provided the stepsize
		$\alpha$ is sufficiently small relative to the gap. Here $\sigma^2$ is
		the gradient-noise variance. By unrolling the disagreement recursion, we establish
		a stochastic bound that averages all nontrivial eigenvalue contributions
		for diagonalizable mixing matrices, up to an eigenvector conditioning
		factor. As a real or complex eigenvalue approaches the unit circle away
		from $1$, its contribution grows as the inverse cube of its distance to
		the circle; growth can be slower near $1$. For symmetric mixing with
		gaps $p$ at $1$ and $\nu$ at $-1$, the weighted sum
		$\theta_1/p+\theta_2/\nu^3$ replaces the prior factor
		$(\Delta_\lambda\nu^2)^{-1}$. For $N$ agents, $N\theta_1$ and
		$N\theta_2$ count the nonnegative and negative nontrivial eigenvalues,
		respectively. On a symmetric family, a finite-horizon lower bound for
		the average iterate matches the worst-case factor $p^{-1}+\nu^{-3}$
		up to universal constants. Numerical experiments confirm the predicted
		scaling.
	\end{abstract}

	\section{Introduction}
	\label{sec:intro}
	We consider a decentralized network of $N$ agents, where each agent can
	communicate only with its direct neighbors. The agents aim to cooperatively
	solve the distributed stochastic optimization problem
	\begin{equation}\label{min_problem}
		\minimize_{x\in\real^d}\quad
		f(x)\define\frac{1}{N}\sum_{i=1}^{N}f_i(x),
		\quad
		f_i(x)\define\Ex[F_i(x;\xi_i)],
	\end{equation}
	where $f_i:\real^d\rightarrow\real$ is the local objective function of
	agent $i$. We assume that each $f_i$ is smooth and strongly convex and is
	given as the expectation of a loss function $F_i$ with respect to the local
	random data $\xi_i$. At each iteration, agent $i$ queries a local stochastic gradient $\grad F_i(x;\xi_i)$ and exchanges iterates with its neighbors through a mixing matrix to solve \eqref{min_problem} without using a central coordinator.
	
	\subsection{Motivation and related work}
	Let $\lambda\in[0,1)$ denote the mixing rate of the network and
	$\Delta_\lambda\define1-\lambda$ the corresponding gap.
	For smooth, strongly convex (or PL) objectives, bounds for decentralized
	stochastic gradient descent (DSGD) \cite{koloskova2020unified},
	gradient tracking (GT) \cite{pu2021distributed,xin2021improved,koloskova2021improved,alghunaim2021unified,
		alghunaim2023enhanced},
	and exact diffusion \cite{yuan2021removing} can attain a leading
	stochastic term of order $\sigma^2/(NK)$ with stepsizes tuned to
	the horizon $K$ or with diminishing stepsizes, up to logarithmic and problem-dependent factors.
	Here $\sigma^2$ bounds the gradient-noise variance at each agent;
	the $1/N$ reduction from averaging independent noise across agents
	gives linear speedup.
	
	At a fixed iteration budget, the network term can mask the $1/N$ gain
	from using more agents. If this term scales as
	$\sigma^2/(\Delta_\lambda^q K^2)$, then $K\gtrsim N\Delta_\lambda^{-q}$
	suffices for it to fall below the leading term $\sigma^2/(NK)$.
	For example, when $\Delta_\lambda=10^{-3}$, the predicted crossover
	is $10^3$ times larger for $q=2$ and $10^6$ times larger for $q=3$
	than for $q=1$. Thus bounds with the same asymptotic linear speedup
	can give very different guarantees on how many iterations suffice
	to see it. The full transient also depends on the initialization
	and the admissible stepsize.

	The network dependence of decentralized gradient descent and DSGD is well understood \cite{yuan2016convergence,lian2017can,yuan2020influence,koloskova2020unified,vogels2023beyond,takezawa2026improved}.
	For symmetric mixing, averaging over the full spectrum sharpens its
	stochastic network bound \cite{takezawa2026improved}. Earlier
	small-stepsize analyses of diffusion LMS obtained full-spectrum
	approximations of the steady-state mean-square deviation for a common
	minimizer and common regressor covariance
	\cite{tu2012diffvscons,tu2013influence}.
	With exact gradients and a sufficiently small constant stepsize,
	GT \cite{di2016next,xu2015augmented} removes the steady-state bias that heterogeneous local objectives
	can cause in decentralized gradient descent
	\cite{xu2015augmented,qu2017harnessing,nedic2017achieving}.
	Its convergence has been extensively studied under various settings; see, e.g.,
	\cite{lu2019gnsd,xin2018linear,pu2020push,xin2020variance,xin2020general}
	and the references therein.

	Despite this vast  literature, the precise dependence of GT on the
	mixing matrix is not yet well established. For smooth
	strongly convex (or PL) objectives, 
	this dependence varies across existing analyses of GT
	\cite{pu2021distributed,xin2021improved,alghunaim2021unified,
		alghunaim2023enhanced,koloskova2021improved}.
	Classical analyses give a stochastic error term of order
	$\alpha^2\Delta_\lambda^{-3}$ \cite{pu2021distributed,xin2021improved},
	that is, a network factor $\Delta_\lambda^{-3}$ in the coefficient of
	the $K^{-2}$ term once $\alpha$ is tuned to $K$. When the mixing matrix is
	symmetric and positive semidefinite, \cite{alghunaim2023enhanced}
	improves this dependence to $\Delta_\lambda^{-1}$. For symmetric mixing
	matrices that are not necessarily positive semidefinite,
	\cite{koloskova2021improved} obtains $\Delta_\lambda^{-1}\nu^{-2}$,
	where $\nu$ measures the distance of the spectrum from $-1$, which implies that eigenvalues near $-1$ worsen the stochastic bound.
	
	This raises three 
	questions. Is the cubic gap dependence unavoidable in the worst case
	for the GT recursion studied here? Is the multiplicative factor
	$\Delta_\lambda^{-1}\nu^{-2}$ necessary for symmetric mixing?
	For complex spectra, which eigenvalue locations give cubic noise
	amplification, and which give milder dependence? We answer these questions through a spectral analysis that covers	general diagonalizable doubly stochastic mixing matrices.
	
	\subsection{Contributions}
	\label{sec:contributions}
	
	Our contributions are the following.
	
	\begin{itemize}
		\item \emph{Sharpness of the cubic network dependence}.
		For identical quadratic objectives on two agents, we compute the
		steady-state disagreement exactly and show that it scales as
		$\alpha^2\sigma^2\Delta_\lambda^{-3}$
		under the stepsize condition \eqref{eq:lb_stepsize}
		(Theorem~\ref{thm:cubic_network_lower_bound}). When the two spectral gaps coincide, Theorem~\ref{thm:jlb_joint}
		gives a finite-horizon lower bound with the same cubic network factor
		for the mean-square error of the average iterate. The construction
		applies to every $N\geq4$ divisible by four, with constants
		independent of $N$.
		
		\item \emph{Exact propagation analysis}.
		We write the disagreement as an exact convolution of the past
		optimization and stochastic-gradient perturbations
		(Lemma~\ref{lemma:gt_nonsym_convolution}) and compute in closed form
		how much noise each (complex) eigenmode accumulates
		(Lemma~\ref{lem:gt_nonsym_modal_energy}). This gives sharper network
		terms while keeping the leading $\sigma^2/(NK)$ term, up to
		logarithmic and problem-dependent factors. Nonorthogonal eigenvectors
		enter only through a condition number $\kappa_W$, which equals one
		for normal matrices (Corollary~\ref{cor:gt_nonsym_normal}).

		\item \emph{Spectral characterization of noise amplification}.
		The closed form shows that the noise contribution of an eigenvalue
		grows cubically in its inverse distance to the unit circle when it
		stays away from $1$. Thus, the cubic penalty is not specific to
		eigenvalues near $-1$: it can arise near any point of the unit circle
		other than $1$. Approaches toward $1$ can give milder amplification. For diagonalizable mixing matrices, the noise bound averages the
		contributions of all nontrivial eigenvalues.
		On the directed ring with fixed laziness, using the average instead
		of the largest contribution at each lag sharpens the bound by a
		factor proportional to the number of agents
		(Example~\ref{ex:gt_nonsym_directed_ring}).
		
		\item \emph{Additive spectral dependence}.
		For real or complex spectra,  two spectral regions contribute
		additively to the stochastic bound
		(Corollary~\ref{cor:gt_nonsym_two_region_rate}). For symmetric mixing matrix $W$, the resulting network factor is
		$\theta_1/p+\theta_2/\nu^3$, where
		$p\define1-\max\{\lambda_2(W),0\}$ and
		$\nu\define1+\min\{\lambda_N(W),0\}$, so that
		$\Delta_\lambda=\min\{p,\nu\}$. Here $N\theta_1$ and $N\theta_2$
		count the nonnegative and negative nontrivial eigenvalues, respectively.
		Thus a single negative eigenvalue contributes only $1/(N\nu^3)$
		to this bound. Even without the weights, the additive factor $p^{-1}+\nu^{-3}$
		is at most twice $(\Delta_\lambda\nu^2)^{-1}$, which is equivalent
		up to absolute constants to the network factor in
		\cite{koloskova2021improved}. The additive factor can be much smaller
		when $p\ll\nu\ll1$ (Lemma~\ref{lemma:cmp_ratio}). Theorem~\ref{thm:jlb_joint} establishes
		a matching finite-horizon lower bound for this additive dependence
		on a family of symmetric networks. (The lower bounds concern algorithm \eqref{gt_atc_semi} under the stated	initialization and constant-stepsize conditions; they are not
		oracle-complexity lower bounds over a class of decentralized algorithms.)
	\end{itemize}

	\paragraph{Organization.}
	Section~\ref{sec:gt_semi_problem} introduces the algorithm and
	error dynamics. Sections~\ref{sec:gt_baseline_result}--%
	\ref{sec:sharp_network_dependence} establish the convergence
	and spectral lower bounds.
	Sections~\ref{sec:comparison} and~\ref{sec:numerics} present
	comparisons and experiments.
	Appendix~\ref{app:shared_estimates} collects shared estimates
	and tuning, Appendices~\ref{app:cubic_network_lower_bound}--%
	\ref{app:finite_lower_bound} prove the new results, and
	Appendix~\ref{sec:gt_two_recursion_appendix} proves the baseline theorem.

	\paragraph{Notation.}
	Vectors are columns. The norms $\|\cdot\|$, $\|\cdot\|_2$, and
	$\|\cdot\|_F$ are Euclidean, spectral, and Frobenius, respectively. The symbols
	$(\cdot)\tran$ and $\Ex$ denote transpose and expectation, respectively,
	while $I$ or $\I$ denotes the identity matrix, $\zero$ denotes
	a zero vector or matrix, and $\one$ denotes the all-ones vector,
	with dimensions determined by context. We use $\otimes$ for the Kronecker
	product and $\col\{\cdot\}$ for column stacking. We write
	$\kappa\define L/\mu$ for the condition number.  We write $a_K=\cO(b_K)$ if $a_K\leq Cb_K$ for some constant $C>0$,
	and $a_K=\Theta(b_K)$ if both $a_K=\cO(b_K)$ and
	$b_K=\cO(a_K)$. The notation $\widetilde\cO(\cdot)$ additionally hides
	polylogarithmic factors. We use $a\lesssim b$ and $a\gtrsim b$ to mean
	$a\leq Cb$ and $a\geq cb$, respectively, for universal constants
	$C,c>0$.
	For symmetric doubly stochastic $W$, its eigenvalues are ordered as
	$1=\lambda_1(W)\geq\lambda_2(W)\geq\cdots\geq\lambda_N(W)$.
	We write $\imath$ for the imaginary unit, with $\imath^2=-1$,
	and use $\operatorname{Re}(z)$ and $\operatorname{Im}(z)$
	for the real and imaginary parts of $z$.
	
	\section{Algorithm, assumptions, and error dynamics}
	\label{sec:gt_semi_problem}
	\label{sec:gt_semi_dynamics_main}
	
	\subsection{Gradient-tracking recursion}
	We describe the gradient-tracking algorithm in network form, which is convenient
	for the subsequent analysis. Let $W=[w_{ij}]\in\real^{N\times N}$ denote the
	mixing matrix associated with the communication graph, where $w_{ij}\geq0$
	determines the weight that agent $i$ assigns to information received from agent
	$j$, and $w_{ij}=0$ whenever agents $i$ and $j$ are not neighbors. Let
	$x_i\in\real^d$ denote the local estimate maintained by agent $i$ of the
	optimization variable $x$ in \eqref{min_problem}.
	
	We introduce the stacked network quantities
	\begin{subequations}\label{gt_semi_stacked_notation}
		\begin{alignat}{2}
			\x&\define\col\{x_i\}_{i=1}^N\in\real^{Nd},&\enspace
			\grad\f(\x)&\define\col\{\grad f_i(x_i)\}_{i=1}^N,\\
			\bxi&\define\col\{\xi_i\}_{i=1}^N,&\enspace
			\grad\F(\x;\bxi)&\define\col\{\grad F_i(x_i;\xi_i)\}_{i=1}^N,\\
			\W&\define W\otimes I_d\in\real^{Nd\times Nd},&&
		\end{alignat}
	\end{subequations}
	Using this notation, the stochastic gradient-tracking method considered in this
	paper takes the form
	\begin{subequations}\label{gt_atc_semi}
		\begin{align}
			\x^{k+1}
			&=
			\W(\x^k-\alpha\g^k),
			\label{gt_atc_semi_x}
			\\
			\g^{k+1}
			&=
			\W\g^k
			+\grad\F(\x^{k+1};\bxi^{k+1})
			-\grad\F(\x^k;\bxi^k),
			\label{gt_atc_semi_g}
		\end{align}
	\end{subequations}
	for $k\geq0$, where $\alpha>0$ is the stepsize and $
	\g^0=\grad\F(\x^0;\bxi^0)$.
	The variable $\g^k$ tracks the aggregate gradient information across the
	network, while the mixing step promotes agreement among the local iterates.
	Recursion \eqref{gt_atc_semi} is the same gradient-tracking formulation
	studied in \cite{koloskova2021improved}.

	\subsection{Assumptions}
	\begin{assumption}[\bfseries \small Smoothness and strong convexity]
		\label{assump:gt_semi_functions} 
		For every $i\in\{1,\ldots,N\}$, the function $f_i$ is differentiable,
		$L$-smooth, and $\mu$-strongly convex: for all $x,z\in\real^d$,
		\begin{subequations}\label{gt_semi_function_assumptions}
			\begin{align}
				&\|\grad f_i(x)-\grad f_i(z)\|\leq L\|x-z\|,\label{gt_semi_smooth}\\
				&f_i(z)\geq f_i(x)+\langle\grad f_i(x),z-x\rangle
				+\frac{\mu}{2}\|z-x\|^2.\label{gt_semi_sc}
			\end{align}
		\end{subequations}
		Consequently, $f$ is $L$-smooth and $\mu$-strongly convex with
		$0<\mu\leq L$.  Hence, the
		solution $x^\star$ of \eqref{min_problem} is unique and satisfies $
		\sum_{i=1}^{N}\grad f_i(x^\star)=0$.
	\end{assumption}
	
	\begin{assumption}[\bfseries \small Mixing matrix] 
		\label{assump:gt_semi_network}
		The matrix $W=[w_{ij}]\in\real^{N\times N}$ is nonnegative and doubly
		stochastic with
		\begin{equation}\label{gt_semi_lambda}
			\lambda\define
			\|W-\tfrac{1}{N}\one\one\tran\|_2<1,
			\qquad
			\Delta_\lambda\define1-\lambda .
		\end{equation}
	\end{assumption}

	\begin{assumption}[\bfseries\small Gradient noise]
		\label{assump:gt_semi_noise}
		Assume that the initial iterate satisfies
		$\Ex\|\x^0\|^2<\infty$. Let $\cF^k$ denote the sigma-algebra generated by the initial iterate
		$\x^0$ and all samples $\{\bxi^t\}_{t=0}^{k-1}$ drawn before iteration
		$k$. Then $\x^k$ is $\cF^k$-measurable. Conditioned on $\cF^k$, the
		samples $\xi_1^k,\ldots,\xi_N^k$ are mutually independent and, for every
		$k\geq0$ and every agent $i$, satisfy
		\begin{subequations}\label{gt_semi_noise_assumption}
			\begin{align}
				\Ex\!\left[
				\grad F_i(x_i^k;\xi_i^k)-\grad f_i(x_i^k)
				\mid\cF^k
				\right]
				&=0,
				\label{gt_semi_unbiased}
				\\
				\Ex\!\left[
				\|\grad F_i(x_i^k;\xi_i^k)-\grad f_i(x_i^k)\|^2
				\mid\cF^k
				\right]
				&\leq\sigma^2.
				\label{gt_semi_variance}
			\end{align}
		\end{subequations}
	\end{assumption}
	
	\subsection{Change of coordinates and centered dynamics}
	Both analyses in this paper rest on a change of coordinates that centers
	the tracking variable, and on the exact average and disagreement
	dynamics that result.
	Define the change of coordinates
	\begin{equation}\label{gt_semi_y_definition}
		\y^k
		\define
		\g^k-\grad\F(\x^k;\bxi^k).
	\end{equation}
	Moreover, $\y^0=\zero$, and, for $k\geq1$,
	\eqref{gt_atc_semi_g} gives
	$\y^k=\W\g^{k-1}-\grad\F(\x^{k-1};\bxi^{k-1})$.
	Thus $\y^k$, and hence $\widetilde{\y}^k$ and
	$\widehat{\e}^k$ below, are $\cF^k$-measurable for every $k\geq0$,
	whereas $\g^k$ is generally not.

	Substituting
	$\g^k=\y^k+\grad\F(\x^k;\bxi^k)$
	into \eqref{gt_atc_semi} gives the equivalent recursion
	\begin{subequations}\label{gt_atc_semi2}
		\begin{align}
			\x^{k+1}
			&=
			\W\left(
			\x^k-\alpha\grad\F(\x^k;\bxi^k)-\alpha\y^k
			\right),
			\label{gt_atc_semi2_x}
			\\
			\y^{k+1}
			&=
			\W\y^k
			-(\I-\W)\grad\F(\x^k;\bxi^k).
			\label{gt_atc_semi2_y}
		\end{align}
	\end{subequations}	
	To separate the average and disagreement dynamics, isolate the
	stochastic-gradient noise, and center the recursion at the optimizer, we
	introduce
	\begin{subequations}\label{gt_semi_definitions}
		\begin{alignat}{2}
			\mathbf\Pi
			&\define
			\frac{1}{N}\one\one\tran\otimes I_d,
			&\qquad
			\widetilde{\W}
			&\define
			\W-\mathbf\Pi,
			\label{gt_semi_projectors}
			\\
			\bar x^k
			&\define
			\frac{1}{N}\sum_{i=1}^{N}x_i^k,
			&\qquad
			\widehat{\x}^k
			&\define
			\x^k-\one\otimes\bar x^k
			=
			(\I-\mathbf\Pi)\x^k,
			\label{gt_semi_x_average_disagreement}
			\\
			\v^k
			&\define
			\grad\F(\x^k;\bxi^k)-\grad\f(\x^k),
			&\qquad
			\bar v^k
			&\define
			\frac{1}{N}\sum_{i=1}^{N}v_i^k,
			\label{gt_semi_noise_definition}
			\\
			\x^\star
			&\define
			\one\otimes x^\star,
			&\qquad
			\widetilde{\y}^k
			&\define
			\y^k+\grad\f(\x^\star),
			\label{gt_semi_xstar_stack}
			\\
			\overline{\grad\f}(\x^k)
			&\define
			\frac{1}{N}\sum_{i=1}^{N}\grad f_i(x_i^k),
			&\qquad
			\mathbf u^k
			&\define
			\grad\f(\x^k)-\grad\f(\x^\star).
			\label{gt_semi_u_definition}
		\end{alignat}
	\end{subequations}
	Here $\mathbf\Pi$ is the consensus projector and $\widetilde{\W}$ the
	disagreement mixing matrix; $\bar x^k$ and $\widehat{\x}^k$ are the
	average iterate and the disagreement; $\v^k$ is the stochastic-gradient
	noise, with network average $\bar v^k$; and $\widetilde{\y}^k$ and $\mathbf u^k$ center the tracking and gradient quantities at the optimizer. Double stochasticity of $W$ gives
	\begin{equation*}
		\mathbf\Pi\W
		=
		\W\mathbf\Pi
		=
		\mathbf\Pi,
		\qquad
		(\I-\mathbf\Pi)\W
		=
		\widetilde{\W},
		\qquad
		\widetilde{\W}\mathbf\Pi
		=
		\zero .
	\end{equation*}
	
	\begin{lemma}[\bfseries\small Average and disagreement dynamics]
		\label{lemma:gt_semi_error_dynamics}
		Let Assumption~\ref{assump:gt_semi_functions} hold, let $W$ be doubly
		stochastic, and let $\g^0=\grad\F(\x^0;\bxi^0)$. Then,
		for all $k\geq0$, $\bar y^k\define(1/N)\sum_{i=1}^Ny_i^k=0$, $\mathbf\Pi\widetilde{\y}^k=\zero$, and
		the average iterate satisfies
		\begin{equation}\label{gt_semi_centroid}
			\bar x^{k+1}
			=
			\bar x^k
			-\alpha\overline{\grad\f}(\x^k)
			-\alpha\bar v^k.
		\end{equation}
		Moreover, the disagreement and centered tracking variable satisfy
		\begin{subequations}\label{gt_semi_centered_recursions}
			\begin{align}
				\widehat{\x}^{k+1}
				&=
				\widetilde{\W}
				\left(
				\widehat{\x}^k
				-\alpha\widetilde{\y}^k
				-\alpha\mathbf u^k
				-\alpha\v^k
				\right),
				\label{gt_semi_centered_x}
				\\
				\widetilde{\y}^{k+1}
				&=
				\widetilde{\W}\widetilde{\y}^k
				-(\I-\W)(\mathbf u^k+\v^k).
				\label{gt_semi_centered_y}
			\end{align}
		\end{subequations}
	\end{lemma}
	
	\begin{proof}
		Applying $\mathbf\Pi$ to \eqref{gt_atc_semi2_y} and using
		$\mathbf\Pi\W=\mathbf\Pi$ gives
		$\mathbf\Pi\y^{k+1}=\mathbf\Pi\y^k-(\mathbf\Pi-\mathbf\Pi\W)\grad\F(\x^k;\bxi^k)
		=\mathbf\Pi\y^k$, so $\mathbf\Pi\y^k=\mathbf\Pi\y^0=\zero$, which is
		$\bar y^k=0$. Since $\sum_{i=1}^{N}\grad f_i(x^\star)=0$, also
		$\mathbf\Pi\grad\f(\x^\star)=\zero$, and hence
		$\mathbf\Pi\widetilde{\y}^k=\zero$. Applying $\mathbf\Pi$ to
		\eqref{gt_atc_semi2_x} and using $\mathbf\Pi\W=\mathbf\Pi$,
		$\mathbf\Pi\y^k=\zero$, and
		$\grad\F(\x^k;\bxi^k)=\grad\f(\x^k)+\v^k$ gives
		\eqref{gt_semi_centroid}.
		
		For the disagreement, apply $\I-\mathbf\Pi$ to
		\eqref{gt_atc_semi2_x} and use $(\I-\mathbf\Pi)\W=\widetilde{\W}$ and
		$\widetilde{\W}\mathbf\Pi=\zero$, so that
		$\widetilde{\W}\x^k=\widetilde{\W}\widehat{\x}^k$:
		\[
		\widehat{\x}^{k+1}
		=
		\widetilde{\W}
		\left(
		\widehat{\x}^k-\alpha\grad\F(\x^k;\bxi^k)-\alpha\y^k
		\right).
		\]
		Substituting
		$\grad\F(\x^k;\bxi^k)=\mathbf u^k+\v^k+\grad\f(\x^\star)$ and
		$\y^k=\widetilde{\y}^k-\grad\f(\x^\star)$, the two terms in
		$\grad\f(\x^\star)$ cancel, which gives
		\eqref{gt_semi_centered_x}. Finally, adding $\grad\f(\x^\star)$ to both sides of
		\eqref{gt_atc_semi2_y} and substituting
		$\y^k=\widetilde{\y}^k-\grad\f(\x^\star)$ and
		$
		\grad\F(\x^k;\bxi^k)
		=
		\mathbf u^k+\v^k+\grad\f(\x^\star)
		$
		gives
		\[
		\widetilde{\y}^{k+1}
		=
		\W\widetilde{\y}^k
		-(\I-\W)(\mathbf u^k+\v^k),
		\]
		where the terms involving $\grad\f(\x^\star)$ cancel.
		Since $\mathbf\Pi\widetilde{\y}^k=\zero$, we have
		$\W\widetilde{\y}^k=\widetilde{\W}\widetilde{\y}^k$, which proves
		\eqref{gt_semi_centered_y}.
	\end{proof}
	
	\section{Baseline convergence analysis}
	\label{sec:gt_baseline_result}
	For completeness, we first establish a baseline convergence guarantee for
	the stochastic gradient-tracking recursion \eqref{gt_atc_semi}, starting
	from the centered dynamics in
	Lemma~\ref{lemma:gt_semi_error_dynamics}. Bounds of this type are known
	\cite{pu2021distributed,xin2021improved}.

	\subsection{Baseline convergence bound}
	The two disagreement variables in
	\eqref{gt_semi_centered_recursions} can be combined into a single scaled
	error. Throughout this subsection, we assume $0<\lambda<1$. Define
	\begin{align} \label{gt_semi_eta_delta_joint}
		\eta&\define\frac{\Delta_\lambda}{2\lambda},
		\quad
		\delta\define\frac{\Delta_\lambda}{2},\quad
		\widehat{\e}^k\define
		\begin{bmatrix}
			\widehat{\x}^k
			\\[1mm]
			\tfrac{\alpha}{\eta}\widetilde{\y}^k
		\end{bmatrix}.
	\end{align}
	The case $\lambda=0$ corresponds to exact averaging and does not need this scaled error. Then \eqref{gt_semi_centered_recursions} gives
	\begin{equation}\label{gt_semi_joint_exact}
		\widehat{\e}^{k+1}
		=
		\C\widehat{\e}^k
		-\alpha\B\mathbf u^k
		-\alpha\B\v^k,
	\end{equation}
	where
	\begin{equation}\label{gt_semi_C_B}
		\C
		\define
		\begin{bmatrix}
			\widetilde{\W} & -\eta\widetilde{\W}
			\\
			\zero & \widetilde{\W}
		\end{bmatrix},
		\qquad
		\B
		\define
		\begin{bmatrix}
			\widetilde{\W}
			\\[1mm]
			\tfrac{1}{\eta}(\I-\W)
		\end{bmatrix}.
	\end{equation}
	Thus, the network part of the GT dynamics is represented by the stable
	linear system \eqref{gt_semi_joint_exact}, driven by the gradient
	perturbation $\mathbf u^k$ and the stochastic-gradient noise $\v^k$.
	
	We now state the convergence result. Define the average optimization error
	and the joint network error by
	\begin{equation}\label{eq:baseline_errors_main}
		\bar e_k
		\define
		\Ex\|\bar x^k-x^\star\|^2,
		\qquad
		b_k
		\define
		\frac{1}{N}\Ex\|\widehat{\e}^k\|^2
		=
		\frac{1}{N}\Ex[
		\|\widehat{\x}^k\|^2
		+\tfrac{\alpha^2}{\eta^2}
		\|\widetilde{\y}^k\|^2
		].
	\end{equation}
	We also define
	\begin{equation}\label{eq:baseline_Q_main}
		Q_\lambda
		\define
		\frac{17\lambda^2}{\Delta_\lambda^2},
		\qquad
		\varsigma_\star^2
		\define
		\frac{1}{N}
		\|\grad\f(\x^\star)\|^2.
	\end{equation}
	
	\begin{theorem}[\bfseries\small Baseline GT convergence]
		\label{thm:gt_baseline_convergence}
		Suppose Assumptions~\ref{assump:gt_semi_functions}--\ref{assump:gt_semi_noise}
		hold, $0<\lambda<1$, the primal variable is initialized consensually,
		$\x^0=\one\otimes x^0$, and $\g^0=\grad\F(\x^0;\bxi^0)$. Suppose
		\begin{equation}\label{eq:baseline_stepsize_main}
			0<\alpha\leq\bar\alpha, \quad \bar\alpha
			\define
			\min\left\{
			\frac{1}{4L},
			\frac{\delta}{L\sqrt{24Q_\lambda}},
			\frac{\delta}{2\mu}
			\right\}.
		\end{equation}
		Then, with
		\begin{equation}\label{eq:baseline_lyapunov_main}
			\mathcal V_k
			\define
			\bar e_k+\frac{3\alpha L}{\delta}b_k,
		\end{equation}
		for every $k\geq0$,
		\begin{equation}\label{eq:baseline_linear_bound_main}
			\mathcal V_k
			\leq
			\left(1-\frac{\mu\alpha}{2}\right)^k\mathcal V_0
			+\frac{2\alpha\sigma^2}{\mu N}
			+\frac{6\alpha^2Q_\lambda L\sigma^2}{\mu\delta},
		\end{equation}
		where $
		\mathcal V_0
		=
		\bar e_0
		+\frac{3\alpha^3L}{\eta^2\delta}\varsigma_\star^2$.
		If $\sigma^2>0$, then for every $K\geq1$ there is an admissible stepsize
		$\alpha=\min\{\bar\alpha,\,2\ln T_K/(\mu K)\}$, with $T_K$ as in
		Lemma~\ref{lem:log_stepsize}, for which
		\begin{equation}\label{eq:baseline_tuned_rate_main}
			\mathcal V_K
			=
			\widetilde\cO\!\left(
			\frac{\sigma^2}{\mu^2NK}
			+
			\frac{L\sigma^2}
			{\mu^3\Delta_\lambda^3K^2}
			\right)
			+
			\left(
			\bar e_0
			+
			\frac{24L\lambda^2\bar\alpha^3}
			{\Delta_\lambda^3}\varsigma_\star^2
			\right)
			\exp\!\left(
			-\frac{\mu\bar\alpha K}{2}
			\right).
		\end{equation}
	\end{theorem}
	\begin{proof}
		See Appendix~\ref{sec:gt_two_recursion_appendix}.
	\end{proof}
	Regarding the network, Theorem~\ref{thm:gt_baseline_convergence}
	requires the norm condition in Assumption~\ref{assump:gt_semi_network}
	but not diagonalizability, unlike the refined analysis of
	Section~\ref{sec:gt_nonsymmetric_refined}.
	For $\sigma^2=0$, taking $\alpha=\bar\alpha$ gives the corresponding
	geometric bound directly. Consensual initialization simplifies
	the displayed initial-error term only. For a constant stepsize,
	\eqref{eq:baseline_linear_bound_main} gives the steady-state bound
	$
	\frac{2\alpha\sigma^2}{\mu N}
	+
	\frac{204\alpha^2\lambda^2L\sigma^2}
	{\mu\Delta_\lambda^3}.
	$
	Thus, the leading term is network independent and enjoys the usual
	$1/N$ variance reduction, while the network appears in the higher-order
	term. With horizon-dependent tuning, these terms become, up to
	logarithmic factors,
	$\sigma^2/(\mu^2NK)$ and
	$L\sigma^2/(\mu^3\Delta_\lambda^3K^2)$.

	\subsection{A two-agent cubic lower bound}
	\label{subsec:lb_two_agent_disagreement}
	The following two-agent construction shows that the coefficient of
	$\alpha^2\sigma^2$ in the steady-state mean-square disagreement can
	scale as $\Delta_\lambda^{-3}$ when the nontrivial eigenvalue
	approaches $-1$.

	\begin{lemma}[\bfseries\small A two-agent mixing matrix with a negative mode]
		\label{lem:negative_mode_matrix}
		For any $\lambda\in[1/2,1)$, let
		\begin{equation}\label{eq:negative_mode_matrix}
			W_\lambda
			\define
			\frac{1}{2}
			\begin{bmatrix}
				\Delta_\lambda & 1+\lambda\\
				1+\lambda & \Delta_\lambda
			\end{bmatrix},
			\qquad
			\Pi
			\define
			\tfrac12\one\one\tran .
		\end{equation}
		Then $W_\lambda$ is symmetric, entrywise positive, and doubly stochastic.
		Its eigenvalues are $1$ and $-\lambda$, and
		\begin{equation}\label{eq:negative_mode_norm}
			W_\lambda-\Pi
			=
			-\lambda vv\tran,
			\qquad
			\|W_\lambda-\Pi\|_2
			=
			\lambda,
			\qquad
			v
			\define
			\tfrac{1}{\sqrt2}
			\begin{bmatrix}
				1\\-1
			\end{bmatrix}.
		\end{equation}
	\end{lemma}
	
	\begin{proof}
		Every row and column of $W_\lambda$ sums to one, and all its entries are
		strictly positive for $\lambda\in[1/2,1)$. Moreover, $W_\lambda\one=\one$ and
		$W_\lambda v=-\lambda v$.
		Since $\{\one/\sqrt2,v\}$ is an orthonormal basis of $\real^2$, the spectral
		decomposition of the symmetric matrix $W_\lambda$ is
		\[
		W_\lambda
		=
		1\cdot
		\frac{\one\one\tran}{2}
		+
		(-\lambda)vv\tran
		=
		\Pi-\lambda vv\tran.
		\]
		Therefore, $W_\lambda-\Pi
		=
		-\lambda vv\tran$.
		Since $v$ is a unit vector, $vv\tran$ is the orthogonal projector onto the
		one-dimensional disagreement subspace and satisfies
		$\|vv\tran\|_2=1$. Hence
		$
		\|W_\lambda-\Pi\|_2
		=
		\lambda\|vv\tran\|_2
		=
		\lambda$,
		which proves \eqref{eq:negative_mode_norm}.
	\end{proof}

	To derive the lower bound, we consider the scalar quadratic instance
	\begin{equation}\label{eq:lb_objective}
		f_i(x)
		=
		\frac{\mu}{2}x^2,
		\qquad
		i\in\{1,2\},
	\end{equation}
	for which $x^\star=0$, with stochastic gradients $\grad F_i(x;\xi_i^k)=
	\mu x+\xi_i^k$ and
	\begin{align}  \label{eq:lb_noise}
		\Ex[\xi_i^k]&=0,
		\quad
		\Ex[|\xi_i^k|^2]=\sigma^2,
	\end{align}
	where the noises $\xi_i^k$ are independent across agents and iterations.
	\begin{theorem}[\bfseries\small Exact steady state and cubic lower bound]
		\label{thm:cubic_network_lower_bound}
		Let $N=2$, $d=1$, $W=W_\lambda$ as in
		\eqref{eq:negative_mode_matrix}, and consider
		\eqref{eq:lb_objective}--\eqref{eq:lb_noise}.
		Run the GT recursion \eqref{gt_atc_semi2} from any deterministic $\x^0$
		with $\y^0=\zero$.
		If
		\begin{equation}\label{eq:lb_stability_stepsize}
			\alpha\mu
			<
			\frac{\Delta_\lambda(1+\lambda)}{\lambda},
		\end{equation}
		then the disagreement dynamics are stable and
		\begin{equation}\label{eq:lb_exact_steady_state}
			\lim_{k\rightarrow\infty}
			\frac{1}{N}\Ex\|\widehat{\x}^k\|^2
			=
			\frac{
				\alpha^2\lambda^2\sigma^2
			}{
				\bigl(1-\lambda^2-\lambda\alpha\mu\bigr)
				\bigl(\Delta_\lambda^2+2\lambda\alpha\mu\bigr)
			}.
		\end{equation}	
		If, in addition,
		\begin{equation}\label{eq:lb_stepsize}
			\alpha\mu
			\leq
			\frac{\Delta_\lambda^2}{16},
		\end{equation}
		then
		\begin{equation}\label{eq:lb_cubic_bound}
			\frac{\alpha^2\sigma^2}
			{9\Delta_\lambda^3}
			\leq
			\lim_{k\rightarrow\infty}
			\frac{1}{N}\Ex\|\widehat{\x}^k\|^2
			\leq
			\frac{2\alpha^2\sigma^2}
			{\Delta_\lambda^3}.
		\end{equation}
	\end{theorem}
	\begin{proof}
		See Appendix \ref{app:cubic_network_lower_bound}.
	\end{proof}
	Condition \eqref{eq:lb_stepsize} is not restrictive: on this instance
	$L=\mu$, and every stepsize admissible in
	Theorem~\ref{thm:gt_baseline_convergence} satisfies
	$\alpha\mu\leq\Delta_\lambda^2/(4\sqrt{102}\,\lambda)\leq\Delta_\lambda^2/16$
	for $\lambda\geq1/2$.
	For $\sigma^2>0$ and every fixed $\lambda\in[1/2,1)$,
	\eqref{eq:lb_exact_steady_state} gives
	\begin{equation}\label{eq:lb_small_stepsize_limit}
		\lim_{\alpha\downarrow0}
		\frac{1}{\alpha^2\sigma^2}
		\left(
		\lim_{k\rightarrow\infty}
		\frac1N\Ex\|\widehat{\x}^k\|^2
		\right)
		=
		\frac{\lambda^2}
		{(1-\lambda^2)\Delta_\lambda^2}
		=
		\frac{\lambda^2}
		{(1+\lambda)\Delta_\lambda^3}.
	\end{equation}
	Hence the coefficient of the $\alpha^2\sigma^2$ steady-state
	disagreement is
	$
	\frac{\lambda^2}{1+\lambda}\,\Delta_\lambda^{-3}.
	$
	Since $\lambda\in[1/2,1)$, the prefactor $\lambda^2/(1+\lambda)$ is
	bounded above and below by positive absolute constants. Therefore, the
	inverse-cubic dependence on $\Delta_\lambda$ cannot be improved
	uniformly over this family in the small-stepsize regime. This behavior
	is caused by the negative eigenvalue $-\lambda$, which approaches $-1$
	as $\lambda\uparrow1$.
	
	Moreover, since $\x^k-\x^\star=\widehat{\x}^k+\one\otimes(\bar x^k-x^\star)$
	with orthogonal terms, under \eqref{eq:lb_stepsize},
	\[
	\liminf_{k\rightarrow\infty}
	\frac1N\Ex\|\x^k-\x^\star\|^2
	=
	\liminf_{k\rightarrow\infty}
	\frac1N\Ex\bigl(\|\widehat{\x}^k\|^2+N\|\bar x^k-x^\star\|^2\bigr)
	\geq
	\frac{\alpha^2\sigma^2}{9\Delta_\lambda^3}.
	\]
	Hence, the cubic dependence on $\Delta_\lambda^{-1}$ cannot, in general, be
	improved in the $\alpha^2\sigma^2$ term for either the disagreement or the
	total primal error.


	\section{Refined analysis for diagonalizable mixing matrices}
	\label{sec:gt_nonsymmetric_refined}
	The baseline analysis of Section~\ref{sec:gt_baseline_result} uses a
	single contraction factor, yielding the network dependence
	$\cO(\Delta_\lambda^{-3})$. We sharpen this bound by
	unrolling the disagreement recursion and controlling the
	stochastic-gradient noise through the squared propagation coefficients.
	
	\begin{assumption}[\bfseries\small Diagonalizable mixing matrix]
		\label{assump:gt_nonsym_matrix}
		The matrix $W\in\real^{N\times N}$ is nonnegative and doubly stochastic.
		With $\Pi\define\frac1N\one\one\tran$ and
		$\widetilde W\define W-\Pi$, assume that $\widetilde W$ is
		diagonalizable over $\mathbb C$ as
		\begin{align}\label{gt_nonsym_diagonalization}
			\widetilde W
			=
			S\diag\{0,\lambda_2,\ldots,\lambda_N\}S^{-1},
			\qquad
			|\lambda_i|<1,
		\end{align}
		for some nonsingular $S\in\mathbb C^{N\times N}$. Define
		\begin{subequations}\label{gt_nonsym_spectrum_interval}
			\begin{alignat}{2}
				\lambda
				&\define
				\max_{2\leq i\leq N}|\lambda_i|,
				&\qquad
				\Delta_\lambda
				&\define
				1-\lambda,
				\\
				\beta_W
				&\define
				\max_{2\leq i\leq N}
				\frac{|1-\lambda_i|}{1-|\lambda_i|},
				&\qquad
				\kappa_W
				&\define
				\|S\|_2\|S^{-1}\|_2 .
			\end{alignat}
		\end{subequations}
		The value of $\kappa_W$ refers to a fixed choice of $S$ with unit-norm
		columns; for normal $W$, including the symmetric case, we take $S$
		unitary, so that $\kappa_W=1$.
	\end{assumption}
	
	Since $W$ is doubly stochastic,
	$\widetilde W\one=\zero$ and $\one\tran\widetilde W=\zero\tran$, so the
	displayed zero eigenvalue in \eqref{gt_nonsym_diagonalization} corresponds to
	the consensus mode. Let $\widetilde{\W}\define\widetilde W\otimes I_d$.
	Diagonalizability gives, for every $m\geq0$,
	\begin{equation}\label{gt_nonsym_power_decay}
		\|\widetilde{\W}^m\|_2
		\leq
		\kappa_W
		\max_{2\leq i\leq N}|\lambda_i|^m
		=
		\kappa_W(1-\Delta_\lambda)^m.
	\end{equation}
	Since $W$ is nonnegative and doubly stochastic, $\|\widetilde{\W}\|_2\leq1$,
	$\|\I-\W\|_2\leq2$.
	\begin{remark}[\bfseries\small Spectral gap in the nonsymmetric case]
		\rm
		For nonsymmetric $W$, the spectral radius
		$\lambda=\max_{i\geq2}|\lambda_i|$ may be strictly smaller than
		$\|\widetilde{\W}\|_2$. Thus, unlike the baseline analysis,
		$\Delta_\lambda$ here is defined using the spectral radius rather than
		the one-step spectral norm. The two definitions coincide when $W$ is
		normal/symmetric. 
	\end{remark}
	
	\begin{remark}[\bfseries\small Interpretation of $\beta_W$]
		\rm
		From $1-|\lambda_i|\leq|1-\lambda_i|\leq2$,
		\[
		1\leq\beta_W\leq\frac{2}{\Delta_\lambda}.
		\]
		If all nontrivial eigenvalues are real and nonnegative, then
		$\beta_W=1$. For symmetric $W$, $\kappa_W=1$; moreover, if
		$\lambda_N(W)<0$ and $\nu\define1+\lambda_N(W)$, then
		\[
		\beta_W
		=
		\frac{1-\lambda_N(W)}{1+\lambda_N(W)}
		=
		\frac{2-\nu}{\nu}.
		\]
		When an eigenvalue $z$ approaches the unit circle away from $1$, its
		ratio $|1-z|/(1-|z|)$ grows as $(1-|z|)^{-1}$.
		Remark~\ref{rem:gt_nonsym_modes} shows why this gives cubic noise
		amplification and why growth can be milder near $1$.
		The factor $\kappa_W$ measures eigenvector nonorthogonality.
	\end{remark}
	
	\subsection{Exact disagreement propagation}
	
	\label{sec:gt_nonsym_filter}
	
	Define the projected perturbation
	\begin{equation}\label{gt_nonsym_h_definition}
		\widehat{\mathbf r}^k
		\define
		(\I-\mathbf\Pi)(\mathbf u^k+\v^k),
	\end{equation}
	where $\mathbf u^k$ is defined in \eqref{gt_semi_u_definition}. Using
	$\widetilde{\W}(\I-\mathbf\Pi)=\widetilde{\W}$ and
	$(\I-\W)(\I-\mathbf\Pi)=\I-\W$, the centered recursions
	\eqref{gt_semi_centered_recursions} become
	\begin{subequations}\label{gt_nonsym_disagreement_system}
		\begin{align}
			\widehat{\x}^{k+1}
			&=
			\widetilde{\W}\widehat{\x}^k
			-\alpha\widetilde{\W}\widetilde{\y}^k
			-\alpha\widetilde{\W}\widehat{\mathbf r}^k,
			\\
			\widetilde{\y}^{k+1}
			&=
			\widetilde{\W}\widetilde{\y}^k
			-(\I-\W)\widehat{\mathbf r}^k .
		\end{align}
	\end{subequations}
	
	\begin{subequations}
		\begin{lemma}[\bfseries\small Exact disagreement representation]
			\label{lemma:gt_nonsym_convolution}
			For every $k\geq0$, the disagreement in
			\eqref{gt_nonsym_disagreement_system} satisfies
			\begin{equation}\label{gt_nonsym_convolution_identity}
				\widehat{\x}^{k}
				=
				\underbrace{
					\widetilde{\W}^{k}\widehat{\x}^0
					-\alpha k\,\widetilde{\W}^{k}\widetilde{\y}^0
				}_{\define\ \T^k}
				+
				\alpha\sum_{t=0}^{k-1}
				\mathbf M_{k-1-t}\widehat{\mathbf r}^{t},
			\end{equation}
			where the propagation matrices are
			\begin{equation}\label{gt_nonsym_B_definition}
				\mathbf M_\ell
				\define
				\ell\,\widetilde{\W}^{\ell}
				-(\ell+1)\widetilde{\W}^{\ell+1},
				\qquad \ell\geq0 .
			\end{equation}
			
		\end{lemma}
		
		\begin{proof}
			The second recursion in \eqref{gt_nonsym_disagreement_system} gives
			\begin{equation}\label{gt_nonsym_y_unrolled}
				\widetilde{\y}^j
				=
				\widetilde{\W}^j\widetilde{\y}^0
				-
				\sum_{t=0}^{j-1}
				\widetilde{\W}^{j-1-t}
				(\I-\W)\widehat{\mathbf r}^t .
			\end{equation}
			For $k=0$, \eqref{gt_nonsym_convolution_identity} holds because the
			sum is empty and $\T^0=\widehat{\x}^0$. For $k\geq1$, unrolling the
			first recursion gives
			\begin{equation}\label{gt_nonsym_x_unrolled}
				\widehat{\x}^{k}
				=
				\widetilde{\W}^{k}\widehat{\x}^0
				-\alpha\sum_{j=0}^{k-1}
				\widetilde{\W}^{k-j}\widetilde{\y}^j
				-\alpha\sum_{t=0}^{k-1}
				\widetilde{\W}^{k-t}\widehat{\mathbf r}^t .
			\end{equation}
			Substituting \eqref{gt_nonsym_y_unrolled} into
			\eqref{gt_nonsym_x_unrolled} gives
			\begin{align}
				\widehat{\x}^{k}
				={}&
				\widetilde{\W}^{k}\widehat{\x}^0
				-\alpha\sum_{j=0}^{k-1}
				\widetilde{\W}^{k-j}\widetilde{\W}^j
				\widetilde{\y}^0
				\nonumber\\
				&\quad
				+\alpha\sum_{j=0}^{k-1}
				\sum_{t=0}^{j-1}
				\widetilde{\W}^{k-j}
				\widetilde{\W}^{j-1-t}
				(\I-\W)\widehat{\mathbf r}^t
				-\alpha\sum_{t=0}^{k-1}
				\widetilde{\W}^{k-t}\widehat{\mathbf r}^t .
				\label{gt_nonsym_x_substituted}
			\end{align}
			Since $
			\sum_{j=0}^{k-1}
			\widetilde{\W}^{k-j}\widetilde{\W}^{j}
			=
			k\widetilde{\W}^{k}$,
			the terms involving $\widetilde{\y}^0$ in
			\eqref{gt_nonsym_x_substituted} reduce to
			\begin{equation}\label{gt_nonsym_y0_term}
				-\alpha\sum_{j=0}^{k-1}
				\widetilde{\W}^{k-j}\widetilde{\W}^{j}
				\widetilde{\y}^0
				=
				-\alpha k\,\widetilde{\W}^{k}\widetilde{\y}^0 .
			\end{equation}
			For the double sum in \eqref{gt_nonsym_x_substituted}, exchanging the order
			of summation gives
			\begin{align}
				\sum_{j=0}^{k-1}\sum_{t=0}^{j-1}
				\widetilde{\W}^{k-j}
				\widetilde{\W}^{j-1-t}
				(\I-\W)\widehat{\mathbf r}^t
				&=
				\sum_{t=0}^{k-1}
				\sum_{j=t+1}^{k-1}
				\widetilde{\W}^{k-1-t}
				(\I-\W)\widehat{\mathbf r}^t
				\nonumber\\
				&=
				\sum_{t=0}^{k-1}
				(k-1-t)\widetilde{\W}^{k-1-t}
				(\I-\W)\widehat{\mathbf r}^t .
				\label{gt_nonsym_double_sum}
			\end{align}
			Indeed, for a fixed perturbation index $t\in\{0,\ldots,k-1\}$, the inner sum runs over
			$j=t+1,\ldots,k-1$ (which is empty when $t=k-1$), giving $k-1-t$ identical terms.
			Setting
			$\ell\define k-1-t$,
			the contribution of $\widehat{\mathbf r}^t$ obtained from
			\eqref{gt_nonsym_double_sum} and the last term of
			\eqref{gt_nonsym_x_substituted} is
			\begin{equation}\label{gt_nonsym_rt_coefficient}
				\alpha
				\left[
				\ell\widetilde{\W}^{\ell}(\I-\W)
				-\widetilde{\W}^{\ell+1}
				\right]
				\widehat{\mathbf r}^t .
			\end{equation}
			By definition, $\widehat{\mathbf r}^t$ lies in the disagreement subspace, so
			\begin{equation}\label{gt_nonsym_r_projected}
				(\I-\mathbf\Pi)\widehat{\mathbf r}^t
				=
				\widehat{\mathbf r}^t,
				\qquad
				\mathbf\Pi\widehat{\mathbf r}^t
				=
				\zero .
			\end{equation}
			Since $\W=\mathbf\Pi+\widetilde{\W}$, \eqref{gt_nonsym_r_projected} implies
			$
			\W\widehat{\mathbf r}^t
			=
			\widetilde{\W}\widehat{\mathbf r}^t$,
			and therefore
			\begin{equation}\label{gt_nonsym_IW_on_r}
				(\I-\W)\widehat{\mathbf r}^t
				=
				(\I-\widetilde{\W})\widehat{\mathbf r}^t .
			\end{equation}
			Left-multiplying \eqref{gt_nonsym_IW_on_r} by
			$\widetilde{\W}^{\ell}$ gives
			\begin{equation}\label{gt_nonsym_power_IW}
				\widetilde{\W}^{\ell}(\I-\W)\widehat{\mathbf r}^t
				=
				(
				\widetilde{\W}^{\ell}
				-\widetilde{\W}^{\ell+1}
				)\widehat{\mathbf r}^t .
			\end{equation}
			Substituting \eqref{gt_nonsym_power_IW} into
			\eqref{gt_nonsym_rt_coefficient} gives
			$\alpha\mathbf M_\ell\widehat{\mathbf r}^t$,
			with $\mathbf M_\ell$ as in \eqref{gt_nonsym_B_definition}.
			Together with \eqref{gt_nonsym_x_substituted} and
			\eqref{gt_nonsym_y0_term}, this proves
			\eqref{gt_nonsym_convolution_identity}.
		\end{proof}
	\end{subequations}
	
	Equation~\eqref{gt_nonsym_convolution_identity} is a discrete convolution:
	each past perturbation $\widehat{\mathbf r}^t$ is propagated to the current
	disagreement through the matrix $\mathbf M_{k-1-t}$. Thus, the sequence
	$\{\mathbf M_\ell\}_{\ell\geq0}$ characterizes how past perturbations affect
	the disagreement. Since $\y^0=\zero$, the homogeneous term $\T^k$ contains
	$\widetilde{\W}^{k}\widehat{\x}^0$ and
	$-\alpha k\widetilde{\W}^{k}\grad\f(\x^\star)$.
	The first term vanishes for consensual initialization.
	The initial noise $\v^0$ appears in the $t=0$ convolution term,
	not in $\T^k$.
	
		For $\omega\in(0,1]$, define
	\begin{equation}\label{gt_nonsym_network_gains}
		\begin{aligned}
			\mathcal G_1(\omega)
			&\define\sum_{\ell\geq0}\omega^{-\ell}\|\mathbf M_\ell\|_2,
			\qquad
			\mathcal G_2\define\sum_{\ell\geq0}\|\mathbf M_\ell\|_2^2,
			\\
			\mathcal G_{2,\mathrm{av}}
			&\define\frac1{Nd}\sum_{\ell\geq0}\|\mathbf M_\ell\|_F^2
			\leq\mathcal G_2.
		\end{aligned}
	\end{equation}
	The gain $\mathcal G_1(\omega)$ controls optimization perturbations
	and the admissible stepsize, while $\mathcal G_{2,\mathrm{av}}$
	controls the accumulated stochastic variance.
	
	\subsection{Spectral bounds for the propagation gains}
	
	\label{sec:gt_nonsym_spectral}
		To make the gains in \eqref{gt_nonsym_network_gains} explicit for the
		analysis, we next bound $\mathbf M_\ell$ through the eigenvalues of $W$. Define
	\begin{align} \label{ml_poly}
		m_\ell(z)
		\define
		\ell z^\ell-(\ell+1)z^{\ell+1}
		=
		z^\ell(\ell(1-z)-z)
	\end{align}
	so that $	\mathbf M_\ell=m_\ell(\widetilde{\W})$.
	
	\begin{lemma}[\bfseries\small Exact squared propagation of one mode]
		\label{lem:gt_nonsym_modal_energy}
		For every $z\in\mathbb C$ with $|z|<1$,
		\begin{equation}\label{gt_nonsym_modal_energy}
			\sum_{\ell=0}^{\infty}|m_\ell(z)|^2
			=
			\frac{2|z|^2|1-z|^2}
			{(1-|z|^2)^3}=	\frac{2|z|^2\beta_z^2}
			{d_z(1+|z|)^3},
		\end{equation}
		where $
		d_z\define1-|z|$ and
		$\beta_z\define\frac{|1-z|}{1-|z|}$.
	\end{lemma}
	\begin{proof}
		See Appendix \ref{app:gt_nonsym_modal_energy}.
	\end{proof}

	\begin{remark}[\bfseries\small Dependence on the location of an eigenvalue]
		\label{rem:gt_nonsym_modes}
		\rm
		Equation~\eqref{gt_nonsym_modal_energy} depends on the location of
		$z$, not only on $|z|$. Writing $z=re^{\imath\theta}$ with $r\in[0,1)$
		and $\theta\in[-\pi,\pi]$, we have
		$|1-z|^2=(1-r)^2+4r\sin^2(\theta/2)$, so that
		\[
		\sum_{\ell\geq0}|m_\ell(z)|^2
		=\frac{2r^2}{(1+r)^3}\left[\frac1{1-r}
		+\frac{4r\sin^2(\theta/2)}{(1-r)^3}\right].
		\]
		The first term is present for every mode and gives
		$\Theta((1-r)^{-1})$ for a positive real mode ($\theta=0$); the second
		gives $\Theta((1-r)^{-3})$ for a negative real mode ($\theta=\pi$), and
		more generally whenever $|1-z|$ is bounded away from zero, that is,
		when $\beta_z=\Theta((1-r)^{-1})$. In between, if $|\theta|\lesssim1-r$
		the energy stays $\Theta((1-r)^{-1})$, and if $1-r\ll|\theta|\ll1$ it is
		$\Theta(\theta^2(1-r)^{-3})$; all $\Theta$ comparisons are as
		$r\uparrow1$.
	\end{remark}
	
	The previous lemma characterizes the propagation associated with a single
	eigenvalue. We next derive bounds that hold uniformly over all nontrivial
	eigenvalues.
	
	\begin{lemma}[\bfseries\small Uniform propagation bounds]
		\label{lem:gt_nonsym_envelope}
		Let Assumption~\ref{assump:gt_nonsym_matrix} hold and define
		\begin{equation}\label{gt_nonsym_omega}
			\omega
			\define
			1-\frac{\Delta_\lambda}{4}.
		\end{equation}
		Then
		\begin{equation}\label{gt_nonsym_envelope_sums}
			\sum_{\ell=0}^{\infty}
			\omega^{-\ell}
			\max_{i\geq2}|m_\ell(\lambda_i)|
			\leq
			\frac{11\beta_W}{\Delta_\lambda},
			\qquad
			\sum_{\ell=0}^{\infty}
			\max_{i\geq2}|m_\ell(\lambda_i)|^2
			\leq
			\frac{8\beta_W^2}{\Delta_\lambda}.
		\end{equation}
	\end{lemma}
	\begin{proof}
		See Appendix \ref{app:gt_nonsym_envelope}.
	\end{proof}

	The next lemma is used in Corollary~\ref{cor:gt_nonsym_two_region}
	for the region where only the eigenvalue magnitudes are bounded.
	
	\begin{lemma}[\bfseries\small Modulus-only propagation bounds]
		\label{lem:gt_nonsym_modulus}
		Let Assumption~\ref{assump:gt_nonsym_matrix} hold and let $\omega$ be
		defined in \eqref{gt_nonsym_omega}. Then
		\begin{equation}\label{gt_nonsym_disc_sums}
			\sum_{\ell=0}^{\infty}
			\omega^{-\ell}
			\max_{i\geq2}|m_\ell(\lambda_i)|
			\leq
			\frac{5}{\Delta_\lambda^2},
			\qquad
			\sum_{\ell=0}^{\infty}
			\max_{i\geq2}|m_\ell(\lambda_i)|^2
			\leq
			\frac{
				2(1-\Delta_\lambda)^2
			}{
				\Delta_\lambda^3(2-\Delta_\lambda)
			}
			\leq
			\frac1{\Delta_\lambda^3}.
		\end{equation}
	\end{lemma}
	\begin{proof}
		See Appendix \ref{app:gt_nonsym_modulus}.
	\end{proof}

	We now translate the modal bounds into bounds on the matrix propagation
	gains.
	
	\begin{lemma}[\bfseries\small Network gains]
		\label{lemma:gt_nonsym_gain_bounds}
		Let Assumption~\ref{assump:gt_nonsym_matrix} hold and let $\omega$ be
		defined in \eqref{gt_nonsym_omega}. Then
		\begin{align}\label{gt_nonsym_gain_bounds}
			\mathcal G_1(\omega)
			&\leq\frac{11\kappa_W\beta_W}{\Delta_\lambda},
			\qquad
			\mathcal G_2
			\leq\frac{8\kappa_W^2\beta_W^2}{\Delta_\lambda}.
		\end{align}
		Moreover,
		\begin{equation}\label{gt_nonsym_avg_gain}
			\mathcal G_{2,\mathrm{av}}
			\leq\frac{2\kappa_W^2}{N}\sum_{i=2}^N
			\frac{|\lambda_i|^2|1-\lambda_i|^2}
			{(1-|\lambda_i|^2)^3},
		\end{equation}
		with equality in \eqref{gt_nonsym_avg_gain} when $W$ is normal.
	\end{lemma}
	
	\begin{proof}
		Let $
		D\define\diag\{0,\lambda_2,\ldots,\lambda_N\}$.
		Since $\mathbf M_\ell=m_\ell(\widetilde{\W})$,
		$
		\mathbf M_\ell
		=
		(S\otimes I_d)
		\bigl(m_\ell(D)\otimes I_d\bigr)
		(S^{-1}\otimes I_d),
		$
		and hence
		\begin{equation}\label{gt_nonsym_B_eig_bound}
			\|\mathbf M_\ell\|_2
			\leq
			\kappa_W
			\max_{2\leq i\leq N}|m_\ell(\lambda_i)|.
		\end{equation}
		Combining \eqref{gt_nonsym_B_eig_bound} with
		\eqref{gt_nonsym_envelope_sums} and the definitions
		\eqref{gt_nonsym_network_gains} gives
		\eqref{gt_nonsym_gain_bounds}.
		The same representation gives
		$\|\mathbf M_\ell\|_F^2
		\leq d\kappa_W^2\sum_{i=2}^N|m_\ell(\lambda_i)|^2$.
		Summing and applying Lemma~\ref{lem:gt_nonsym_modal_energy}
		proves \eqref{gt_nonsym_avg_gain}, with equality for normal $W$
		by unitary invariance of the Frobenius norm.
		The comparison $\mathcal G_{2,\mathrm{av}}\leq\mathcal G_2$
		follows from
		$\|\mathbf M_\ell\|_F^2\leq Nd\|\mathbf M_\ell\|_2^2$.
	\end{proof}
	
	The bounds in Lemma~\ref{lemma:gt_nonsym_gain_bounds} show how the
	eigenvalue locations affect the network dependence. If
	$\beta_W=\cO(1)$, then
	$\mathcal G_1(\omega)=\cO(\kappa_W/\Delta_\lambda)$ and
	$\mathcal G_2=\cO(\kappa_W^2/\Delta_\lambda)$. In contrast, the
	general bound $\beta_W\leq2/\Delta_\lambda$ recovers the worst-case
	scaling
	$\mathcal G_2=\cO(\kappa_W^2/\Delta_\lambda^3)$.
	Thus the improvement depends on the eigenvalue locations, while
	$\kappa_W$ captures the effect of nonnormality.
	
	The bounds on $\mathcal G_1$ and $\mathcal G_2$ use global spectral
	parameters, which can be conservative when different parts of
	the spectrum exhibit different geometry. We therefore split the
	spectrum into two useful regions.
	
	\begin{corollary}[\bfseries\small Two-region refinement]
		\label{cor:gt_nonsym_two_region}
		Under Assumption~\ref{assump:gt_nonsym_matrix},
		suppose that $\{2,\ldots,N\}$ can be partitioned into
		$\mathcal I_1$ and $\mathcal I_2$ such that, for some
		$p,\nu\in(0,1]$ and $c\geq1$,
		\begin{equation}\label{gt_nonsym_two_region_conditions}
			\begin{array}{lll}
				1-|\lambda_i|\geq p,
				&
				\displaystyle
				\frac{|1-\lambda_i|}{1-|\lambda_i|}\leq c,
				&
				i\in\mathcal I_1,
				\\[2mm]
				1-|\lambda_i|\geq\nu,
				&
				&
				i\in\mathcal I_2.
			\end{array}
		\end{equation}
		Set
		\begin{equation}\label{gt_nonsym_two_region_omega}
			\omega_{\mathrm{reg}}
			\define
			1-\frac{\min\{p,\nu\}}{4}.
		\end{equation}
		Then
		\begin{equation}\label{gt_nonsym_two_region_gains}
			\mathcal G_1(\omega_{\mathrm{reg}})
			\leq
			\kappa_W
			\left(
			\frac{11c}{p}
			+\frac5{\nu^2}
			\right),
			\qquad
			\mathcal G_2
			\leq
			\kappa_W^2
			\left(
			\frac{8c^2}{p}
			+\frac1{\nu^3}
			\right).
		\end{equation}
		With $\theta_j\define|\mathcal I_j|/N$, $j=1,2$,
		\begin{equation}\label{gt_nonsym_two_region_avg_gain}
			\mathcal G_{2,\mathrm{av}}
			\leq\kappa_W^2
			\left(
			\frac{\theta_1c^2}{4p}
			+\frac{\theta_2}{\nu^3}
			\right).
		\end{equation}
		For real spectra in $[-1+\nu,1-p]$, the nonnegative/negative
		partition gives $c=1$; for symmetric $W$, also $\kappa_W=1$.
	\end{corollary}
	\begin{proof}
		See Appendix \ref{app:gt_nonsym_two_region}.
	\end{proof}

\begin{figure}[t]
		\centering
		\tikzexternaldisable
		\resizebox{0.6\linewidth}{!}{%
			\begin{tikzpicture}[>=Latex,scale=3.55,every node/.style={font=\footnotesize}]

\def\pGap{0.12}
\def\nuGap{0.45}
\def\cLoc{2.2}
\pgfmathsetmacro{\rP}{1-\pGap}
\pgfmathsetmacro{\rNu}{1-\nuGap}

\colorlet{myblue}{blue!75!black}
\colorlet{myorange}{orange!85!black}
\colorlet{mygray}{gray!65}

\newcommand{\radc}[1]{%
  ((\cLoc*\cLoc-cos(#1))
   -sqrt((\cLoc*\cLoc-cos(#1))*(\cLoc*\cLoc-cos(#1))
         -(\cLoc*\cLoc-1)*(\cLoc*\cLoc-1)))
  /(\cLoc*\cLoc-1)%
}

\fill[myorange!20] (0,0) circle (\rNu);
\path[fill=myblue!20,draw=none]
  plot[domain=-180:180,samples=361,variable=\t]
  ({min(\rP,\radc{\t})*cos(\t)},{min(\rP,\radc{\t})*sin(\t)}) -- cycle;

\draw[myorange,thick] (0,0) circle (\rNu);
\draw[myblue,thick]
  plot[domain=-180:180,samples=361,variable=\t]
  ({min(\rP,\radc{\t})*cos(\t)},{min(\rP,\radc{\t})*sin(\t)}) -- cycle;
\draw[mygray,densely dashed] (0,0) circle (1);
\draw[myblue!55,densely dashed] (0,0) circle (\rP);

\draw[->,thick] (-1.12,0) -- (1.22,0) node[anchor=north east,inner sep=1pt] {$\operatorname{Re}$};
\draw[->,thick] (0,-1.08) -- (0,1.14) node[left] {$\operatorname{Im}$};

\fill (0,0) circle (0.016);
\fill[myorange] (\rNu,0) circle (0.013);
\fill[myblue] (\rP,0) circle (0.013);
\fill (1,0) circle (0.016);
\node[anchor=north east,inner sep=1.5pt] at (0,0) {$0$};
\node[below=2pt] at (\rNu,0) {$1-\nu$};
\node[anchor=south east,inner sep=1pt,fill=white,fill opacity=.85,text opacity=1] at (\rP-0.01,0.03) {$1-p$};
\node[anchor=north west,inner sep=1.5pt] at (1,0) {$1$};

\draw[myblue!70,densely dotted] (\rP,0.015) -- (\rP,0.30);
\draw[black!50,densely dotted] (1,0.015) -- (1,0.30);
\draw[{Latex[length=1.2mm,width=1mm]}-{Latex[length=1.2mm,width=1mm]},myblue,thick] (\rP,0.27) -- (1,0.27)
  node[midway,above=1pt,fill=white,inner sep=1pt,text=black] {$p$};
\draw[myorange!75,densely dotted] (\rNu,-0.015) -- (\rNu,-0.34);
\draw[black!50,densely dotted] (1,-0.015) -- (1,-0.34);
\draw[{Latex[length=1.5mm,width=1.2mm]}-{Latex[length=1.5mm,width=1.2mm]},myorange,thick] (\rNu,-0.31) -- (1,-0.31)
  node[midway,below=1pt,fill=white,inner sep=1pt,text=black] {$\nu$};

\node[anchor=south west,inner sep=1pt] at (0.72,0.74) {$|z|=1$};
\node[anchor=south east,inner sep=1pt] at (-0.60,0.78) {$|z|=1-p$};
\draw[myblue!55,densely dashed,->] (-0.62,0.77) -- ({\rP*cos(132)},{\rP*sin(132)});

\node[align=center,fill=white,fill opacity=.85,text opacity=1,inner sep=1.5pt]
  at (-0.36,0.30) {\textcolor{myorange}{Region 2}\\$|z|\le 1-\nu$};
\node[align=center,fill=white,fill opacity=.85,text opacity=1,inner sep=1.5pt]
  at (0.56,-0.58) {\textcolor{myblue}{Region 1}\\$|z|\le1-p$,\\$|1-z|\le c(1-|z|)$};

\end{tikzpicture}
		}
		\tikzexternalenable
		\caption{Illustration of the two-region spectral refinement
			in the complex plane, with $p=0.12$, $\nu=0.45$, and $c=2.2$.
			The geometric regions may overlap; the partition is over
			eigenvalue indices.}
		\label{fig:two_region_complex}
	\end{figure}
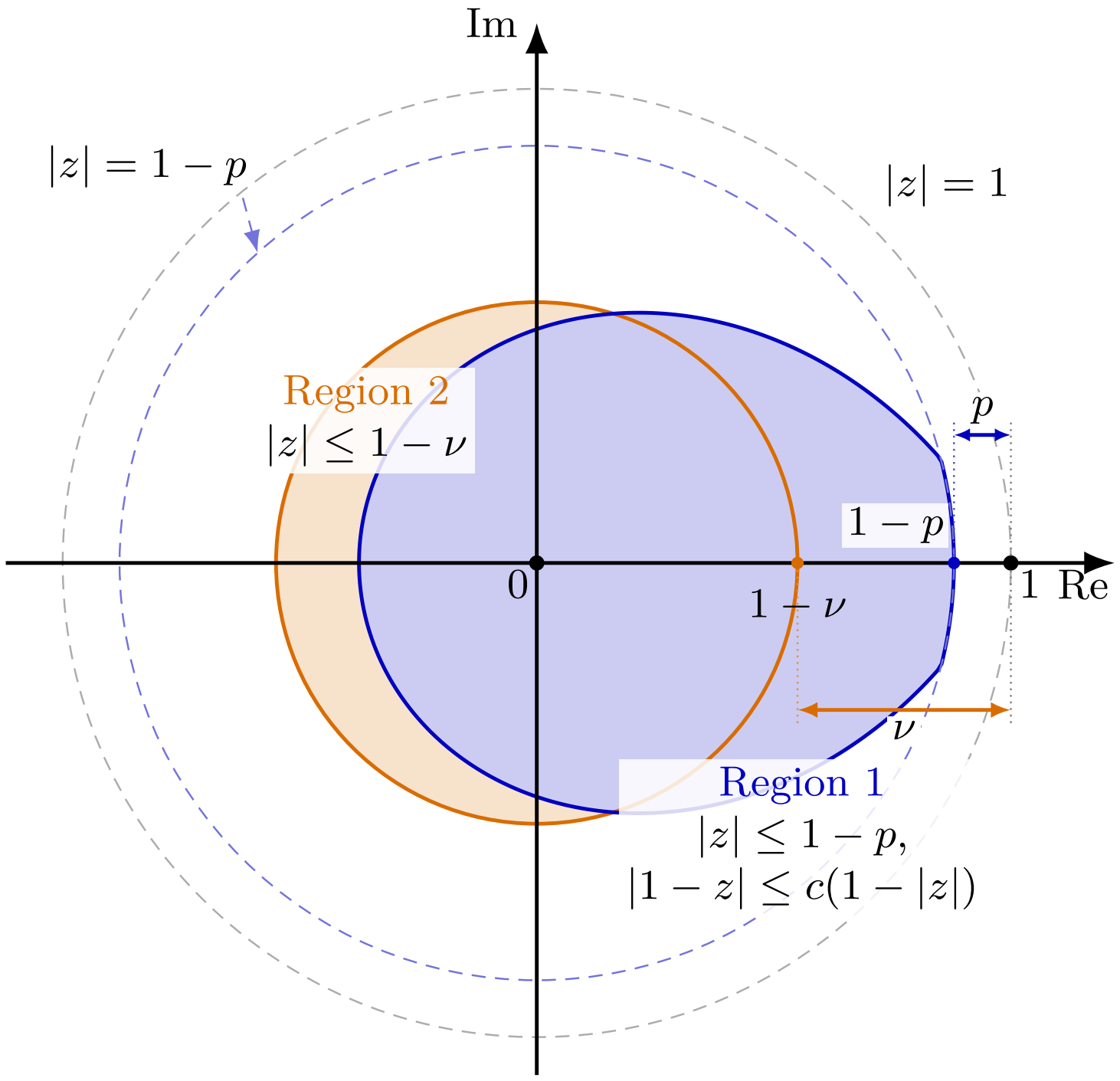

	The ratio
	$
	\frac{|1-\lambda_i|}{1-|\lambda_i|}
	$
	measures how unfavorable the location of $\lambda_i$ is relative to
	its distance from the unit circle. It equals one for
	$\lambda_i\in[0,1)$, but can be large for negative or complex
	eigenvalues near the unit circle.  However, for $i\in\mathcal I_1$, the condition
	$1-|\lambda_i|\geq p$ ensures
	\[
	\frac{|1-\lambda_i|}{1-|\lambda_i|}
	\leq
	\frac{1+|\lambda_i|}{1-|\lambda_i|}
	\leq
	\frac{2-p}{p},
	\]
	so $c$ is finite. The two-region decomposition therefore keeps the
	location-dependent factor moderate on $\mathcal I_1$ and handles
	spectrally unfavorable modes in $\mathcal I_2$ using the
	modulus-only bound. See Figure~\ref{fig:two_region_complex} for illustration of the two regions.
	
	\subsection{Weighted optimization and disagreement bounds}
	
	\label{sec:gt_nonsym_weighted}
	
	The following bounds use the average-recursion and gradient-perturbation
	estimates in Lemma~\ref{lemma:gt_shared_estimates} of
	Appendix~\ref{app:shared_estimates}.
	We next combine the exact representation
	\eqref{gt_nonsym_convolution_identity} with exponentially weighted sums.
	The function-value gap is used as an auxiliary quantity to control the
	optimization perturbation through cocoercivity without introducing an
	additional factor of $1/\mu$.
	
	Fix $K\geq1$ and define
	\begin{subequations}\label{gt_nonsym_gamma}
		\begin{align}
			\gamma&\define1-\frac{\mu\alpha}{2},
			\qquad
			\tilde f_k\define
			\Ex[f(\bar x^k)-f(x^\star)],
			\label{gt_nonsym_gamma_a}\\
			E_k&\define\frac1N\Ex\|\widehat{\x}^k\|^2,
			\qquad
			\mathcal D_K
			\define
			\sum_{k=0}^{K}\gamma^{-k}E_k.
			\label{gt_nonsym_weighted_sums}
		\end{align}
	\end{subequations}
	
	\begin{lemma}[\bfseries\small Weighted function-gap inequality]
		\label{lemma:gt_nonsym_weighted_gap}
		Let Assumptions~\ref{assump:gt_semi_functions} and
		\ref{assump:gt_semi_noise} hold, let $W$ be doubly stochastic, and let
		$0<\alpha\leq1/(4L)$.  Then
		\begin{equation}\label{gt_nonsym_weighted_gap}
			\sum_{k=0}^{K}\gamma^{-k}\tilde f_k
			\leq
			\frac{\bar e_0}{\alpha}
			+3L\mathcal D_K
			+\frac{4\sigma^2}{\mu N}\gamma^{-K}.
		\end{equation}
	\end{lemma}
	\begin{proof}
		See Appendix \ref{app:gt_nonsym_weighted_gap}.
	\end{proof}

	\begin{lemma}[\bfseries\small Weighted disagreement inequality]
		\label{lemma:gt_nonsym_weighted_disagreement}
		Let Assumptions~\ref{assump:gt_semi_functions},
		\ref{assump:gt_semi_noise}, and
		\ref{assump:gt_nonsym_matrix} hold. 		Let $\g^0=\grad\F(\x^0;\bxi^0)$, so that
		$\widetilde{\y}^0=\grad\f(\x^\star)$, suppose $0<\alpha\leq1/(4L)$ and
		\begin{equation}\label{gt_nonsym_alpha_g}
			\mu\alpha\leq\frac{\Delta_\lambda}{2}.
		\end{equation}
		Then $\gamma\geq\omega\define1-\Delta_\lambda/4$, and
		\begin{equation}\label{gt_nonsym_weighted_disagreement}
			\begin{aligned}
				\mathcal D_K
				&\leq 3\mathcal H_0
				+6\alpha^2\mathcal G_1(\omega)^2
				\left(
				2L^2\mathcal D_K
				+4L\sum_{k=0}^{K}\gamma^{-k}\tilde f_k
				\right)
				+
				\frac{6\alpha\sigma^2\mathcal G_{2,\mathrm{av}}}{\mu}
				\gamma^{-K},
			\end{aligned}
		\end{equation}
		where
		\begin{equation}\label{gt_nonsym_H0}
			\mathcal H_0
			\define
			2\kappa_W^2
			\left(
			\frac{E_0}{\Delta_\lambda}
			+\frac{2\alpha^2\varsigma_\star^2}{\Delta_\lambda^3}
			\right).
		\end{equation}
	\end{lemma}
	\begin{proof}
		See Appendix \ref{app:gt_nonsym_weighted_disagreement}.
	\end{proof}

	\subsection{Refined finite-horizon rates}
	\label{sec:gt_nonsym_rate}
	
Combining Lemmas~\ref{lemma:gt_nonsym_weighted_gap} and~\ref{lemma:gt_nonsym_weighted_disagreement} yields the finite-horizon convergence bound. Throughout, the quantity $1/(\sqrt{168}\,L\mathcal G_1(\omega))$ is taken to be $+\infty$ when $\mathcal G_1(\omega)=0$.
	
	\begin{subequations}
		\begin{theorem}[\bfseries\small Refined finite-horizon bound]
			\label{thm:gt_nonsym_refined}
			Let Assumptions~\ref{assump:gt_semi_functions},
			\ref{assump:gt_semi_noise}, and \ref{assump:gt_nonsym_matrix} hold.
			Initialize
			$\g^0=\grad\F(\x^0;\bxi^0)$, so that
			$\y^0=\zero$ and
			$\widetilde{\y}^0=\grad\f(\x^\star)$, and suppose $\omega\define1-\frac{\Delta_\lambda}{4}$ and
			\begin{equation}\label{gt_nonsym_stepsize_conditions}
				0<\alpha\leq
				\min\left\{
				\frac{1}{4L},
				\frac{\Delta_\lambda}{2\mu},
				\frac{1}{\sqrt{168}\,L\mathcal G_1(\omega)}
				\right\}.
			\end{equation}
			Then, for every $K\geq1$,
			\begin{equation}\label{gt_nonsym_main_neighborhood}
				\bar e_K
				\leq
				\exp\!\left(-\frac{\mu\alpha K}{2}\right)
				\bigl(\mathcal I_0+\widetilde c_2\alpha^3\bigr)
				+c_1\alpha+c_2\alpha^2,
			\end{equation}
			where
			\begin{equation}\label{gt_nonsym_constants}
				\mathcal I_0
				\define
				2\bar e_0+\frac{9\kappa_W^2}{\Delta_\lambda}E_0,
				\qquad
				\widetilde c_2
				\define
				\frac{72L\kappa_W^2\varsigma_\star^2}{\Delta_\lambda^3},
			\end{equation}
			and $c_1\define\frac{5\sigma^2}{\mu N}$,
			$c_2\define\frac{36L\sigma^2\mathcal G_{2,\mathrm{av}}}{\mu}$.
		\end{theorem}
		\begin{proof}
			See Appendix \ref{app:gt_nonsym_refined}.
		\end{proof}
	\end{subequations}
	
	For normal $W$, $\mathcal G_2$ sums the largest squared modal
	coefficient at each lag, while $\mathcal G_{2,\mathrm{av}}$
	averages these coefficients over the spectrum. Thus
	$\mathcal G_2/N\leq\mathcal G_{2,\mathrm{av}}\leq\mathcal G_2$:
	retaining modal multiplicities can reduce the stochastic network
	bound by up to a factor $N$, as illustrated by
	Example~\ref{ex:gt_nonsym_directed_ring}.
	The stepsize and initialization bounds are unchanged.
	On $W_{p,\nu}$ of Section~\ref{sec:sharp_network_dependence},
	a constant fraction of the modes lies at each spectral end,
	so $\mathcal G_{2,\mathrm{av}}\asymp\mathcal G_2$ and the stochastic
	network dependence remains sharp.
	
	\begin{remark}[\bfseries\small Non-diagonalizable mixing matrices]
		\rm
		Diagonalizability is used only to bound $\|\widetilde{\W}^k\|_2$ and
		$\mathbf M_\ell$ through the spectrum;
		\eqref{gt_nonsym_convolution_identity} does not require it. If only
		$\rho(\widetilde W)<1$, then $\|\widetilde{\W}^k\|_2\le C_rr^k$ for any
		$r\in(\rho(\widetilde W),\omega)$, so $\mathcal G_1(\omega)$,
		$\mathcal G_2$, $Q_0\define\sum_k\omega^{-k}\|\widetilde{\W}^k\|_2^2$, and
		$Q_2\define\sum_kk^2\omega^{-k}\|\widetilde{\W}^k\|_2^2$ are finite, and
		Lemma~\ref{lemma:gt_nonsym_weighted_disagreement} and
		Theorem~\ref{thm:gt_nonsym_refined} hold with $\kappa_W^2/\Delta_\lambda$
		and $2\kappa_W^2/\Delta_\lambda^3$ replaced by $Q_0$ and $Q_2$ in
		$\mathcal H_0$, $\mathcal I_0$, and $\widetilde c_2$; explicit spectral
		bounds then require estimates of the Jordan blocks or of nonnormality.
		Theorem~\ref{thm:gt_baseline_convergence} already covers
		non-diagonalizable $W$ with $\|W-\Pi\|_2<1$ at the cubic level.
	\end{remark}
	\begin{corollary}[\bfseries\small Refined neighborhood]
		\label{cor:gt_nonsym_neighborhood}
		Under the conditions of Theorem~\ref{thm:gt_nonsym_refined},
		the mean-square error of the average iterate is bounded by a term
		decaying geometrically with factor $1-\mu\alpha/2$ plus a constant.
		In particular,
		\begin{equation}\label{gt_nonsym_explicit_neighborhood}
			\limsup_{k\rightarrow\infty}
			\Ex\|\bar x^k-x^\star\|^2
			=
			\cO\!\left(
			\frac{\alpha\sigma^2}{\mu N}
			+
			\frac{\alpha^2L\sigma^2}{\mu}\mathcal G_{2,\mathrm{av}}
			\right).
		\end{equation}
	\end{corollary}
	
	\begin{corollary}[\bfseries\small Tuned finite-horizon rate]
		\label{cor:gt_nonsym_linear_speedup}
		Let the conditions of Theorem~\ref{thm:gt_nonsym_refined} hold
		and assume $\sigma^2>0$. Define
		\begin{equation}\label{gt_nonsym_alpha_bar}
			\bar\alpha
			\define
			\min\left\{
			\frac{1}{4L},\,
			\frac{\Delta_\lambda}{2\mu},\,
			\frac{1}
			{\sqrt{168}\,L\mathcal G_1(\omega)}
			\right\}.
		\end{equation}
		For $K\geq1$, choose $\alpha$ according to
		\eqref{eq:tuning_alpha}. Then
		\begin{equation}\label{gt_nonsym_big_o_neighborhood}
			\bar e_K
			=
			\widetilde\cO\!\left(
			\frac{\sigma^2}{\mu^2NK}
			+
			\frac{L\sigma^2\mathcal G_{2,\mathrm{av}}}{\mu^3K^2}
			\right)
			+
			C_{\mathrm{init}}
			\exp\!\left(
			-\frac{\mu\bar\alpha K}{2}
			\right),
		\end{equation}
		where $
		C_{\mathrm{init}}
		\define
		2\bar e_0
		+\frac{9\kappa_W^2}{\Delta_\lambda}E_0
		+\frac{
			72L\kappa_W^2\varsigma_\star^2
		}{
			\Delta_\lambda^3
		}
		\bar\alpha^3$.
	\end{corollary}
	\begin{proof}
		See Appendix \ref{app:gt_nonsym_linear_speedup}.
	\end{proof}

	\begin{corollary}[\bfseries\small Two-region finite-horizon rate]
		\label{cor:gt_nonsym_two_region_rate}
		Let the conditions of Theorem~\ref{thm:gt_nonsym_refined} and
		Corollary~\ref{cor:gt_nonsym_two_region} hold, and assume
		$\sigma^2>0$ with $\bar\alpha$ replaced by
		\begin{equation}\label{gt_nonsym_two_region_alpha_bar}
			\bar\alpha
			\define
			\min\left\{
			\frac{1}{4L},\,
			\frac{\min\{p,\nu\}}{2\mu},\,
			\frac{1}{
				\sqrt{168}\,L\kappa_W
				\left(
				11c/p+5/\nu^2
				\right)
			}
			\right\}.
		\end{equation}
		For $K\geq1$, choose $\alpha$ according to
		\eqref{eq:tuning_alpha}. Then
		\begin{equation}\label{gt_nonsym_two_region_rate}
			\bar e_K
			=
			\widetilde\cO\!\left(
			\frac{\sigma^2}{\mu^2NK}
			+
			\frac{
				\kappa_W^2L\sigma^2
			}{
				\mu^3K^2
			}
			\left(
			\frac{\theta_1c^2}{p}
			+
			\frac{\theta_2}{\nu^3}
			\right)
			\right)
			+
			C_{\mathrm{init}}
			\exp\!\left(
			-\frac{\mu\bar\alpha K}{2}
			\right).
		\end{equation}
		Moreover, the network-dependent $K^{-2}$ term is dominated by the
		leading $K^{-1}$ term, up to logarithmic factors, once
		\begin{equation}\label{gt_nonsym_two_region_total_transient}
			K
			\geq
			K_{\mathrm{trans}}
			\define
			\widetilde\cO\!\left(
			N\kappa\kappa_W^2
			\left(
			\frac{\theta_1c^2}{p}
			+
			\frac{\theta_2}{\nu^3}
			\right)
			\right).
		\end{equation}
	\end{corollary}
	
	\begin{proof}
		Let $\omega_0=1-\Delta_\lambda/4$. For each nonempty
		$\mathcal I_j$, define
		\[
		d_j=\min_{i\in\mathcal I_j}(1-|\lambda_i|),
		\qquad
		a_\ell^{(j)}=\max_{i\in\mathcal I_j}|m_\ell(\lambda_i)|.
		\]
		For an empty region, set $d_j=1$ and $a_\ell^{(j)}=0$.
		Then $d_1\geq p$, $d_2\geq\nu$, and
		$\Delta_\lambda\leq d_j$. Consequently,
		$\omega_0\geq1-d_j/4$.
		
		Applying the scalar envelope arguments of
		Lemmas~\ref{lem:gt_nonsym_envelope} and
		\ref{lem:gt_nonsym_modulus} to the two regions separately gives
		\[
		\sum_{\ell\geq0}\omega_0^{-\ell}a_\ell^{(1)}
		\leq \frac{11c}{d_1}\leq\frac{11c}{p},
		\qquad
		\sum_{\ell\geq0}\omega_0^{-\ell}a_\ell^{(2)}
		\leq \frac{5}{d_2^2}\leq\frac{5}{\nu^2}.
		\]
		Since
		$\|\mathbf M_\ell\|_2
		\leq\kappa_W(a_\ell^{(1)}+a_\ell^{(2)})$, it follows that
		\[
		\mathcal G_1(\omega_0)
		\leq\kappa_W\left(\frac{11c}{p}+\frac{5}{\nu^2}\right).
		\]
		Together with $\Delta_\lambda\geq\min\{p,\nu\}$, this shows that
		\eqref{gt_nonsym_two_region_alpha_bar} is no larger than the
		stepsize cap in \eqref{gt_nonsym_alpha_bar}.
		Applying Corollary~\ref{cor:gt_nonsym_linear_speedup} with
		\eqref{gt_nonsym_two_region_avg_gain} gives
		\eqref{gt_nonsym_two_region_rate}.
		Comparing its network-dependent $K^{-2}$ term with the leading
		$K^{-1}$ term gives
		\eqref{gt_nonsym_two_region_total_transient}.
	\end{proof}
	
	The two spectral regions contribute \emph{additively}, weighted
	by their spectral proportions. The stochastic $K^{-2}$ term
	depends on $\theta_1c^2/p+\theta_2/\nu^3$, so a region contributes
	less when it contains fewer eigenvalues. Discarding the weights
	recovers the envelope $c^2/p+\nu^{-3}$. The stepsize remains
	governed by the operator-norm gain $\mathcal G_1$.
	
	\begin{corollary}[\bfseries\small Symmetric two-region specialization]
		\label{cor:gt_symmetric_two_region}
		Under the assumptions and stepsize choice of
		Corollary~\ref{cor:gt_nonsym_two_region_rate}, suppose that $W$ is
		symmetric and its nontrivial spectrum lies in $[-1+\nu,1-p]$. Take
		$\mathcal I_1=\{i\geq2:\lambda_i\geq0\}$,
		$\mathcal I_2=\{i\geq2:\lambda_i<0\}$, and
		$\theta_j=|\mathcal I_j|/N$.
		Then $\kappa_W=1$ and $c=1$ in
		Corollary~\ref{cor:gt_nonsym_two_region}. Consequently,
		\begin{equation}\label{gt_symmetric_two_region_rate}
			\begin{aligned}
				\bar e_K
				&=
				\widetilde\cO\!\left(
				\tfrac{\sigma^2}{\mu^2NK}
				+
				\tfrac{L\sigma^2}{\mu^3K^2}
				\left(
				\tfrac{\theta_1}{p}
				+\tfrac{\theta_2}{\nu^3}
				\right)
				\right)
				+
				C_{\mathrm{init}}
				\exp\!\left(
				-\tfrac{\mu\bar\alpha K}{2}
				\right).
			\end{aligned}
		\end{equation}
		In particular, $K_{\mathrm{trans}}
		=
		\widetilde\cO\!\left(
		N\kappa
		\left(
		\frac{\theta_1}{p}
		+\frac{\theta_2}{\nu^3}
		\right)
		\right)$.
	\end{corollary}
	The separation between the two bounds is largest when $p\ll\nu$.
	If $p=1-\lambda_2(W)<\nu$, then $\Delta_\lambda=p$,
	so a bound based only on $\Delta_\lambda$ can scale as $p^{-3}$,
	whereas \eqref{gt_symmetric_two_region_rate} has stochastic factor
	$\theta_1/p+\theta_2/\nu^3\leq p^{-1}+\nu^{-3}$.
	Thus, if $\nu=\Theta(1)$, the stochastic network dependence is
	at most linear in the inverse positive spectral gap.
	The negative-end penalty is likewise weighted by the fraction
	of negative eigenvalues.
	
	For symmetric mixing, replacing $W$ by $(W+I)/2$ makes all
	eigenvalues nonnegative, so $\nu=1$, while decreasing $p$
	by at most a factor of two.
	
	The eigenvector-conditioning penalty also disappears for the broader
	class of normal mixing matrices, which includes nonsymmetric circulant
	matrices.
	\begin{corollary}[\bfseries\small Normal mixing matrices]
		\label{cor:gt_nonsym_normal}
		Suppose $W$ is nonnegative, primitive, doubly stochastic, and normal,
		i.e., $WW\tran=W\tran W$. Then
		Assumption~\ref{assump:gt_nonsym_matrix} holds with $\kappa_W=1$, so
		that, for $\omega=1-\Delta_\lambda/4$,
		\begin{equation}\label{gt_nonsym_normal_gains}
			\mathcal G_1(\omega)
			\leq
			\frac{11\beta_W}{\Delta_\lambda},
			\qquad
			\mathcal G_2
			\leq
			\frac{8\beta_W^2}{\Delta_\lambda},
		\end{equation}
		and, under the conditions of
		Corollary~\ref{cor:gt_nonsym_linear_speedup} with
		$\bar\alpha=\min\{\frac{1}{4L},\frac{\Delta_\lambda}{2\mu},
		\frac{\Delta_\lambda}{11\sqrt{168}\,L\beta_W}\}$,
		\begin{equation}\label{gt_nonsym_normal_rate}
			\bar e_K
			=
			\widetilde\cO\!\left(
			\frac{\sigma^2}{\mu^2NK}
			+
			\frac{\beta_W^2L\sigma^2}
			{\mu^3\Delta_\lambda K^2}
			\right)
			+
			C_{\mathrm{init}}
			\exp\!\left(
			-\frac{\mu\bar\alpha K}{2}
			\right).
		\end{equation}
	\end{corollary}
	\begin{proof}
		Double stochasticity gives $\widetilde W\widetilde W\tran=WW\tran-\Pi$
		and $\widetilde W\tran\widetilde W=W\tran W-\Pi$, so $\widetilde W$ is
		normal and unitarily diagonalizable; choosing $S$ unitary in
		\eqref{gt_nonsym_diagonalization} gives $\kappa_W=1$, and primitivity
		gives $|\lambda_i|<1$. Then \eqref{gt_nonsym_normal_gains} follows from
		Lemma~\ref{lemma:gt_nonsym_gain_bounds}, and
		\eqref{gt_nonsym_normal_rate} from
		Corollary~\ref{cor:gt_nonsym_linear_speedup}.
	\end{proof}
	Thus normality removes the eigenvector-conditioning penalty entirely,
	and when $\beta_W=\cO(1)$ the $K^{-2}$ term has the linear-gap
	dependence $\widetilde\cO(\Delta_\lambda^{-1}K^{-2})$.
	
	\begin{remark}[\bfseries\small Transient time for linear speedup]
		\label{rem:gt_nonsym_transient}
		\rm
		In this remark, we treat $L$ and $\mu$ as fixed and suppress their
		dependence in the rates and transient bounds. In \eqref{gt_nonsym_big_o_neighborhood}, the $K^{-2}$ term is
		dominated by the leading $K^{-1}$ term after
		$\widetilde\cO(N\mathcal G_{2,\mathrm{av}})$ iterations.
		The exponential term has decay scale $(\mu\bar\alpha)^{-1}$.
		Thus a sufficient transient time for the rate
		$\widetilde\cO(\sigma^2/(NK))$ is
		$\widetilde\cO(\max\{N\mathcal G_{2,\mathrm{av}},
		(\mu\bar\alpha)^{-1}\})$.
		The dependence on $\varsigma_\star^2$ enters only logarithmically
		through $C_{\mathrm{init}}=\mathcal I_0+
		72L\kappa_W^2\bar\alpha^3\varsigma_\star^2/\Delta_\lambda^3$.
		
		To make the spectral dependence explicit, let
		\[
		S_W\define\sum_{i=2}^N
		\frac{2|\lambda_i|^2|1-\lambda_i|^2}{(1-|\lambda_i|^2)^3}.
		\]
		By \eqref{gt_nonsym_avg_gain},
		$N\mathcal G_{2,\mathrm{av}}\leq\kappa_W^2 S_W$,
		with equality for normal $W$ when $\kappa_W=1$.
		The general cap \eqref{gt_nonsym_alpha_bar} and the bound on
		$\mathcal G_1(\omega)$ in \eqref{gt_nonsym_gain_bounds} give
		\begin{equation}\label{gt_nonsym_transient_total}
			K_{\mathrm{trans}}
			=\widetilde\cO\!\left(
			\kappa_W\max\!\left\{
			\kappa_W S_W,\frac{\beta_W}{\Delta_\lambda}
			\right\}\right).
		\end{equation}	
		Under the two-region conditions, the stepsize cap
		\eqref{gt_nonsym_two_region_alpha_bar} instead gives
		\begin{equation}\label{gt_nonsym_two_region_transient_cap}
			K_{\mathrm{trans}}
			=\widetilde\cO\!\left(
			\kappa_W\max\!\left\{
			\kappa_W S_W,\frac{c}{p},\frac{1}{\nu^2}
			\right\}\right).
		\end{equation}
		Here $\max\{c/p,\nu^{-2}\}$ controls the stepsize since
		$11c/p+5/\nu^2\leq16\max\{c/p,\nu^{-2}\}$.
		If only the region sizes and bounds are known, the modal estimates
		in Corollary~\ref{cor:gt_nonsym_two_region} give
		\[
		S_W\leq N\left(\frac{\theta_1c^2}{4p}
		+\frac{\theta_2}{\nu^3}\right).
		\]
		Thus the stochastic contribution retains the number of modes in
		each region. For symmetric $W$ with nontrivial spectrum in
		$[-1+\nu,1-p]$, take $\kappa_W=c=1$.
	\end{remark}
	
	\begin{example}[\bfseries\small Lazy directed ring]
		\label{ex:gt_nonsym_directed_ring}
		\rm
		Let $P$ denote the cyclic-shift permutation matrix on $N$ nodes,
		$P_{i,i+1}=1$ for $i<N$, $P_{N,1}=1$, and zero otherwise, so that each
		agent receives information from one neighbor along a directed ring, and
		define
		\begin{equation}\label{gt_nonsym_ring}
			W=(1-\epsilon)P+\epsilon I,
			\qquad
			\epsilon\in(0,1).
		\end{equation}
		Then $W$ is nonnegative, primitive, doubly stochastic, and normal,
		with eigenvalues
		\[
		z_j
		=
		\epsilon
		+(1-\epsilon)
		\exp\!\left(\frac{2\pi\imath j}{N}\right),
		\qquad
		j=0,\ldots,N-1.
		\]
		Here $z_0=1$ is the consensus eigenvalue.
		Let
		\[
		\theta\define\frac{2\pi}{N},
		\qquad
		r_\theta
		\define
		\sqrt{
			1-2\epsilon(1-\epsilon)(1-\cos\theta)
		}.
		\]
		The largest nontrivial eigenvalue modulus is $r_\theta$, and hence
		\begin{equation}\label{gt_nonsym_ring_gap_exact}
			\Delta_\lambda
			=
			1-r_\theta
			=
			\frac{
				2\epsilon(1-\epsilon)(1-\cos\theta)
			}{
				1+r_\theta
			}.
		\end{equation}
		The largest modulus and the worst spectral-location ratio are attained at $j=1$ and $j=N-1$, giving
		\begin{equation}\label{gt_nonsym_ring_beta_exact}
			\beta_W
			=
			\frac{
				(1-\epsilon)|1-\exp(\imath\theta)|
			}{
				\Delta_\lambda
			}
			=
			\frac{
				1+r_\theta
			}{
				2\epsilon\sin(\theta/2)
			}.
		\end{equation}
		For fixed $\epsilon\in(0,1)$ and $N\rightarrow\infty$,
		\begin{equation}\label{gt_nonsym_ring_params}
			\Delta_\lambda
			\sim
			\frac{
				2\pi^2\epsilon(1-\epsilon)
			}{
				N^2
			},
			\qquad
			\beta_W
			\sim
			\frac{N}{\pi\epsilon}.
		\end{equation}
		Therefore,
		\begin{equation}\label{gt_nonsym_ring_gain}
			\mathcal G_2
			=
			\cO\!\left(
			\frac{\beta_W^2}{\Delta_\lambda}
			\right)
			=
			\cO\!\left(
			\frac{N^4}{\epsilon^3(1-\epsilon)}
			\right).
		\end{equation}
		For fixed $\epsilon\in(0,1)$ and $N\to\infty$,
		\eqref{gt_nonsym_ring_gain} and
		Lemma~\ref{lem:gt_nonsym_modal_energy} give
		$\mathcal G_2=\Theta(N^4)=\Theta(\Delta_\lambda^{-2})$.
		However, \eqref{gt_nonsym_avg_gain} yields
		$\mathcal G_{2,\mathrm{av}}
		=\Theta(N^3)=\Theta(\Delta_\lambda^{-3/2})$.
		Indeed, the energy of mode $z_j$ equals
		$|z_j|^2/[8\epsilon^3(1-\epsilon)\sin^4(\pi j/N)]$.
		Its sum over $j$ is $\Theta(N^4)$, by comparison with
		$\sum_{j\geq1}j^{-4}$, with the matching lower bound supplied
		by $j=1$. Dividing by $N$ gives the claimed averaged gain.
		Thus spectral averaging improves the stochastic network factor
		by order $N$. Constants may depend on $\epsilon$;
		the cubic lower bound remains a worst-case statement.
		
		In contrast, for fixed $N\geq2$ and $\epsilon\downarrow0$,
		$\Delta_\lambda=\Theta(\epsilon)$ and the dominant mode
		approaches $e^{2\pi\imath/N}\neq1$.
		The same modal identity and the upper bound then give
		$\mathcal G_2=\Theta(\Delta_\lambda^{-3})$,
		with implicit constants that may depend on $N$,
		recovering the worst-case cubic dependence. For odd $N$ all nontrivial eigenvalues are nonreal (for $N=3$ they
		approach $e^{\pm2\pi\imath/3}$), so negative real eigenvalues are not
		necessary for the cubic gain.
	\end{example}

	\section{Finite-horizon lower bound}
	\label{sec:sharp_network_dependence}
	
	For quadratics with the same curvature, disagreement cancels from the
	averaged gradient. We strengthen the lower bound in Section \ref{subsec:lb_two_agent_disagreement}
	with a finite-horizon symmetric construction that captures both ends of
	the spectrum and applies directly to the average iterate. The symmetric family below is sufficient to establish worst-case spectral lower bounds over the larger mixing-matrix class.
	Throughout, $d=1$, $L=3\mu$, $\Pi=\one\one\tran/N$, and
	$\widetilde W=W-\Pi$. The stepsize is constant within each run.
	
	\paragraph{Construction.}
	Let $N\geq4$ be a multiple of four and let
	$P_{\mathrm{pair}}\define I_{N/2}\otimes\bigl[\begin{smallmatrix}0&1\\1&0\end{smallmatrix}\bigr]$.
	Then $P_{\mathrm{pair}}x=x$ if $x_{2j-1}=x_{2j}$ for all $j$, and
	$P_{\mathrm{pair}}x=-x$ if $x_{2j-1}=-x_{2j}$ for all $j$. For
	$0<p\leq\nu\leq1/4$, we define
	\begin{equation}\label{eq:jlb_mixing}
		W_{p,\nu}=(1-\nu)P_{\mathrm{pair}}
		+\frac{\nu-p}{2}(I_N+P_{\mathrm{pair}})+p\Pi,
		\qquad s^{\mathrm{lb}}=\col\{1,1,-1,-1,\ldots\}.
	\end{equation}
	Write $W=W_{p,\nu}$. This matrix is symmetric, every entry is at least
	$p/N>0$, and doubly stochastic. The consensus vector $\one$ has eigenvalue $1$. On vectors that are constant within each pair,
	$x_{2j-1}=x_{2j}$ for all $j$, and whose entries sum to zero,
	$W$ has eigenvalue $z_+=1-p$. On vectors that have opposite values
	within each pair, $x_{2j-1}=-x_{2j}$ for all $j$, $W$ has eigenvalue
	$z_-=-1+\nu$.
	Since $p\leq\nu$, the largest nontrivial eigenvalue modulus is
	$1-p$, so $\Delta_\lambda=p$. Symmetry gives $\kappa_W=1$, and
	the chosen vector $s^{\mathrm{lb}}$ satisfies
	$Ws^{\mathrm{lb}}=(1-p)s^{\mathrm{lb}}$.

	For given $\mu>0$ and $\varsigma_\star\geq0$, we consider the local costs
	\begin{equation}\label{eq:glb_objectives}
		f_i(t)=\frac\mu2t^2+\mu(t_+)^2+\varsigma_\star s_i^{\mathrm{lb}}t,
		\quad
		f_i'(t)=2\mu t+\mu|t|+\varsigma_\star s_i^{\mathrm{lb}},
	\end{equation}
	where $t_+=\max\{t,0\}$.
	These costs are continuously differentiable, $\mu$-strongly convex, and
	$3\mu$-smooth. Because $\one\tran s^{\mathrm{lb}}=0$ and
	$\|s^{\mathrm{lb}}\|^2=N$, the minimizer of the aggregate cost is
	$x^\star=0$ and $N^{-1}\|\grad\f(\x^\star)\|^2=\varsigma_\star^2$. For
	$\sigma,\bar e_0\geq 0$, we run \eqref{gt_atc_semi} with noises that are
	independent across agents and iterations,
	\begin{equation}\label{eq:glb_noise_initialization}
		\grad F_i(t;\xi_i^k)=f_i'(t)+\xi_i^k,\quad
		\xi_i^k\sim\mathcal N(0,\sigma^2),\quad
		x_i^0=-\sqrt{\bar e_0},\quad \g^0=\grad\F(\x^0;\bxi^0).
	\end{equation}
	We also write $\rho_{\mathrm{lb}}=1-2\alpha\mu$ and
	$T_a=\lceil(4a)^{-1}\rceil$ for $a\in\{p,\nu\}$. 
	
	\begin{theorem}[\bfseries\small Joint finite-horizon lower bound]
		\label{thm:jlb_joint}
		For \eqref{eq:jlb_mixing}--\eqref{eq:glb_noise_initialization}, suppose
		\begin{equation}\label{eq:jlb_stepsize}
			0<\alpha L\leq\frac{\min\{p,\nu^2\}}4,\qquad L=3\mu.
		\end{equation}
		Then, for every $K\geq T_p+1$, we have
		\begin{align}\label{eq:jlb_main}
			\bar e_K\gtrsim{}&
			\bar e_0\rho_{\mathrm{lb}}^{2K}
			+\frac{\alpha\sigma^2}{\mu N}(1-\rho_{\mathrm{lb}}^{2K})
			+\alpha^2\sigma^2
			\left(1-\rho_{\mathrm{lb}}^{K-T_p}\right)^2
			\left(
			\frac1p+\frac1{\nu^3}
			\right)
			+\frac{\alpha^4\mu^2\varsigma_\star^2}{p^4}
			\rho_{\mathrm{lb}}^{2(K-T_p-1)}.
		\end{align}
	\end{theorem}
	\begin{proof}
		See Appendix~\ref{app:finite_lower_bound}.
	\end{proof}

	\begin{corollary}[\bfseries\small Horizon-dependent joint lower bound]
		\label{cor:jlb_tuned}
		Under the assumptions of Theorem~\ref{thm:jlb_joint}, for every
		$K\geq2T_p$ and every $\alpha$ satisfying \eqref{eq:jlb_stepsize},
		\begin{equation}\label{eq:jlb_tuned}
			\begin{aligned}
				\bar e_K\gtrsim{}&
				\min\left\{\bar e_0,\;
				\frac{\sigma^2}{\mu^2NK}
				+\frac{\sigma^2}{\mu^2K^2}
				\left(\frac1p+\frac1{\nu^3}\right)\right\}
				+\frac{\alpha^4\mu^2\varsigma_\star^2}{p^4}
				\rho_{\mathrm{lb}}^{2(K-T_p-1)}.
			\end{aligned}
		\end{equation}
	\end{corollary}
	
	\begin{proof}
		We retain the heterogeneity term in \eqref{eq:jlb_main} and bound the
		remaining terms. Since \eqref{eq:jlb_stepsize} gives
		$0<2\alpha\mu\leq\nu^2/6\leq1/96$, the elementary inequalities
		$e^{-2x}\leq1-x\leq e^{-x}$, valid for $0\leq x\leq1/2$, yield
		\[
		e^{-4\alpha\mu}\leq\rho_{\mathrm{lb}}\leq e^{-2\alpha\mu}.
		\]
		If $\alpha\mu K\leq1$, then
		$\rho_{\mathrm{lb}}^{2K}\geq e^{-8}$, so the initialization term in
		\eqref{eq:jlb_main} is at least a constant times $\bar e_0$, and the
		nonnegative noise terms may be dropped.
		Otherwise, $\alpha\geq1/(\mu K)$ and
		$1-\rho_{\mathrm{lb}}^{2K}\geq1-e^{-4}$. Moreover, $p\leq\nu$ gives
		$T_\nu\leq T_p\leq K/2$, so, for $a\in\{p,\nu\}$,
		\[
		1-\rho_{\mathrm{lb}}^{K-T_a}
		\geq1-e^{-2\alpha\mu(K-T_a)}
		\geq1-e^{-\alpha\mu K}
		\geq1-e^{-1}.
		\]
		Substituting into \eqref{eq:jlb_main} and using
		$\alpha\geq1/(\mu K)$, the stochastic terms are bounded below by a
		constant times
		\[
		\frac{\sigma^2}{\mu^2NK}
		+\frac{\sigma^2}{\mu^2K^2}
		\left(\frac1p+\frac1{\nu^3}\right).
		\]
		The two cases give the minimum in \eqref{eq:jlb_tuned}, while the
		heterogeneity term is preserved in both.
	\end{proof}
	
	\paragraph{Discussion.}
	On the family \eqref{eq:jlb_mixing}, take any stepsize satisfying
	\eqref{gt_nonsym_stepsize_conditions}. As shown below, it also
	satisfies \eqref{eq:jlb_stepsize}. Once
	$\alpha\mu(K-T_p)\gtrsim1$, the stochastic terms in the upper and
	lower bounds agree up to universal constants and have the form
	\[
	\frac{\alpha\sigma^2}{\mu N}
	+\alpha^2\sigma^2\left(\frac1p+\frac1{\nu^3}\right).
	\]
	On $W_{p,\nu}$, the eigenvalues $1-p$ and $-1+\nu$ have
	multiplicities $N/2-1$ and $N/2$, respectively. Thus
	$\theta_1=1/2-1/N\geq1/4$ and $\theta_2=1/2$, so the weighted
	upper factor $\theta_1/p+\theta_2/\nu^3$ is within a factor of four
	of $p^{-1}+\nu^{-3}$. Hence the lower bound matches the stochastic
	network dependence on this family. In particular,
	$p=\nu=\Delta_\lambda$ gives the worst-case cubic dependence
	$\Delta_\lambda^{-3}$.
	
	For brevity, write $m_\ell^\pm\define m_\ell(z_\pm)$.
	The lower-bound stepsize range also contains the range used by the refined
	upper bound on this family. Indeed,    since $z_\pm$ are eigenvalues of $\widetilde W$, $\|\mathbf M_\ell\|_2\geq|m_\ell^\pm|$
	and $\omega\leq1$, so $\mathcal G_1(\omega)\geq\sum_{\ell\geq0}|m_\ell^-|
	=2(1-\nu)/\nu^2\geq3/(2\nu^2)$ and $\mathcal G_1(\omega)\geq\sum_{\ell\geq0}|m_\ell^+|
	=2n_p(1-p)^{n_p}\geq1/(3p)$ with $n_p=\lfloor1/p\rfloor$, the last identity
	because $\sum_{\ell\geq0}m_\ell^+=0$ and $m_\ell^+\leq0$ exactly for $\ell<n_p$,
	and the last inequality from $n_pp\geq3/4$ and $(1-p)^{1/p}\geq(3/4)^4$ for $p\leq1/4$. Hence,
	\[
	\mathcal G_1(\omega)
	\geq
	\max\left\{
	\frac{1}{3p},
	\frac{3}{2\nu^2}
	\right\},
	\]
	so \eqref{gt_nonsym_stepsize_conditions} implies
	\[
	\alpha L
	\leq
	\min\left\{
	\frac{3p}{\sqrt{168}},
	\frac{2\nu^2}{3\sqrt{168}}
	\right\}
	\leq
	\frac{\min\{p,\nu^2\}}4,
	\]
	which is \eqref{eq:jlb_stepsize}.
	
	Likewise, \eqref{eq:jlb_tuned} matches the two polynomial stochastic
	terms in \eqref{gt_symmetric_two_region_rate}, up to logarithmic
	factors whenever the initialization is not smaller than these terms.
	The heterogeneity transients are not matched: the upper and lower bounds
	scale as $\alpha^3L\varsigma_\star^2/p^3$ and
	$\alpha^4\mu^2\varsigma_\star^2/p^4$, respectively. We therefore make
	no sharpness claim for the heterogeneity transient, condition-number
	dependence, or geometric decay rate.

	\section{Comparison with prior analyses}
	\label{sec:comparison}
	
	We compare our strongly convex results with
	\cite{alghunaim2023enhanced,koloskova2021improved} and use the lower
	bounds of Section~\ref{sec:sharp_network_dependence} to assess the
	sharpness of the network factors. Logarithmic factors and constants
	depending only on the optimization parameters are suppressed. 	Table~\ref{tab:gt_rate_comparison} compares the GT bounds.

	For symmetric mixing, let
	$q=1-\max_{i\geq2}\lambda_i^2=\Delta_\lambda(2-\Delta_\lambda)$.
	The DSGD result in \cite[Theorem~1]{takezawa2026improved} assumes
	smooth sample losses and bounds an exponentially weighted average
	of expected objective gaps. Its bound has leading term $\sigma_\star^2/(NK)$, initialization term
	$(\bar e_0/q)e^{-c_0qK/\kappa}$, and network terms
	\[
	\frac{1}{K^2}\left(
	\frac{\sigma_\star^2}{N}\sum_{i=2}^N
	\frac{\lambda_i^2}{1-\lambda_i^2}
	+\frac{(1-q)\varsigma_\star^2}{q^2}
	\right),
	\]
	up to logarithmic and problem-dependent factors. Here $\sigma_\star^2$
	bounds the noise variance at the optimum and $c_0>0$ is universal.
	The GT noise weight in \eqref{gt_nonsym_avg_gain} differs from the
	DSGD weight by a factor $2/(1+\lambda_i)^2$. Thus the weights are
	comparable for nonnegative eigenvalues, but near $-1$ they grow as
	$(1+\lambda_i)^{-1}$ for DSGD and $(1+\lambda_i)^{-3}$ for GT.
	For small $\Delta_\lambda$, the DSGD heterogeneity term scales as
	$\varsigma_\star^2/(\Delta_\lambda^2K^2)$, whereas in our GT bound
	heterogeneity enters only through the exponentially decaying term.

	\newcommand{\rowlabel}[2]{%
		{\renewcommand{\arraystretch}{1.05}%
			\begin{tabular}[c]{@{}l@{}}#1\\{\scriptsize #2}\end{tabular}}}
	
	\begin{table*}[t]
		\centering
		\caption{Stochastic and heterogeneity contributions to strongly convex
			GT bounds under consensual initialization.}
		\label{tab:gt_rate_comparison}
		\small
		\renewcommand{\arraystretch}{1.35}
		\setlength{\tabcolsep}{3pt}
		\begin{tabular*}{\linewidth}{@{\extracolsep{\fill}}lccc@{}}
			\toprule
			Result and spectral setting
			& Stochastic terms
			& 	\shortstack{Coefficient of\\
				$\varsigma_\star^2 e^{-\mu\bar\alpha K/2}$}
			& 
			\shortstack{Network scaling of\\
				an admissible cap $\bar\alpha$}
			\\
			\midrule
			\multicolumn{4}{@{}l}{\textbf{Existing benchmarks}}\\[-1mm]
			
			\rowlabel{\cite{alghunaim2023enhanced}}
			{symmetric $W\succeq0$}
			& $\displaystyle
			\frac{1}{NK}+\frac{1}{pK^2}+\frac{1}{p^2NK^3}$
			& $\displaystyle
			\frac{\bar\alpha^2}{p^2}$
			& $p$
			\\
			\addlinespace[2pt]
			
			\rowlabel{\cite{koloskova2021improved}}
			{symmetric $W$}
			& $\displaystyle
			\frac{1}{NK}+\frac{1}{\Delta_\lambda\nu^2K^2}$
			& $\displaystyle
			\frac{\bar\alpha}{\Delta_\lambda^2}$
			& $\displaystyle
			\frac{\Delta_\lambda\nu}{\ell_\lambda}$
			\\
			\midrule
			\multicolumn{4}{@{}l}{\textbf{This work}}\\[-1mm]
			
			\rowlabel{Full spectrum}{Cor.~\ref{cor:gt_nonsym_linear_speedup}}
			& $\displaystyle
			\frac{1}{NK}
			+\frac{\kappa_W^2}{NK^2}\sum_{i=2}^N
			\frac{|\lambda_i|^2|1-\lambda_i|^2}
			{(1-|\lambda_i|^2)^3}$
			& $\displaystyle
			\frac{\kappa_W^2\bar\alpha^3}{\Delta_\lambda^3}$
			& $\displaystyle
			\frac{\Delta_\lambda}{\kappa_W\beta_W}$
			\\
			\addlinespace[2pt]
			
			\rowlabel{Two-region}{Cor.~\ref{cor:gt_nonsym_two_region_rate}}
			& $\displaystyle
			\frac{1}{NK}+\frac{\kappa_W^2}{K^2}
			\left(\frac{\theta_1c^2}{p}+\frac{\theta_2}{\nu^3}\right)$
			& $\displaystyle
			\frac{\kappa_W^2\bar\alpha^3}{\Delta_\lambda^3}$
			& $\displaystyle
			\frac{\min\{p/c,\,\nu^2\}}{\kappa_W}$
			\\
			\addlinespace[2pt]
			
			\rowlabel{Two-region, symmetric}{Cor.~\ref{cor:gt_symmetric_two_region}}
			& $\displaystyle
			\frac{1}{NK}+\frac{1}{K^2}
			\left(\frac{\theta_1}{p}+\frac{\theta_2}{\nu^3}\right)$
			& $\displaystyle
			\frac{\bar\alpha^3}{\Delta_\lambda^3}$
			& $\min\{p,\,\nu^2\}$
			\\
			\midrule
			\multicolumn{4}{@{}l}{\textbf{This work: lower bound}}\\[-1mm]
			
			\rowlabel{Family $W_{p,\nu}$}{Cor.~\ref{cor:jlb_tuned}}
			& \multicolumn{3}{l@{}}{%
				$\displaystyle
				\bar e_K\gtrsim
				\min\!\left\{\bar e_0,\,
				\frac{\sigma^2}{NK}
				+\frac{\sigma^2}{K^2}
				\left(\frac1p+\frac1{\nu^3}\right)\right\}
				+\frac{\alpha^4\varsigma_\star^2}{p^4}
				\rho_{\mathrm{lb}}^{\,2(K-T_p-1)}$}
			\\
			\bottomrule
		\end{tabular*}
		
		\par\smallskip
		\begin{minipage}{\linewidth}
			\footnotesize
			\textit{Notes.}
			The stochastic entries in the upper-bound rows omit the common
			factor $\sigma^2$.
			The third column gives the heterogeneity-dependent part of the
			transient, suppressing logarithmic factors and constants depending
			only on the optimization parameters; other initialization terms
			are omitted. The coefficients for the prior bounds follow from
			\cite[proof of Theorem~2]{alghunaim2023enhanced} and
			\cite[App.~B, Lemmas~23--24]{koloskova2021improved}. The last column gives
			the order of an admissible cap.
			For \cite{koloskova2021improved},
			$\ell_\lambda=1+\ln(1/\Delta_\lambda)$ and
			$K\gtrsim\ell_\lambda/\Delta_\lambda$.
			
			The full-spectrum and general two-region rows allow diagonalizable
			$W$ with complex eigenvalues. Here
			$\Delta_\lambda=1-\max_{i\geq2}|\lambda_i|$,
			$\beta_W=\max_{i\geq2}|1-\lambda_i|/(1-|\lambda_i|)$,
			and $\kappa_W$ is the eigenvector condition number.
			For normal $W$, the modal sum in \eqref{gt_nonsym_avg_gain} is exact.
			The regions are defined in \eqref{gt_nonsym_two_region_conditions},
			with $\theta_j=|\mathcal I_j|/N$. For symmetric $W$, they contain
			the nonnegative and negative nontrivial eigenvalues, respectively;
			$p=1-\max\{\lambda_2(W),0\}$,
			$\nu=1+\min\{\lambda_N(W),0\}$, and
			$\Delta_\lambda=\min\{p,\nu\}$.
			
			The lower row uses $W_{p,\nu}$ with $4\mid N$, $N\geq4$,
			$L=3\mu$, $0<p\leq\nu\leq1/4$, and $K\geq2T_p$.
			Here $\alpha$ is the actual stepsize satisfying
			\eqref{eq:jlb_stepsize}, and $\rho_{\mathrm{lb}}=1-2\mu\alpha$.
			On this family, $\theta_1=1/2-1/N$ and $\theta_2=1/2$.
		\end{minipage}
	\end{table*}

	\subsection{Comparison for symmetric mixing matrices}
	\label{sec:comparison_dictionary}
	
	For symmetric $W$, with $p$ and $\nu$ as in Section~\ref{sec:contributions}
	and $\Delta_\lambda=\min\{p,\nu\}$, the graph parameters	$p_{\mathrm{ref}},c_{\mathrm{ref}}$ of
	\cite{koloskova2021improved} satisfy $p_{\mathrm{ref}}=\Delta_\lambda(2-\Delta_\lambda)$,
	$c_{\mathrm{ref}}=\nu(2-\nu)$.
	Thus its $K^{-2}$ network factor
	$(p_{\mathrm{ref}}c_{\mathrm{ref}}^2)^{-1}$ is equivalent, up to
	absolute constants, to $(\Delta_\lambda\nu^2)^{-1}$, whereas our
	averaged symmetric factor is
	$\theta_1/p+\theta_2/\nu^3\leq p^{-1}+\nu^{-3}$.
	For $W\succeq0$, $\nu=1$, and our average-iterate bound recovers
	the linear dependence $p^{-1}$ of \cite{alghunaim2023enhanced}.
	
	The work \cite{alghunaim2023enhanced} analyzes the ATC variant that also
	mixes the stochastic-gradient difference. We compare the spectral-gap
	dependence of the bounds for the respective GT variants, rather than two
	analyses of an identical recursion. Table~\ref{tab:gt_rate_comparison} compares the rates; its last row shows
	that the additive factor $p^{-1}+\nu^{-3}$ is matched, up to constants and
	logarithmic factors, on $W_{p,\nu}$.

	\begin{remark}[\bfseries\small Deterministic case]
		\label{rem:deterministic}
		\rm
		Table~\ref{tab:gt_rate_comparison} displays the
		heterogeneity-dependent initialization prefactors and omits
		the other initialization terms. For $\sigma=0$ and
		consensual initialization, Theorem~\ref{thm:gt_nonsym_refined} and,
		on $W_{p,\nu}$ under its conditions, Theorem~\ref{thm:jlb_joint} give
		\[
		\begin{aligned}
			\bar e_K&\leq\Bigl(2\bar e_0+
			\frac{72L\kappa_W^2\varsigma_\star^2\alpha^3}{\Delta_\lambda^3}\Bigr)
			e^{-\mu\alpha K/2},\\
			\bar e_K&\gtrsim\bar e_0\rho_{\mathrm{lb}}^{2K}
			+\frac{\alpha^4\mu^2\varsigma_\star^2}{p^4}
			\rho_{\mathrm{lb}}^{2(K-T_p-1)}.
		\end{aligned}
		\]
		These heterogeneity prefactors are not matched.
		On $W_{p,\nu}$, $\mathcal G_1(\omega)=\Theta(p^{-1}+\nu^{-2})$
		and $\bar\alpha=\Theta(\min\{p,\nu^2\}/L)$. When
		$p=\cO(\nu^2)$, an admissible choice $\alpha=\Theta(p/L)$
		cancels their explicit $p$-powers. The decay remains
		network-dependent: with $\alpha=\bar\alpha$, accuracy $\epsilon$
		is guaranteed within $\widetilde\cO(\max\{\kappa,\Delta_\lambda^{-1},
		\kappa\mathcal G_1(\omega)\}\ln(1/\epsilon))$ iterations
		(Remark~\ref{rem:gt_nonsym_transient}), which is the centralized order
		$\widetilde\cO(\kappa\ln(1/\epsilon))$ when
		$\Delta_\lambda^{-1},\mathcal G_1(\omega)=\cO(1)$, and larger by the
		factor $p^{-1}+\nu^{-2}$ on $W_{p,\nu}$.
	\end{remark}

	\begin{lemma}[\bfseries\small Additive versus multiplicative dependence]
		\label{lemma:cmp_ratio}
		For $j\in\{1,2\}$, let
		\begin{equation}\label{cmp_Rj}
			\mathcal N_j
			\define
			\frac1p+\frac1{\nu^{j+1}},
			\qquad
			\mathcal N_j^\star
			\define
			\frac{1}{\Delta_\lambda\nu^j}.
		\end{equation}
		Then
		\begin{equation}\label{cmp_Rj_bound}
			\frac{\mathcal N_j}{\mathcal N_j^\star}
			=
			\frac{\Delta_\lambda\nu^j}{p}
			+\frac{\Delta_\lambda}{\nu}
			\leq2.
		\end{equation}
		If $p\leq\nu$, this ratio equals
		$\nu^j+p/\nu$.
	\end{lemma}
	
	\begin{proof}
		The result follows from
		$\Delta_\lambda\leq p$, $\Delta_\lambda\leq\nu$, and
		$\Delta_\lambda=p$ when $p\leq\nu$.
	\end{proof}
	
	For $j=2$, Lemma~\ref{lemma:cmp_ratio} compares our unweighted
	stochastic envelope with the prior bound; $j=1$ compares the
	geometric time scales up to logarithmic factors. For example, if
	$p=\nu^3$, the prior stochastic factor is $\nu^{-5}$, our envelope
	is $2\nu^{-3}$, and our averaged factor is
	$(\theta_1+\theta_2)\nu^{-3}=(1-1/N)\nu^{-3}$.
	The geometric time scale improves from order $\nu^{-4}$ to
	$\nu^{-3}$, up to logarithmic factors.

	\subsection{Beyond symmetric mixing}
	\label{sec:comparison_nonsymmetric}
	Both \cite{alghunaim2023enhanced,koloskova2021improved} consider
	symmetric mixing, whereas our analysis also covers diagonalizable
	nonsymmetric matrices with real or complex spectra. Nonnormality enters
	through $\kappa_W$, which equals one for normal matrices. For example,
	the lazy directed ring of Example~\ref{ex:gt_nonsym_directed_ring}
	satisfies
	$\mathcal G_2=\cO(\Delta_\lambda^{-2})$, showing that the exact
	propagation analysis can capture spectral regimes between the linear
	and worst-case cubic inverse-gap dependences.
	
	\section{Numerical validation}
	\label{sec:numerics}
	
	We provide three numerical experiments that isolate the spectral effects
	studied in the paper. The first validates the cubic lower bound, while
	the other two evaluate the propagation gain $\mathcal G_2$ in
	\eqref{gt_nonsym_network_gains}, which bounds the network-dependent
	stochastic term. For the two normal-matrix experiments, we evaluate the
	propagation gain by the truncated modal sum
	$
	\mathcal G_{2,L_{\max}}
	=
	\sum_{\ell=0}^{L_{\max}}
	\max_{i\geq2}|m_\ell(\lambda_i)|^2$,
	$L_{\max}=\left\lceil20/\Delta_\lambda\right\rceil$.
	The plotted numerical gains use this finite-sum evaluation.
	For the cycle and directed-ring grids below, doubling the cutoff
	to $\lceil40/\Delta_\lambda\rceil$ changes every computed gain
	by less than $10^{-12}$ relative and leaves the fitted slopes
	unchanged to the reported precision.
	
	\subsection{Cubic lower-bound scaling}
	
\begin{figure}[h!]
		\centering
		\includegraphics[width=0.68\linewidth]{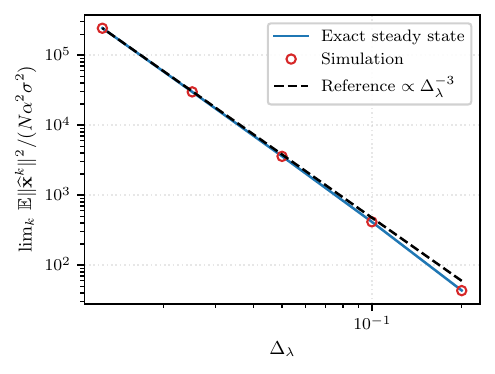}
		\caption{Normalized steady-state disagreement for the two-agent
			negative-mode construction. Sample averages over $4000$ noise
			realizations agree with the exact steady state and follow the predicted
			$\Delta_\lambda^{-3}$ scaling.}
		\label{fig:gt_cubic_lower_bound}
	\end{figure}

	We first consider the two-agent quadratic construction of
	Theorem~\ref{thm:cubic_network_lower_bound}. We set
	$\mu=\sigma^2=1$ and, for
	$\Delta_\lambda\in\{0.2,0.1,0.05,0.025,0.0125\}$, choose
	$\alpha=\Delta_\lambda^2/64$, which satisfies
	\eqref{eq:lb_stepsize}. We simulate the two-dimensional modal
	recursion \eqref{eq:lb_scalar_recursion}, initialized at $u_0=w_0=0$,
	with $4000$ independent Gaussian noise sequences. A single NumPy
	generator, initialized by \texttt{default\_rng(7)}, is used for
	the gaps in the listed order. We estimate
	$\Ex\|\widehat{\x}^k\|^2/(N\alpha^2\sigma^2)$ by averaging over the
	realizations and the last $\lceil20/\Delta_\lambda\rceil$ iterations,
	after a burn-in of $\lceil40/\Delta_\lambda\rceil$ iterations.
	Estimated standard errors, computed across the independent
	per-trajectory time averages, range from $0.72\%$ to $0.78\%$
	of the corresponding estimates.
	Figure~\ref{fig:gt_cubic_lower_bound} shows that the sample averages
	differ from the exact values in \eqref{eq:lb_exact_steady_state}
	by less than $0.6\%$. A log--log fit
	gives exponent $-3.11$ for both the simulated and exact values,
	consistent with the predicted cubic dependence. The small deviation
	from $-3$ over this finite range comes from the gap-dependent prefactor
	in the exact expression.
	
	\subsection{Positive--negative spectral separation}
	
\begin{figure}[h!]
		\centering
		\includegraphics[width=0.68\linewidth]{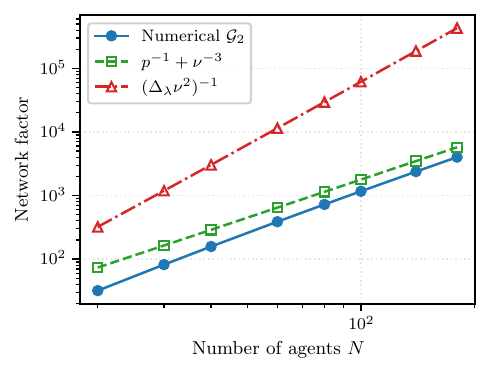}
		\caption{Propagation gain for the symmetric cycle family with
			$\epsilon_N=N^{-2/3}$. The numerically evaluated $\mathcal G_2$ closely follows the
			additive factor $p^{-1}+\nu^{-3}$, while
			$(\Delta_\lambda\nu^2)^{-1}$ grows substantially faster.}
		\label{fig:gt_two_region}
	\end{figure}

	We next illustrate the benefit of separating the two ends of a symmetric
	spectrum. Let $P$ denote the cyclic shift on an even $N$-node cycle,
	with $N\geq4$, and set
	\[
	A=\frac12(P+P\tran),
	\qquad
	W_N=\epsilon_N I+(1-\epsilon_N)A,
	\qquad
	\epsilon_N=N^{-2/3}.
	\]
	Then $W_N$ is symmetric, nonnegative, and doubly stochastic, with
	\[
	p=\Theta(N^{-2}),
	\qquad
	\nu=2\epsilon_N=\Theta(N^{-2/3}),
	\qquad
	\Delta_\lambda=p.
	\]
	Hence
	$p^{-1}+\nu^{-3}=\Theta(N^2)$, whereas the multiplicative factor of
	\cite{koloskova2021improved} satisfies
	$(\Delta_\lambda\nu^2)^{-1}=\Theta(N^{10/3})$.
	For $N\in\{20,30,40,60,80,100,140,180\}$, we evaluate
	$\mathcal G_{2,L_{\max}}$ from the eigenvalues using \eqref{ml_poly}. Figure~\ref{fig:gt_two_region} shows log--log slopes
	$2.20$ for the numerically evaluated gain and $1.98$
	for the additive factor, with their ratio between $0.43$ and $0.70$. In contrast, the multiplicative proxy
	has slope $3.28$ and grows much faster. The size-dependent
	laziness is used only to produce a controlled family in which the two
	spectral ends deteriorate at different rates.
	
	\subsection{Nonsymmetric normal mixing}
	
\begin{figure}[h!]
		\centering
		\includegraphics[width=0.68\linewidth]{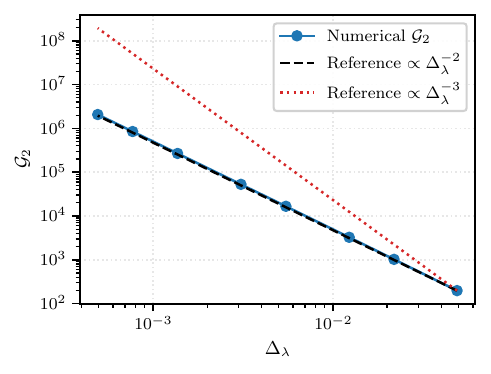}
		\caption{Propagation gain $\mathcal G_2$ for the lazy directed ring
			$W=(I+P)/2$. The observed scaling follows
			$\Delta_\lambda^{-2}$ and is substantially sharper than the
			modulus-only cubic dependence.}
		\label{fig:gt_directed_ring_gain}
	\end{figure}

	Finally, we consider the lazy directed ring
	$W=\frac12(I+P)$ of Example~\ref{ex:gt_nonsym_directed_ring} with $N\in\{10,15,20,30,\allowbreak40,60,80,100\}$. This matrix is nonsymmetric
	but normal, so $\kappa_W=1$. For each network size, we evaluate
	$\mathcal G_{2,L_{\max}}$ directly from its complex eigenvalues. As shown in
	Figure~\ref{fig:gt_directed_ring_gain}, a log--log fit gives exponent
	$-2.01$, showing a $\Delta_\lambda^{-2}$ dependence. For the dominant mode
	$z_1=\cos(\pi/N)e^{\imath\pi/N}$,
	Lemma~\ref{lem:gt_nonsym_modal_energy} gives
	\[
	\Delta_\lambda^2
	\sum_{\ell=0}^{\infty}|m_\ell(z_1)|^2
	=
	\frac{2\cos^2(\pi/N)}{(1+\cos(\pi/N))^2}
	\longrightarrow\frac12.
	\]
	Since $\|\mathbf M_\ell\|_2\geq|m_\ell(z_1)|$ for normal $W$, this
	single-mode value bounds $\Delta_\lambda^2\mathcal G_2$ from below,
	consistent with Example~\ref{ex:gt_nonsym_directed_ring}.
	
	\section{Conclusion}
	\label{sec:conclusion}

	This work shows that the noise amplification of stochastic GT is
	governed by the location of the eigenvalues, not by the spectral gap
	alone.			As an eigenvalue $z$ approaches the unit circle away from $1$, its
	weight in the GT noise bound grows as $(1-|z|)^{-3}$. For normal
	mixing, the bound averages these weights. On the directed ring with
	fixed laziness, this improves the network factor from
	$\Delta_\lambda^{-2}$ to $\Delta_\lambda^{-3/2}$. For symmetric
	mixing, the additive factor $p^{-1}+\nu^{-3}$ is sharp up to
	logarithmic factors on a family of arbitrarily large networks.
	Exact diffusion \cite{yuan2020influence} relies on a
	positive semidefinite combination matrix, and for this case its
	comparison with GT follows from existing bounds
	\cite{alghunaim2021unified}; modifying the algorithm itself so that it
	applies to mixing matrices with complex spectra is an interesting
	direction.

	\newpage
	\appendix
	\section{Shared estimates and stepsize tuning}
	\label{app:shared_estimates}
	This appendix collects two classical estimates and a stepsize-tuning
	lemma used in Appendices~\ref{app:refined_convergence}
	and~\ref{sec:gt_two_recursion_appendix}. The estimates require only
	double stochasticity of $W$ and not the norm condition of
	Assumption~\ref{assump:gt_semi_network}, so that the refined analysis
	of Section~\ref{sec:gt_nonsymmetric_refined} does not depend on it.
	
	\begin{lemma}[\bfseries\small Average recursion and gradient perturbation]
		\label{lemma:gt_shared_estimates}
		Let Assumptions~\ref{assump:gt_semi_functions} and
		\ref{assump:gt_semi_noise} hold, let $W$ be doubly stochastic, let
		$0<\alpha\leq1/(4L)$, and write $\widetilde x^k\define\bar x^k-x^\star$.
		Then, for every $k\geq0$,
		\begin{align}
			\Ex\!\left[\|\widetilde x^{k+1}\|^2\,\middle|\,\cF^k\right]
			&\leq(1-\mu\alpha)\|\widetilde x^k\|^2
			-2\alpha(1-2\alpha L)\bigl(f(\bar x^k)-f(x^\star)\bigr)
			\nonumber\\&\quad+\frac{\alpha L}{N}(1+2\alpha L)\|\widehat{\x}^k\|^2
			+\frac{\alpha^2\sigma^2}{N},
			\label{gt_semi_average_before_step}
		\end{align}
		and, with $\mathbf u^k$ as in \eqref{gt_semi_u_definition},
		\begin{align}
			\|\mathbf u^k\|^2
			&\leq2L^2\|\widehat{\x}^k\|^2
			+2\sum_{i=1}^{N}\|\grad f_i(\bar x^k)-\grad f_i(x^\star)\|^2,
			\label{eq:appendix_u_split}\\
			\sum_{i=1}^{N}\|\grad f_i(\bar x^k)-\grad f_i(x^\star)\|^2
			&\leq2LN\bigl(f(\bar x^k)-f(x^\star)\bigr).
			\label{eq:appendix_cocoercive}
		\end{align}
	\end{lemma}
	\begin{proof}
		Let $\widetilde x^k\define\bar x^k-x^\star$. From
		\eqref{gt_semi_centroid}, conditional unbiasedness gives
		\begin{align}
			\Ex\!\left[
			\|\widetilde x^{k+1}\|^2
			\,\middle|\,
			\cF^k
			\right]
			&=
			\|\widetilde x^k
			-\alpha\overline{\grad\f}(\x^k)\|^2
			+\alpha^2
			\Ex\!\left[
			\|\bar v^k\|^2
			\,\middle|\,
			\cF^k
			\right].
			\label{eq:appendix_average_expansion}
		\end{align}
		By conditional independence across agents and
		\eqref{gt_semi_variance},
		$\Ex[\|\bar v^k\|^2\mid\cF^k]\leq\sigma^2/N$.
		For each agent, strong convexity and smoothness give
		\[
		\begin{aligned}
			\langle\widetilde x^k,\grad f_i(x_i^k)\rangle
			&\geq
			f_i(\bar x^k)-f_i(x^\star)
			+\frac{\mu}{2}\|x_i^k-x^\star\|^2
			-\frac{L}{2}\|x_i^k-\bar x^k\|^2.
		\end{aligned}
		\]
		Summing over $i$, dividing by $N$, using
		$\sum_i\|x_i^k-x^\star\|^2
		=\|\widehat{\x}^k\|^2+N\|\widetilde x^k\|^2$,
		and dropping the nonnegative
		$\mu\|\widehat{\x}^k\|^2/(2N)$ term yields
		\begin{equation}\label{eq:appendix_average_inner}
			\left\langle
			\widetilde x^k,\overline{\grad\f}(\x^k)
			\right\rangle
			\geq
			f(\bar x^k)-f(x^\star)
			+\frac{\mu}{2}\|\widetilde x^k\|^2
			-\frac{L}{2N}\|\widehat{\x}^k\|^2.
		\end{equation}
		Write
		\[
		\overline{\grad\f}(\x^k)
		=
		\grad f(\bar x^k)+r^k,
		\qquad
		r^k
		\define
		\frac1N\sum_i
		\bigl(\grad f_i(x_i^k)-\grad f_i(\bar x^k)\bigr).
		\]
		By smoothness and Jensen's inequality,
		$\|r^k\|^2\leq L^2\|\widehat{\x}^k\|^2/N$. Since
		$\|\grad f(\bar x^k)\|^2
		\leq2L(f(\bar x^k)-f(x^\star))$,
		\begin{equation}\label{eq:appendix_average_gradient}
			\|\overline{\grad\f}(\x^k)\|^2
			\leq
			4L\bigl(f(\bar x^k)-f(x^\star)\bigr)
			+\frac{2L^2}{N}\|\widehat{\x}^k\|^2.
		\end{equation}
		Substituting
		\eqref{eq:appendix_average_inner},
		\eqref{eq:appendix_average_gradient}, and
		$\Ex[\|\bar v^k\|^2\mid\cF^k]\leq\sigma^2/N$
		into \eqref{eq:appendix_average_expansion} yields
		\eqref{gt_semi_average_before_step}.
		For $\mathbf u^k$ in \eqref{gt_semi_u_definition}, splitting at $\bar x^k$ and applying smoothness gives \eqref{eq:appendix_u_split}.
		For every convex $L$-smooth function,
		\begin{equation}
			\|\grad f_i(y)-\grad f_i(z)\|^2
			\leq
			2L\Bigl(
			f_i(y)-f_i(z)
			-\langle\grad f_i(z),y-z\rangle
			\Bigr).
		\end{equation}
		Taking $y=\bar x^k$, $z=x^\star$, summing over $i$, and using
		$\sum_{i=1}^{N}\grad f_i(x^\star)=\zero$, we obtain \eqref{eq:appendix_cocoercive}.
	\end{proof}
	
	The following elementary stepsize lemma is used for
	Theorem~\ref{thm:gt_baseline_convergence} and
	Corollary~\ref{cor:gt_nonsym_linear_speedup}.
	
	\begin{lemma}[\bfseries\small Logarithmic stepsize choice]
		\label{lem:log_stepsize}
		Let $K\geq1$, $\mu>0$, and $\bar\alpha>0$. Suppose that, for every
		$\alpha\in(0,\bar\alpha]$, a nonnegative quantity $\mathcal R_K(\alpha)$
		satisfies
		\begin{equation}
			\mathcal R_K(\alpha)
			\leq
			\exp\!\left(
			-\frac{\mu\alpha K}{2}
			\right)
			\left(
			c_0+\widetilde c_2\alpha^3
			\right)
			+c_1\alpha+c_2\alpha^2,
			\label{eq:tuning_recursion}
		\end{equation}
		where $c_0,\widetilde c_2,c_2\geq0$, $c_1>0$, and
		$A\define c_0+\widetilde c_2\bar\alpha^3\geq 0$. Define
		\begin{equation}
			T_K
			\define
			\max\left\{
			\exp(1),
			\frac{\mu A K}{c_1}
			\right\},
		\end{equation}
		and choose
		\begin{equation}
			\alpha
			=
			\min\left\{
			\frac{2\ln T_K}{\mu K},
			\bar\alpha
			\right\}.
			\label{eq:tuning_alpha}
		\end{equation}
		Then
		\begin{align}
			\mathcal R_K(\alpha)
			&\leq
			A
			\exp\!\left(
			-\frac{\mu\bar\alpha K}{2}
			\right)
			+\frac{(1+2\ln T_K)c_1}{\mu K}
			+\frac{4c_2\ln^2T_K}{\mu^2K^2}.
			\label{eq:tuning_bound_exact}
		\end{align}
		Consequently,
		\begin{equation}
			\mathcal R_K
			=
			\widetilde\cO\!\left(
			\frac{c_1}{\mu K}
			+
			\frac{c_2}{\mu^2K^2}
			\right)
			+
			A
			\exp\!\left(
			-\frac{\mu\bar\alpha K}{2}
			\right).
			\label{eq:tuning_bound}
		\end{equation}
	\end{lemma}
	
	\begin{proof}
		If
		$\alpha=2\ln T_K/(\mu K)\leq\bar\alpha$, then
		$e^{-\mu\alpha K/2}=T_K^{-1}$ and, by the definition of $T_K$,
		\[
		Ae^{-\mu\alpha K/2}
		\leq
		\frac{c_1}{\mu K}.
		\]
		Moreover,
		\[
		c_1\alpha
		=
		\frac{2c_1\ln T_K}{\mu K},
		\qquad
		c_2\alpha^2
		=
		\frac{4c_2\ln^2T_K}{\mu^2K^2}.
		\]
		This gives \eqref{eq:tuning_bound_exact}.
		
		If instead $\alpha=\bar\alpha$, then
		$\bar\alpha\leq2\ln T_K/(\mu K)$, so
		\[
		c_1\bar\alpha
		\leq
		\frac{2c_1\ln T_K}{\mu K},
		\qquad
		c_2\bar\alpha^2
		\leq
		\frac{4c_2\ln^2T_K}{\mu^2K^2},
		\]
		while the transient term is
		$A\exp(-\mu\bar\alpha K/2)$.
		Thus \eqref{eq:tuning_bound_exact} holds in both cases.
		The bound \eqref{eq:tuning_bound} follows by suppressing logarithmic
		factors.
	\end{proof}

	\section{Proof of Theorem~\ref{thm:cubic_network_lower_bound}}
	\label{app:cubic_network_lower_bound}
	\begin{proof}[Proof of Theorem \ref{thm:cubic_network_lower_bound}]
		By Lemma~\ref{lem:negative_mode_matrix}, the disagreement subspace is
		one-dimensional and is spanned by $v$. Define
		\begin{equation}\label{eq:lb_coordinates}
			u_k
			\define
			v\tran\x^k,
			\qquad
			w_k
			\define
			\alpha v\tran\y^k,
			\qquad
			\varepsilon_k
			\define
			v\tran\bxi^k
			=
			\frac{\xi_1^k-\xi_2^k}{\sqrt2},
		\end{equation}
		and let
		\begin{equation}\label{eq:lb_beta}
			\beta\define\alpha\mu.
		\end{equation}
		By \eqref{eq:lb_noise}, $\Ex[\varepsilon_k]=0$,
		$\Ex[|\varepsilon_k|^2]=\sigma^2$,
		and $\{\varepsilon_k\}$ is independent across iterations.	
		Multiplying \eqref{gt_atc_semi2} by $v\tran$ on the left and using
		$v\tran W_\lambda=-\lambda v\tran$ gives
		\begin{subequations}\label{eq:lb_scalar_recursion}
			\begin{align}
				u_{k+1}
				&=
				-\lambda(1-\beta)u_k
				+\lambda w_k
				+\alpha\lambda\varepsilon_k,
				\\
				w_{k+1}
				&=
				-\beta(1+\lambda)u_k
				-\lambda w_k
				-\alpha(1+\lambda)\varepsilon_k.
			\end{align}
		\end{subequations}
		Equivalently,
		\begin{equation}\label{eq:lb_state_recursion}
			\theta_{k+1}
			=
			A\theta_k+B\varepsilon_k,
			\qquad
			\theta_k
			\define
			\col\{u_k,w_k\},
		\end{equation}
		where
		\begin{equation}\label{eq:lb_AB}
			A
			\define
			\begin{bmatrix}
				-\lambda(1-\beta) & \lambda\\
				-\beta(1+\lambda) & -\lambda
			\end{bmatrix},
			\qquad
			B
			\define
			\alpha
			\begin{bmatrix}
				\lambda\\
				-(1+\lambda)
			\end{bmatrix}.
		\end{equation}
		Since $I-\Pi=vv\tran$, the disagreement lies entirely along $v$:
		\[
		\widehat{\x}^k
		=
		(I-\Pi)\x^k
		=
		vv\tran\x^k
		=
		u_k v.
		\]
		Hence, using $\|v\|=1$ and $N=2$,
		\begin{equation}\label{eq:lb_disagreement_mode}
			\frac1N\|\widehat{\x}^k\|^2
			=
			\frac{|u_k|^2}{2}.
		\end{equation}
		The characteristic polynomial of $A$ is
		\begin{align}
			\det(zI-A)
			&=
			z^2-\operatorname{tr}(A)z+\det(A) =z^2+\lambda(2-\beta)z+\lambda^2+\lambda\beta.
		\end{align}
		To verify stability, we use the Jury criterion for a second-order polynomial
		\cite{jury1964}. In particular, a monic quadratic with real coefficients
		$
		z^2+bz+c
		$
		has both roots strictly inside the unit disk if and only if
		\[
		1+b+c>0,
		\qquad
		1-b+c>0,
		\qquad
		1-c>0.
		\]
		For the characteristic polynomial above,
		$b=\lambda(2-\beta)$ and $c=\lambda^2+\lambda\beta$, so these three
		conditions become
		\begin{subequations}
			\begin{align}
				1+\lambda(2-\beta)+\lambda^2+\lambda\beta
				&=(1+\lambda)^2>0,
				\\
				1-\lambda(2-\beta)+\lambda^2+\lambda\beta
				&=\Delta_\lambda^2+2\lambda\beta>0,
				\\
				1-\lambda^2-\lambda\beta
				&=\Delta_\lambda(1+\lambda)-\lambda\beta.
			\end{align}
		\end{subequations}
		The first two Jury conditions hold automatically. The third requires
		\[
		1-\lambda^2-\lambda\beta
		=
		\Delta_\lambda(1+\lambda)-\lambda\beta
		>0,
		\]
		which is equivalent to
		$
		\beta
		<
		\frac{\Delta_\lambda(1+\lambda)}{\lambda}.
		$
		Thus \eqref{eq:lb_stability_stepsize} implies $\rho(A)<1$. Since $\rho(A)<1$, the linear recursion \eqref{eq:lb_state_recursion} is
		Schur stable. Moreover, $\varepsilon_k$ is zero-mean white noise with variance
		$\sigma^2$ and is independent of $\theta_k$. Thus, defining
		$
		P_k\define\Ex[\theta_k\theta_k\tran],
		$
		we obtain
		\[
		P_{k+1}
		=
		AP_kA\tran+\sigma^2BB\tran.
		\]
		For a Schur-stable matrix $A$, this second-moment recursion converges to the
		unique solution of the discrete-time Lyapunov equation \cite{chen1999linear}:
		\begin{equation}\label{eq:lb_lyapunov}
			P
			=
			APA\tran+\sigma^2BB\tran.
		\end{equation}
		Equivalently, 	$P
		=
		\sum_{j=0}^{\infty}
		A^j(\sigma^2BB\tran)(A\tran)^j$. Solving the discrete-time Lyapunov equation gives
		\begin{equation}\label{eq:lb_covariance}
			P
			=
			\frac{\alpha^2\sigma^2}
			{
				(1-\lambda^2-\lambda\beta)
				(\Delta_\lambda^2+2\lambda\beta)
			}
			\begin{bmatrix}
				2\lambda^2
				&
				-\lambda(1-\lambda^2)
				\\
				-\lambda(1-\lambda^2)
				&
				(1+\lambda)^2
				(\Delta_\lambda^2+\lambda\beta)
			\end{bmatrix},
		\end{equation}
		which can be verified by direct substitution into
		\eqref{eq:lb_lyapunov}. In particular,
		\[
		\lim_{k\rightarrow\infty}\Ex[|u_k|^2]
		=
		P_{11}
		=
		\frac{
			2\alpha^2\lambda^2\sigma^2
		}{
			(1-\lambda^2-\lambda\beta)
			(\Delta_\lambda^2+2\lambda\beta)
		}.
		\]
		Combining this identity with \eqref{eq:lb_disagreement_mode} and
		$\beta=\alpha\mu$ yields \eqref{eq:lb_exact_steady_state}. 
		Under \eqref{eq:lb_stepsize}, and since $\lambda\in[1/2,1)$ and
		$\beta\leq\Delta_\lambda^2/16$,
		\begin{align}
			1-\lambda^2-\lambda\beta
			&\leq
			2\Delta_\lambda,
			\\
			\Delta_\lambda^2+2\lambda\beta
			&\leq
			\frac98\Delta_\lambda^2,
		\end{align}
		while $\lambda^2\geq1/4$. Substituting these inequalities into
		\eqref{eq:lb_exact_steady_state} gives
		\[
		\lim_{k\rightarrow\infty}
		\frac1N\Ex\|\widehat{\x}^k\|^2
		\geq
		\frac{\alpha^2\sigma^2}{9\Delta_\lambda^3}.
		\]
		For the opposite direction, since $1+\lambda\geq\frac32$ and
		$\lambda\beta\leq\Delta_\lambda^2/16\leq\Delta_\lambda$,
		\[
		1-\lambda^2-\lambda\beta
		\geq
		\frac{\Delta_\lambda}{2},
		\qquad
		\Delta_\lambda^2+2\lambda\beta
		\geq
		\Delta_\lambda^2,
		\]
		and $\lambda^2\leq1$, which yields
		$
		\lim_{k\rightarrow\infty}
		\frac1N\Ex\|\widehat{\x}^k\|^2
		\leq
		\frac{2\alpha^2\sigma^2}{\Delta_\lambda^3}.
		$
		This proves \eqref{eq:lb_cubic_bound}.
	\end{proof}
	
	\section{Proofs of the spectral propagation bounds}
	
	\subsection{Proof of Lemma \ref{lem:gt_nonsym_modal_energy}} \label{app:gt_nonsym_modal_energy}
	\begin{proof}[Proof of Lemma \ref{lem:gt_nonsym_modal_energy}]
		Let $r\define|z|$ and $s\define r^2$. Recall that
		$m_\ell(z)=z^\ell(\ell(1-z)-z)$. Using
		$|a-b|^2=|a|^2+|b|^2-2\operatorname{Re}(a\bar b)$ with
		$a=\ell(1-z)$ and $b=z$, where $\operatorname{Re}(\cdot)$ denotes the real part and
		$\bar z$ the complex conjugate, we obtain
		\[
		|\ell(1-z)-z|^2
		=
		\ell^2|1-z|^2+r^2
		-2\ell\operatorname{Re}\!\left((1-z)\bar z\right).
		\]
		Since $
		\operatorname{Re}\!\left((1-z)\bar z\right)
		=
		\operatorname{Re}(z)-r^2$,
		it follows that
		\[
		|m_\ell(z)|^2
		=
		s^\ell
		\left[
		\ell^2|1-z|^2
		+2\ell(r^2-\operatorname{Re} z)
		+r^2
		\right].
		\]
		Since $|z|<1$, we have $0\leq s<1$, so the series below converge. Using
		\[
		\sum_{\ell\geq0}s^\ell
		=
		\frac1{1-s},
		\qquad
		\sum_{\ell\geq0}\ell s^\ell
		=
		\frac{s}{(1-s)^2},
		\qquad
		\sum_{\ell\geq0}\ell^2s^\ell
		=
		\frac{s(1+s)}{(1-s)^3},
		\]
		we obtain
		\begin{align*}
			\sum_{\ell=0}^{\infty}|m_\ell(z)|^2
			&=
			|1-z|^2\frac{s(1+s)}{(1-s)^3}
			+2(s-\operatorname{Re} z)\frac{s}{(1-s)^2}
			+\frac{s}{1-s}
			\\
			&=
			\frac{s}{(1-s)^3}
			\left[
			|1-z|^2(1+s)
			+2(s-\operatorname{Re} z)(1-s)
			+(1-s)^2
			\right].
		\end{align*}
		Using $|1-z|^2=1+s-2\operatorname{Re} z$, the bracketed term becomes
		\begin{align*}
			&(1+s-2\operatorname{Re} z)(1+s)
			+2(s-\operatorname{Re} z)(1-s)
			+(1-s)^2
			\\
			&=
			\bigl(1+2s+s^2-2(1+s)\operatorname{Re} z\bigr)
			+\bigl(2s-2s^2-2(1-s)\operatorname{Re} z\bigr)
			+\bigl(1-2s+s^2\bigr)
			\\
			&=
			2+2s-4\operatorname{Re} z
			=
			2|1-z|^2.
		\end{align*}
		Therefore, $
		\sum_{\ell=0}^{\infty}|m_\ell(z)|^2
		=
		\frac{2s|1-z|^2}{(1-s)^3}$,
		which gives the first equality in \eqref{gt_nonsym_modal_energy}.
		The second follows from $|1-z|=d_z\beta_z$ and
		$1-s=d_z(1+|z|)$.
	\end{proof}

	\subsection{Proof of Lemma \ref{lem:gt_nonsym_envelope}}
	\label{app:gt_nonsym_envelope}
	\begin{proof}[Proof of Lemma \ref{lem:gt_nonsym_envelope}]
		For every $i\geq2$, define
		\begin{equation}\label{gt_nonsym_Delta_i}
			\Delta_i\define1-|\lambda_i|.
		\end{equation}
		By the definitions of $\Delta_\lambda$ and $\beta_W$ in
		\eqref{gt_nonsym_spectrum_interval},
		\begin{equation}\label{gt_nonsym_spectral_bounds}
			\Delta_i\geq\Delta_\lambda,
			\qquad
			|1-\lambda_i|
			\leq
			\beta_W\Delta_i.
		\end{equation}
		Using the representation of $m_\ell$ in \eqref{ml_poly},
		\begin{align}
			|m_\ell(\lambda_i)|
			&=
			|\lambda_i|^\ell
			\bigl|\ell(1-\lambda_i)-\lambda_i\bigr|
			\nonumber\\
			&\leq
			\ell|\lambda_i|^\ell|1-\lambda_i|
			+
			|\lambda_i|^{\ell+1}
			\nonumber\\
			&\overset{\eqref{gt_nonsym_spectral_bounds}}{\leq}
			\beta_W\,\ell\Delta_i(1-\Delta_i)^\ell
			+(1-\Delta_i)^{\ell+1}.
			\label{gt_nonsym_mode_envelope}
		\end{align}	
		We next bound the first term in \eqref{gt_nonsym_mode_envelope}
		uniformly over the spectrum. By \eqref{gt_nonsym_Delta_i} and
		\eqref{gt_nonsym_spectral_bounds},
		\begin{equation}\label{gt_nonsym_ri}
			r_i\define1-\Delta_i=|\lambda_i|
			\in[0,1-\Delta_\lambda].
		\end{equation}
		Thus,
		\begin{equation}\label{gt_nonsym_Delta_to_r}
			\ell\Delta_i(1-\Delta_i)^\ell
			=
			\ell r_i^\ell(1-r_i).
		\end{equation}
		Set $h_0(r)\define0$. For each $\ell\geq1$, consider
		\begin{equation}\label{gt_nonsym_scalar_h_definition}
			h_\ell(r)
			\define
			\ell r^\ell(1-r),
			\qquad
			r\in[0,1-\Delta_\lambda].
		\end{equation}
		Its derivative is
		\begin{equation}\label{gt_nonsym_h_derivative}
			h_\ell'(r)
			=
			\ell r^{\ell-1}
			\bigl(\ell-(\ell+1)r\bigr).
		\end{equation}
		Hence $h_\ell$ is increasing for
		$r<\ell/(\ell+1)$ and decreasing for
		$r>\ell/(\ell+1)$. Its unconstrained maximizer on $[0,1]$ is therefore
		\begin{equation}\label{gt_nonsym_h_unconstrained_max}
			r_\ell^\star
			=
			\frac{\ell}{\ell+1}.
		\end{equation}	
		There are two cases. First,
		$\ell/(\ell+1)\leq1-\Delta_\lambda$ if and only if
		$\ell\leq\Delta_\lambda^{-1}-1$. In this case, the maximizer lies in
		the admissible interval and
		\begin{equation}\label{gt_nonsym_h_first_case}
			\max_{0\leq r\leq1-\Delta_\lambda}h_\ell(r)
			=
			h_\ell\!\left(\frac{\ell}{\ell+1}\right)
			=
			\left(\frac{\ell}{\ell+1}\right)^{\ell+1}.
		\end{equation}
		Second, if $\ell>\Delta_\lambda^{-1}-1$, then
		$r_\ell^\star>1-\Delta_\lambda$. Since $h_\ell$ is increasing
		throughout $[0,1-\Delta_\lambda]$, its maximum is attained at the
		right endpoint:
		\begin{equation}\label{gt_nonsym_h_second_case}
			\max_{0\leq r\leq1-\Delta_\lambda}h_\ell(r)
			=
			h_\ell(1-\Delta_\lambda)
			=
			\ell\Delta_\lambda(1-\Delta_\lambda)^\ell.
		\end{equation}
		Accordingly, for $\ell\geq1$, define
		\begin{equation}\label{gt_nonsym_psi_definition}
			\psi_\ell
			\define
			\max_{0\leq r\leq1-\Delta_\lambda}
			\ell r^\ell(1-r)
			=
			\begin{cases}
				\displaystyle
				\left(\frac{\ell}{\ell+1}\right)^{\ell+1},
				&
				\ell\leq\Delta_\lambda^{-1}-1,
				\\[3mm]
				\displaystyle
				\ell\Delta_\lambda
				(1-\Delta_\lambda)^\ell,
				&
				\ell>\Delta_\lambda^{-1}-1,
			\end{cases}
			\qquad
			\psi_0\define0.
		\end{equation}
		Since \eqref{gt_nonsym_ri} shows that
		$r_i\in[0,1-\Delta_\lambda]$, the definition
		\eqref{gt_nonsym_psi_definition} yields
		\begin{equation}\label{gt_nonsym_Delta_psi_bound}
			\ell\Delta_i(1-\Delta_i)^\ell
			\overset{\eqref{gt_nonsym_Delta_to_r}}{=}
			h_\ell(r_i)
			\overset{\eqref{gt_nonsym_psi_definition}}{\leq}
			\psi_\ell.
		\end{equation}
		Combining \eqref{gt_nonsym_mode_envelope},
		\eqref{gt_nonsym_Delta_psi_bound}, and
		\eqref{gt_nonsym_spectral_bounds} gives
		\begin{align}
			\max_{i\geq2}|m_\ell(\lambda_i)|
			&\overset{\eqref{gt_nonsym_mode_envelope},
				\eqref{gt_nonsym_Delta_psi_bound}}{\leq}
			\beta_W\psi_\ell
			+(1-\Delta_\lambda)^{\ell+1}
			\nonumber\\
			&\leq
			\beta_W\psi_\ell
			+(1-\Delta_\lambda)^\ell.
			\label{gt_nonsym_mode_envelope_uniform}
		\end{align}
		Let
		\begin{equation}\label{gt_nonsym_LDelta}
			L_\Delta
			\define
			\left\lfloor
			\Delta_\lambda^{-1}-1
			\right\rfloor.
		\end{equation}
		By \eqref{gt_nonsym_psi_definition}, the sums involving $\psi_\ell$
		split into the ranges $1\leq\ell\leq L_\Delta$ and
		$\ell\geq L_\Delta+1$. We bound these two ranges separately.
		
		For $1\leq\ell\leq L_\Delta$, the first branch of
		\eqref{gt_nonsym_psi_definition} gives
		\begin{equation}\label{gt_nonsym_psi_first_bound}
			\psi_\ell
			=
			\left(\frac{\ell}{\ell+1}\right)^{\ell+1}
			=
			\left(1-\frac1{\ell+1}\right)^{\ell+1}
			\leq
			e^{-1}.
		\end{equation}
		Since $L_\Delta\leq\Delta_\lambda^{-1}$, it follows that
		\begin{equation}\label{gt_nonsym_psi_first_square}
			\sum_{\ell=1}^{L_\Delta}\psi_\ell^2
			\overset{\eqref{gt_nonsym_psi_first_bound}}{\leq}
			L_\Delta e^{-2}
			\leq
			\frac{e^{-2}}{\Delta_\lambda}.
		\end{equation}
		For the second range, \eqref{gt_nonsym_psi_definition} gives
		$\psi_\ell=\ell\Delta_\lambda(1-\Delta_\lambda)^\ell$ for
		$\ell\geq L_\Delta+1$. Extending the sum to all $\ell\geq1$ and
		using the standard identity
		$\sum_{\ell=1}^{\infty}\ell^2x^\ell
		=x(1+x)/(1-x)^3$ for $|x|<1$, with
		$x=(1-\Delta_\lambda)^2$, gives
		\begin{align}
			\sum_{\ell=L_\Delta+1}^{\infty}\psi_\ell^2
			&\overset{\eqref{gt_nonsym_psi_definition}}{\leq}
			\Delta_\lambda^2
			\sum_{\ell=1}^{\infty}
			\ell^2(1-\Delta_\lambda)^{2\ell}
			\nonumber\\
			&=
			\Delta_\lambda^2
			\frac{
				(1-\Delta_\lambda)^2
				\left[1+(1-\Delta_\lambda)^2\right]
			}{
				\left[1-(1-\Delta_\lambda)^2\right]^3
			}.
			\label{gt_nonsym_psi_tail_square}
		\end{align}
		Since
		\begin{equation}\label{gt_nonsym_square_denominator}
			1-(1-\Delta_\lambda)^2
			=
			\Delta_\lambda(2-\Delta_\lambda)
			\geq
			\Delta_\lambda,
		\end{equation}
		and
		\begin{equation}\label{gt_nonsym_square_numerator}
			(1-\Delta_\lambda)^2
			\left[1+(1-\Delta_\lambda)^2\right]
			\leq2,
		\end{equation}
		we obtain
		\begin{equation}\label{gt_nonsym_psi_tail_square_bound}
			\sum_{\ell=L_\Delta+1}^{\infty}\psi_\ell^2
			\overset{\eqref{gt_nonsym_psi_tail_square},
				\eqref{gt_nonsym_square_denominator},
				\eqref{gt_nonsym_square_numerator}}{\leq}
			\frac{2}{\Delta_\lambda}.
		\end{equation}
		Combining this with \eqref{gt_nonsym_psi_first_square} gives
		\begin{equation}\label{gt_nonsym_psi_square_sum}
			\sum_{\ell=0}^{\infty}\psi_\ell^2
			\overset{\eqref{gt_nonsym_psi_first_square},
				\eqref{gt_nonsym_psi_tail_square_bound}}{\leq}
			\frac{2+e^{-2}}{\Delta_\lambda}
			\leq
			\frac{2.6}{\Delta_\lambda}.
		\end{equation}
		
		We next consider the discounted sum. On the first range,
		$\ell\leq L_\Delta\leq\Delta_\lambda^{-1}$. We use the generalized
		Bernoulli inequality $(1+x)^r\geq1+rx$, valid for $x\geq-1$ and
		$r\geq1$. Since $0<\Delta_\lambda\leq1$, we may take
		$x=-\Delta_\lambda/4$ and $r=1/\Delta_\lambda\geq1$. Recalling
		$\omega$ from \eqref{gt_nonsym_omega},
		\begin{equation}\label{gt_nonsym_omega_Bernoulli}
			\omega^{1/\Delta_\lambda}
			=
			\left(1-\frac{\Delta_\lambda}{4}\right)^{1/\Delta_\lambda}
			\geq
			1-\frac{1}{\Delta_\lambda}\frac{\Delta_\lambda}{4}
			=
			\frac34.
		\end{equation}
		Hence, for $1\leq\ell\leq L_\Delta$,
		\begin{equation}\label{gt_nonsym_omega_first_range}
			\omega^{-\ell}
			\leq
			\omega^{-1/\Delta_\lambda}
			\overset{\eqref{gt_nonsym_omega_Bernoulli}}{\leq}
			\frac43,
		\end{equation}
		and therefore
		\begin{equation}\label{gt_nonsym_psi_first_weighted}
			\sum_{\ell=1}^{L_\Delta}
			\omega^{-\ell}\psi_\ell
			\overset{\eqref{gt_nonsym_psi_first_bound},
				\eqref{gt_nonsym_omega_first_range}}{\leq}
			\frac{4}{3e}L_\Delta
			\leq
			\frac{4}{3e\,\Delta_\lambda}.
		\end{equation}	
		For the second range, define
		\begin{equation}\label{gt_nonsym_rho}
			\rho
			\define
			\frac{1-\Delta_\lambda}{\omega}
			=
			\frac{1-\Delta_\lambda}
			{1-\Delta_\lambda/4}.
		\end{equation}
		Since $0\leq1-\Delta_\lambda<\omega<1$, we have
		$\rho\in[0,1)$. Moreover,
		\begin{equation}\label{gt_nonsym_one_minus_rho}
			1-\rho
			=
			1-\frac{1-\Delta_\lambda}{\omega}
			=
			\frac{3\Delta_\lambda}{4\omega}.
		\end{equation}
		For $\ell\geq L_\Delta+1$, the second branch of
		\eqref{gt_nonsym_psi_definition} gives
		\begin{equation}\label{gt_nonsym_weighted_psi_tail_identity}
			\omega^{-\ell}\psi_\ell
			=
			\ell\Delta_\lambda
			\left(
			\frac{1-\Delta_\lambda}{\omega}
			\right)^\ell
			\overset{\eqref{gt_nonsym_rho}}{=}
			\Delta_\lambda\ell\rho^\ell.
		\end{equation}
		Therefore,
		\begin{align}
			\sum_{\ell=L_\Delta+1}^{\infty}
			\omega^{-\ell}\psi_\ell
			&\overset{\eqref{gt_nonsym_weighted_psi_tail_identity}}{=}
			\Delta_\lambda
			\sum_{\ell=L_\Delta+1}^{\infty}
			\ell\rho^\ell
			\leq
			\Delta_\lambda
			\sum_{\ell=1}^{\infty}
			\ell\rho^\ell,
			\label{gt_nonsym_psi_tail_weighted_start}
		\end{align}
		where the inequality follows because $\rho\in[0,1)$ and hence all
		terms $\ell\rho^\ell$ are nonnegative. Using
		$\sum_{\ell=1}^{\infty}\ell x^\ell=x/(1-x)^2$ for $|x|<1$, we obtain
		\begin{align}
			\sum_{\ell=L_\Delta+1}^{\infty}
			\omega^{-\ell}\psi_\ell
			&\overset{\eqref{gt_nonsym_psi_tail_weighted_start}}{\leq}
			\Delta_\lambda
			\sum_{\ell=1}^{\infty}\ell\rho^\ell
			\nonumber\\
			&=
			\Delta_\lambda
			\frac{\rho}{(1-\rho)^2}
			\nonumber\\
			&\overset{\eqref{gt_nonsym_one_minus_rho}}{=}
			\frac{16\rho\omega^2}
			{9\Delta_\lambda}
			\leq
			\frac{16}{9\Delta_\lambda}.
			\label{gt_nonsym_psi_tail_weighted}
		\end{align}
		Consequently,
		\begin{equation}\label{gt_nonsym_psi_weighted_sum}
			\sum_{\ell=0}^{\infty}
			\omega^{-\ell}\psi_\ell
			\overset{\eqref{gt_nonsym_psi_first_weighted},
				\eqref{gt_nonsym_psi_tail_weighted}}{\leq}
			\left(
			\frac{4}{3e}
			+\frac{16}{9}
			\right)
			\frac1{\Delta_\lambda}
			\leq
			\frac{8.6}{\Delta_\lambda}.
		\end{equation}
		Thus,
		\begin{equation}\label{gt_nonsym_psi_sums}
			\sum_{\ell=0}^{\infty}\psi_\ell^2
			\overset{\eqref{gt_nonsym_psi_square_sum}}{\leq}
			\frac{2.6}{\Delta_\lambda},
			\qquad
			\sum_{\ell=0}^{\infty}
			\omega^{-\ell}\psi_\ell
			\overset{\eqref{gt_nonsym_psi_weighted_sum}}{\leq}
			\frac{8.6}{\Delta_\lambda}.
		\end{equation}
		
		It remains to bound the geometric part of
		\eqref{gt_nonsym_mode_envelope_uniform}. Since
		$0<\Delta_\lambda\leq1$, we have
		$0\leq1-\Delta_\lambda<1$. Hence, using the geometric-series identity
		$\sum_{\ell=0}^{\infty}x^\ell=1/(1-x)$ for $|x|<1$, with
		$x=(1-\Delta_\lambda)^2$,
		\begin{align}
			\sum_{\ell=0}^{\infty}
			(1-\Delta_\lambda)^{2\ell}
			&=
			\frac1{1-(1-\Delta_\lambda)^2}
			=
			\frac1{\Delta_\lambda(2-\Delta_\lambda)}
			\leq
			\frac1{\Delta_\lambda}.
			\label{gt_nonsym_geometric_square}
		\end{align}
		Moreover, recalling \eqref{gt_nonsym_rho},
		\begin{align}
			\sum_{\ell=0}^{\infty}
			\omega^{-\ell}(1-\Delta_\lambda)^\ell
			&=
			\sum_{\ell=0}^{\infty}
			\left(
			\frac{1-\Delta_\lambda}{\omega}
			\right)^\ell
			\nonumber\\
			&\overset{\eqref{gt_nonsym_rho}}{=}
			\sum_{\ell=0}^{\infty}\rho^\ell
			=
			\frac1{1-\rho}
			\nonumber\\
			&\overset{\eqref{gt_nonsym_one_minus_rho}}{=}
			\frac{4\omega}{3\Delta_\lambda}
			\leq
			\frac{4}{3\Delta_\lambda}
			\leq
			\frac{2}{\Delta_\lambda}.
			\label{gt_nonsym_geometric_weighted}
		\end{align}
		Therefore,
		\begin{equation}\label{gt_nonsym_geometric_sums}
			\sum_{\ell=0}^{\infty}
			(1-\Delta_\lambda)^{2\ell}
			\overset{\eqref{gt_nonsym_geometric_square}}{\leq}
			\frac1{\Delta_\lambda},
			\qquad
			\sum_{\ell=0}^{\infty}
			\omega^{-\ell}(1-\Delta_\lambda)^\ell
			\overset{\eqref{gt_nonsym_geometric_weighted}}{\leq}
			\frac2{\Delta_\lambda}.
		\end{equation}	
		Combining
		\eqref{gt_nonsym_mode_envelope_uniform},
		\eqref{gt_nonsym_psi_sums}, and
		\eqref{gt_nonsym_geometric_sums} gives
		\begin{align}
			\sum_{\ell=0}^{\infty}
			\omega^{-\ell}
			\max_{i\geq2}|m_\ell(\lambda_i)|
			&\overset{\eqref{gt_nonsym_mode_envelope_uniform}}{\leq}
			\beta_W
			\sum_{\ell=0}^{\infty}\omega^{-\ell}\psi_\ell
			+
			\sum_{\ell=0}^{\infty}
			\omega^{-\ell}(1-\Delta_\lambda)^\ell
			\nonumber\\
			&\overset{\eqref{gt_nonsym_psi_sums},
				\eqref{gt_nonsym_geometric_sums}}{\leq}
			\frac{8.6\beta_W+2}{\Delta_\lambda}
			\leq
			\frac{11\beta_W}{\Delta_\lambda},
			\label{gt_nonsym_envelope_weighted_final}
		\end{align}
		where the last inequality uses $\beta_W\geq1$, which holds because $|1-\lambda_i|\geq1-|\lambda_i|$.	Similarly, using $(a+b)^2\leq2a^2+2b^2$ together with
		\eqref{gt_nonsym_mode_envelope_uniform},
		\begin{align}
			\sum_{\ell=0}^{\infty}
			\max_{i\geq2}|m_\ell(\lambda_i)|^2
			&\overset{\eqref{gt_nonsym_mode_envelope_uniform}}{\leq}
			2\beta_W^2
			\sum_{\ell=0}^{\infty}\psi_\ell^2
			+
			2\sum_{\ell=0}^{\infty}
			(1-\Delta_\lambda)^{2\ell}
			\nonumber\\
			&\overset{\eqref{gt_nonsym_psi_sums},
				\eqref{gt_nonsym_geometric_sums}}{\leq}
			\frac{5.2\beta_W^2+2}{\Delta_\lambda}
			\leq
			\frac{8\beta_W^2}{\Delta_\lambda},
			\label{gt_nonsym_envelope_square_final}
		\end{align}
		where again $\beta_W\geq1$. This proves
		\eqref{gt_nonsym_envelope_sums}.
	\end{proof}
	
	\subsection{Proof of Lemma \ref{lem:gt_nonsym_modulus}}
	\label{app:gt_nonsym_modulus}
	\begin{proof}[Proof of Lemma \ref{lem:gt_nonsym_modulus}]
		Let $r=|\lambda_i|$. Using \eqref{ml_poly},
		\begin{align}
			|m_\ell(\lambda_i)|
			&=
			r^\ell
			\bigl|\ell-(\ell+1)\lambda_i\bigr|
			\nonumber\\
			&\leq
			r^\ell\bigl[\ell+(\ell+1)r\bigr]
			=
			|m_\ell(-r)|.
			\label{gt_nonsym_modulus_fixed_r}
		\end{align}
		The inequality is tight for the negative real point $\lambda_i=-r$.
		Moreover, the function
		$r^\ell[\ell+(\ell+1)r]$ is increasing in $r\in[0,1]$.
		Since $r\leq1-\Delta_\lambda$, it follows that
		\begin{equation}\label{gt_nonsym_modulus_envelope}
			\max_{i\geq2}|m_\ell(\lambda_i)|
			\leq
			|m_\ell(-(1-\Delta_\lambda))|
			=
			(1-\Delta_\lambda)^\ell
			\bigl[
			\ell+(\ell+1)(1-\Delta_\lambda)
			\bigr].
		\end{equation}
		Applying the exact modal identity
		\eqref{gt_nonsym_modal_energy} at
		$z=-(1-\Delta_\lambda)$ gives
		\[
		\sum_{\ell=0}^{\infty}
		\max_{i\geq2}|m_\ell(\lambda_i)|^2
		\leq
		\frac{
			2(1-\Delta_\lambda)^2
		}{
			\Delta_\lambda^3(2-\Delta_\lambda)
		}
		\leq
		\frac1{\Delta_\lambda^3}.
		\]
		Moreover, \eqref{gt_nonsym_modulus_envelope} implies
		\[
		\max_{i\geq2}|m_\ell(\lambda_i)|
		\leq
		(1-\Delta_\lambda)^\ell(2\ell+1).
		\]
		Let
		\begin{equation}\label{gt_nonsym_modulus_rho}
			\rho
			\define
			\frac{1-\Delta_\lambda}{\omega}.
		\end{equation}
		Since $\omega=1-\Delta_\lambda/4$ by \eqref{gt_nonsym_omega},
		\begin{equation}\label{gt_nonsym_modulus_rho_gap}
			1-\rho
			=
			\frac{\omega-(1-\Delta_\lambda)}{\omega}
			=
			\frac{3\Delta_\lambda}{4\omega}
			\geq
			\frac{3\Delta_\lambda}{4},
		\end{equation}
		where the last inequality follows from $\omega\leq1$.	Moreover, \eqref{gt_nonsym_modulus_envelope} and
		$1-\Delta_\lambda\leq1$ give
		\begin{align}
			\omega^{-\ell}
			\max_{i\geq2}|m_\ell(\lambda_i)|
			&\overset{\eqref{gt_nonsym_modulus_envelope}}{\leq}
			\omega^{-\ell}(1-\Delta_\lambda)^\ell
			\bigl[
			\ell+(\ell+1)(1-\Delta_\lambda)
			\bigr]
			\nonumber\\
			&\leq
			(2\ell+1)
			\left(
			\frac{1-\Delta_\lambda}{\omega}
			\right)^\ell
			\overset{\eqref{gt_nonsym_modulus_rho}}{=}
			(2\ell+1)\rho^\ell.
			\label{gt_nonsym_modulus_weighted_envelope}
		\end{align}
		Since $\rho\in[0,1)$, using the standard identities
		$\sum_{\ell=0}^{\infty}\rho^\ell=(1-\rho)^{-1}$ and
		$\sum_{\ell=0}^{\infty}\ell\rho^\ell
		=\rho(1-\rho)^{-2}$, we obtain
		\begin{align}
			\sum_{\ell=0}^{\infty}
			\omega^{-\ell}
			\max_{i\geq2}|m_\ell(\lambda_i)|
			&\overset{\eqref{gt_nonsym_modulus_weighted_envelope}}{\leq}
			\sum_{\ell=0}^{\infty}(2\ell+1)\rho^\ell
			\nonumber\\
			&=
			\frac{2\rho}{(1-\rho)^2}
			+\frac1{1-\rho}
			\nonumber\\
			&\overset{\eqref{gt_nonsym_modulus_rho_gap}}{\leq}
			\frac{32}{9\Delta_\lambda^2}
			+\frac{4}{3\Delta_\lambda}
			\nonumber\\
			&\leq
			\left(\frac{32}{9}+\frac43\right)
			\frac1{\Delta_\lambda^2}
			=
			\frac{44}{9\Delta_\lambda^2}
			\leq
			\frac5{\Delta_\lambda^2},
			\label{gt_nonsym_modulus_weighted_sum}
		\end{align}
		where we used $0<\Delta_\lambda\leq1$.
		This proves \eqref{gt_nonsym_disc_sums}.
	\end{proof}
	
	\subsection{Proof of Corollary \ref{cor:gt_nonsym_two_region}}
	\label{app:gt_nonsym_two_region}
	\begin{proof}[Proof of Corollary \ref{cor:gt_nonsym_two_region}]
		For $j\in\{1,2\}$, define
		\begin{equation}\label{gt_nonsym_two_region_al}
			a_\ell^{(j)}
			\define
			\max_{i\in\mathcal I_j}|m_\ell(\lambda_i)|,
		\end{equation}
		with the convention $a_\ell^{(j)}=0$ if $\mathcal I_j=\emptyset$.
		Since $\mathcal I_1$ and $\mathcal I_2$ partition the nontrivial
		eigenvalues,
		\begin{equation}\label{gt_nonsym_two_region_split}
			\max_{i\geq2}|m_\ell(\lambda_i)|
			=
			\max\{a_\ell^{(1)},a_\ell^{(2)}\}
			\leq
			a_\ell^{(1)}+a_\ell^{(2)},
		\end{equation}
		and similarly,
		\begin{equation}\label{gt_nonsym_two_region_split_square}
			\max_{i\geq2}|m_\ell(\lambda_i)|^2
			\leq
			\bigl(a_\ell^{(1)}\bigr)^2
			+
			\bigl(a_\ell^{(2)}\bigr)^2.
		\end{equation}
		For the modes in $\mathcal I_1$, the conditions in
		\eqref{gt_nonsym_two_region_conditions} are precisely those used in
		Lemma~\ref{lem:gt_nonsym_envelope}, with $\Delta_\lambda$ and $\beta_W$ replaced by $p$ and $c$,
		respectively; the assumption $c\geq1$ plays the role of
		$\beta_W\geq1$. Moreover, \eqref{gt_nonsym_two_region_omega} gives
		$\omega_{\mathrm{reg}}\geq1-p/4$. Hence, if $\omega_p\define1-p/4$, then
		$\omega_{\mathrm{reg}}^{-\ell}\leq\omega_p^{-\ell}$ for every $\ell\geq0$, so the
		same proof yields
		\begin{equation}\label{gt_nonsym_two_region_one_bounds}
			\sum_{\ell\geq0}\omega_{\mathrm{reg}}^{-\ell}a_\ell^{(1)}
			\leq
			\frac{11c}{p},
			\qquad
			\sum_{\ell\geq0}\bigl(a_\ell^{(1)}\bigr)^2
			\leq
			\frac{8c^2}{p}.
		\end{equation}
		For the modes in $\mathcal I_2$, only the modulus bound
		$1-|\lambda_i|\geq\nu$ is used. Thus, repeating the proof of
		Lemma~\ref{lem:gt_nonsym_modulus} with
		$\Delta_\lambda$ replaced by $\nu$ gives the corresponding bounds.
		Indeed, \eqref{gt_nonsym_two_region_omega} implies
		$\omega_{\mathrm{reg}}\geq1-\nu/4$, and therefore
		\begin{equation}\label{gt_nonsym_two_region_two_bounds}
			\sum_{\ell\geq0}\omega_{\mathrm{reg}}^{-\ell}a_\ell^{(2)}
			\leq
			\frac5{\nu^2},
			\qquad
			\sum_{\ell\geq0}\bigl(a_\ell^{(2)}\bigr)^2
			\leq
			\frac1{\nu^3}.
		\end{equation}	
		Using \eqref{gt_nonsym_B_eig_bound},
		\eqref{gt_nonsym_two_region_split}, and
		\eqref{gt_nonsym_two_region_one_bounds}--%
		\eqref{gt_nonsym_two_region_two_bounds}, we obtain
		\begin{align}
			\mathcal G_1(\omega_{\mathrm{reg}})
			&\leq
			\kappa_W
			\sum_{\ell\geq0}
			\omega_{\mathrm{reg}}^{-\ell}
			\max_{i\geq2}|m_\ell(\lambda_i)|
			\leq
			\kappa_W
			\left(
			\frac{11c}{p}
			+
			\frac5{\nu^2}
			\right).
			\label{gt_nonsym_two_region_G1_proof}
		\end{align}
		Likewise, by \eqref{gt_nonsym_two_region_split_square},
		\begin{align}
			\mathcal G_2
			&\leq
			\kappa_W^2
			\sum_{\ell\geq0}
			\max_{i\geq2}|m_\ell(\lambda_i)|^2
			\leq
			\kappa_W^2
			\left(
			\frac{8c^2}{p}
			+
			\frac1{\nu^3}
			\right).
			\label{gt_nonsym_two_region_G2_proof}
		\end{align}
		Together with \eqref{gt_nonsym_network_gains}, these estimates prove
		\eqref{gt_nonsym_two_region_gains}.

	Let $r=|z|<1$. Since $1-r^2=(1-r)(1+r)$ and
	$|1-z|=\beta_z(1-r)$, Lemma~\ref{lem:gt_nonsym_modal_energy} gives
	\[
	\begin{aligned}
		\sum_{\ell\geq0}|m_\ell(z)|^2
		&=\frac{2r^2}{(1+r)^3}\frac{\beta_z^2}{1-r}\\
		&\leq\frac{\beta_z^2}{4(1-r)}
		\leq\frac{1}{(1-r)^3}.
	\end{aligned}
	\]
	The first inequality uses $(1+r)^2\geq4r$ and $1+r\geq2r$,
	so $(1+r)^3\geq8r^2$. The second uses
	$\beta_z=|1-z|/(1-r)\leq2/(1-r)$.	For $i\in\mathcal I_1$, we know both $\beta_{\lambda_i}\leq c$
	and $1-|\lambda_i|\geq p$, so the first bound gives $c^2/(4p)$
	per eigenvalue. For $i\in\mathcal I_2$, we only use
	$1-|\lambda_i|\geq\nu$, so the second bound gives $1/\nu^3$
	per eigenvalue. Since $\theta_j=|\mathcal I_j|/N$, summing
	these contributions in \eqref{gt_nonsym_avg_gain} gives
	\[
	\mathcal G_{2,\mathrm{av}}
	\leq \kappa_W^2
	\left(\frac{\theta_1c^2}{4p}+\frac{\theta_2}{\nu^3}\right),
	\]
	as claimed.
		
		For a real spectrum, choose
		\begin{equation}\label{gt_nonsym_real_partition}
			\mathcal I_1
			=
			\{i\in\{2,\ldots,N\}:\lambda_i\geq0\},
			\qquad
			\mathcal I_2
			=
			\{i\in\{2,\ldots,N\}:\lambda_i<0\}.
		\end{equation}
		If $\lambda_i\in[0,1-p]$, then
		$1-|\lambda_i|\geq p$ and
		\[
		\frac{|1-\lambda_i|}{1-|\lambda_i|}
		=
		\frac{1-\lambda_i}{1-\lambda_i}
		=
		1.
		\]
		Thus the first region satisfies
		\eqref{gt_nonsym_two_region_conditions} with $c=1$.
		For $\lambda_i\in[-1+\nu,0)$,
		\[
		1-|\lambda_i|
		=
		1+\lambda_i
		\geq
		\nu,
		\]
		so the second-region condition also holds. Hence, if the nontrivial spectrum lies in $[-1+\nu,1-p]$, this
		partition gives $c=1$. If $W$ is symmetric, an orthonormal
		eigenbasis also gives $\kappa_W=1$.
	\end{proof}

	\section{Proofs of the refined convergence results}
	\label{app:refined_convergence}
	
	\subsection{Proof of Lemma \ref{lemma:gt_nonsym_weighted_gap}}
	\label{app:gt_nonsym_weighted_gap}
	\begin{proof}[Proof of Lemma \ref{lemma:gt_nonsym_weighted_gap}]
		Since
		\begin{equation}\label{gt_nonsym_gamma_contraction}
			\gamma^2
			=
			1-\mu\alpha+\frac{\mu^2\alpha^2}{4}
			\geq
			1-\mu\alpha,
		\end{equation}
		we have
		\begin{equation}\label{gt_nonsym_gamma_ratio}
			\frac{1-\mu\alpha}{\gamma}
			\overset{\eqref{gt_nonsym_gamma_contraction}}{\leq}
			\gamma.
		\end{equation}
		Moreover, since $1-\gamma=\mu\alpha/2$, the geometric-series identity
		gives
		\begin{align}
			\sum_{k=0}^{K}\gamma^{-k}
			&=
			\frac{\gamma^{-K}-\gamma}{1-\gamma}
			\leq
			\frac{\gamma^{-K}}{1-\gamma}
			=
			\frac{2\gamma^{-K}}{\mu\alpha}.
			\label{gt_nonsym_geometric_sum}
		\end{align}
		Finally, $\mu\leq L$ and $\alpha\leq1/(4L)$ imply
		$\mu\alpha\leq1/4$, and hence
		\begin{equation}\label{gt_nonsym_gamma_inverse}
			1\leq\gamma^{-1}
			=
			\frac{1}{1-\mu\alpha/2}
			\leq2.
		\end{equation}	
		Taking total expectations in \eqref{gt_semi_average_before_step} and using
		$2(1-2\alpha L)\geq1$ and
		$1+2\alpha L\leq3/2$, we obtain
		\begin{equation}\label{gt_nonsym_average_bound_simple}
			\bar e_{k+1}
			\leq
			(1-\mu\alpha)\bar e_k
			-\alpha\tilde f_k
			+\frac{3\alpha L}{2}E_k
			+\frac{\alpha^2\sigma^2}{N}.
		\end{equation}
		Multiplying \eqref{gt_nonsym_average_bound_simple} by
		$\gamma^{-(k+1)}$ gives
		\begin{align}
			\gamma^{-(k+1)}\bar e_{k+1}
			\leq{}&
			\frac{1-\mu\alpha}{\gamma}\,
			\gamma^{-k}\bar e_k
			-\alpha\gamma^{-1}\gamma^{-k}\tilde f_k
			\nonumber\\
			&+
			\gamma^{-1}\gamma^{-k}
			\left(
			\frac{3\alpha L}{2}E_k
			+\frac{\alpha^2\sigma^2}{N}
			\right).
		\end{align}
		Since $\tilde f_k\geq0$, \eqref{gt_nonsym_gamma_ratio} and
		$\gamma^{-1}\geq1$ imply
		\begin{align}
			\gamma^{-(k+1)}\bar e_{k+1}
			\leq{}&
			\gamma\,\gamma^{-k}\bar e_k
			-\alpha\gamma^{-k}\tilde f_k
			\nonumber\\
			&+
			\gamma^{-1}\gamma^{-k}
			\left(
			\frac{3\alpha L}{2}E_k
			+\frac{\alpha^2\sigma^2}{N}
			\right).
			\label{gt_nonsym_weighted_average_step}
		\end{align}
		Summing \eqref{gt_nonsym_weighted_average_step} over
		$k=0,\ldots,K$ and collecting the terms involving $\bar e_k$,
		each intermediate term $\gamma^{-k}\bar e_k$, $1\leq k\leq K$,
		appears once with coefficient $\gamma$ from iteration $k$ and once
		with coefficient $-1$ from iteration $k-1$. Hence its net coefficient
		is $\gamma-1=-(1-\gamma)$. The initial and terminal terms remain
		$\gamma\bar e_0$ and $-\gamma^{-(K+1)}\bar e_{K+1}$, respectively.
		Therefore,
		\begin{align}
			\alpha\sum_{k=0}^{K}\gamma^{-k}\tilde f_k
			\leq{}&
			\gamma\bar e_0
			-(1-\gamma)
			\sum_{k=1}^{K}\gamma^{-k}\bar e_k
			-\gamma^{-(K+1)}\bar e_{K+1}
			\nonumber\\
			&+
			\gamma^{-1}
			\left(
			\frac{3\alpha L}{2}\mathcal D_K
			+\frac{\alpha^2\sigma^2}{N}
			\sum_{k=0}^{K}\gamma^{-k}
			\right).
			\label{gt_nonsym_weighted_gap_sum}
		\end{align}
		The second and third terms on the right-hand side are nonpositive.
		Dropping them and using $\gamma\leq1$,
		\eqref{gt_nonsym_gamma_inverse}, and
		\eqref{gt_nonsym_geometric_sum}, we obtain
		\begin{align}
			\alpha\sum_{k=0}^{K}\gamma^{-k}\tilde f_k
			&\leq
			\bar e_0
			+3\alpha L\mathcal D_K
			+\frac{4\alpha\sigma^2}{\mu N}\gamma^{-K}.
		\end{align}
		Dividing by $\alpha$ proves
		\eqref{gt_nonsym_weighted_gap}.
	\end{proof}
	
	\subsection{Proof of Lemma \ref{lemma:gt_nonsym_weighted_disagreement}}
	\label{app:gt_nonsym_weighted_disagreement}
	\begin{proof}[Proof of Lemma \ref{lemma:gt_nonsym_weighted_disagreement}]
		By \eqref{gt_nonsym_alpha_g},
		$\gamma=1-\mu\alpha/2\geq1-\Delta_\lambda/4=\omega$.
		Since $\mathcal G_1(\cdot)$ is nonincreasing in its argument, this
		implies
		$\mathcal G_1(\gamma)\leq\mathcal G_1(\omega)$ and
		$\mathcal G_1(1)\leq\mathcal G_1(\omega)$.
		
		Using
		$\widehat{\mathbf r}^t
		=(\I-\mathbf\Pi)(\mathbf u^t+\v^t)$ and $\ell\define k-1-t$,
		the representation \eqref{gt_nonsym_convolution_identity} 
		becomes
		$\widehat{\x}^{k}
		=\T^k+\alpha\D^k+\alpha\mathbf S^k$, where $\T^k=	\widetilde{\W}^{k}\widehat{\x}^0
		-\alpha k\,\widetilde{\W}^{k}\widetilde{\y}^0$,
		\[
		\D^k
		\define
		\sum_{\ell=0}^{k-1}
		\mathbf M_\ell(\I-\mathbf\Pi)
		\mathbf u^{k-1-\ell},
		\qquad
		\mathbf S^k
		\define
		\sum_{\ell=0}^{k-1}
		\mathbf M_\ell(\I-\mathbf\Pi)
		\v^{k-1-\ell}.
		\]
		Hence,
		\begin{equation}\label{gt_nonsym_three_way}
			E_k
			\leq
			\frac3N\Ex\|\T^k\|^2
			+\frac{3\alpha^2}{N}\Ex\|\D^k\|^2
			+\frac{3\alpha^2}{N}\Ex\|\mathbf S^k\|^2.
		\end{equation}
		Multiplying \eqref{gt_nonsym_three_way} by $\gamma^{-k}$ and summing
		over $k=0,\ldots,K$ gives
		\begin{align}
			\mathcal D_K
			&=
			\sum_{k=0}^{K}\gamma^{-k}E_k
			\nonumber\\
			&\overset{\eqref{gt_nonsym_three_way}}{\leq}
			3\sum_{k=0}^{K}
			\gamma^{-k}\frac1N\Ex\|\T^k\|^2
			+
			3\alpha^2\sum_{k=0}^{K}
			\gamma^{-k}\frac1N\Ex\|\D^k\|^2
			+
			3\alpha^2\sum_{k=0}^{K}
			\gamma^{-k}\frac1N\Ex\|\mathbf S^k\|^2.
			\label{gt_nonsym_three_way_weighted}
		\end{align}
		We bound the three terms separately.
		
		\emph{(i) Initialization term.}
		By \eqref{gt_nonsym_power_decay},
		\[
		\frac1N\Ex\|\T^k\|^2
		\leq
		2\kappa_W^2(1-\Delta_\lambda)^{2k}
		\left(
		E_0+\alpha^2k^2\varsigma_\star^2
		\right).
		\]
		Since $\gamma\geq1-\Delta_\lambda$, we have
		$(1-\Delta_\lambda)^2/\gamma\leq1-\Delta_\lambda$. Therefore,
		\begin{align}
			\sum_{k=0}^{K}
			\gamma^{-k}\frac1N\Ex\|\T^k\|^2
			&\leq
			2\kappa_W^2
			\left(
			\frac{E_0}{\Delta_\lambda}
			+\frac{2\alpha^2\varsigma_\star^2}
			{\Delta_\lambda^3}
			\right)
			\nonumber\\
			&=
			\mathcal H_0,
			\label{gt_nonsym_initialization_sum}
		\end{align}
		where we used
		$\sum_{k\geq0}(1-\Delta_\lambda)^k
		=1/\Delta_\lambda$ and
		$\sum_{k\geq0}k^2(1-\Delta_\lambda)^k
		\leq2/\Delta_\lambda^3$.
		
		\emph{(ii) Optimization term.}
		By the triangle inequality and $\|\I-\mathbf\Pi\|_2=1$,
		\[
		\|\D^k\|
		\leq
		\sum_{\ell=0}^{k-1}
		\|\mathbf M_\ell\|_2
		\|\mathbf u^{k-1-\ell}\|.
		\]
		Applying the weighted Cauchy--Schwarz inequality
		\[
		\left(\sum_{\ell}a_\ell b_\ell\right)^2
		\leq
		\left(\sum_{\ell}a_\ell\right)
		\left(\sum_{\ell}a_\ell b_\ell^2\right),
		\qquad a_\ell\geq0,
		\]
		with
		$a_\ell=\|\mathbf M_\ell\|_2$ and
		$b_\ell=\|\mathbf u^{k-1-\ell}\|$, gives
		\begin{align}
			\|\D^k\|^2
			&\leq
			\left(
			\sum_{\ell=0}^{k-1}\|\mathbf M_\ell\|_2
			\right)
			\left(
			\sum_{\ell=0}^{k-1}
			\|\mathbf M_\ell\|_2
			\|\mathbf u^{k-1-\ell}\|^2
			\right)
			\nonumber\\
			&\leq
			\mathcal G_1(1)
			\sum_{\ell=0}^{k-1}
			\|\mathbf M_\ell\|_2
			\|\mathbf u^{k-1-\ell}\|^2,
		\end{align}
		where the last inequality follows from
		$\sum_{\ell=0}^{k-1}\|\mathbf M_\ell\|_2
		\leq\sum_{\ell=0}^{\infty}\|\mathbf M_\ell\|_2
		=\mathcal G_1(1)$. Since $\D^0=0$ by definition, we have $
		\sum_{k=0}^{K}
		\gamma^{-k}\frac1N\Ex\|\D^k\|^2
		=
		\sum_{k=1}^{K}
		\gamma^{-k}\frac1N\Ex\|\D^k\|^2$.
		Multiplying the preceding bound by $\gamma^{-k}/N$, taking
		expectations, and summing over $k=1,\ldots,K$ gives
		\begin{align}
			\sum_{k=0}^{K}
			\gamma^{-k}\frac1N\Ex\|\D^k\|^2
			&\leq
			\mathcal G_1(1)
			\sum_{k=1}^{K}
			\sum_{\ell=0}^{k-1}
			\gamma^{-k}\|\mathbf M_\ell\|_2
			\frac1N\Ex
			\|\mathbf u^{k-1-\ell}\|^2.
			\label{gt_nonsym_D_double_sum}
		\end{align}
		To reindex the double sum, set $t=k-1-\ell$, so that
		$k=t+\ell+1$. The original range
		$1\leq k\leq K$, $0\leq\ell\leq k-1$ is equivalent to
		$0\leq\ell\leq K-1$, $0\leq t\leq K-1-\ell$.
		Moreover,
		$\gamma^{-k}=\gamma^{-(\ell+1)}\gamma^{-t}$. Hence,
		\begin{align}
			\sum_{k=0}^{K}
			\gamma^{-k}\frac1N\Ex\|\D^k\|^2
			&\overset{\eqref{gt_nonsym_D_double_sum}}{\leq}
			\mathcal G_1(1)
			\sum_{\ell=0}^{K-1}
			\gamma^{-(\ell+1)}
			\|\mathbf M_\ell\|_2
			\sum_{t=0}^{K-1-\ell}
			\gamma^{-t}
			\frac1N\Ex\|\mathbf u^t\|^2
			\nonumber\\
			&\leq
			\gamma^{-1}\mathcal G_1(1)
			\left(
			\sum_{\ell=0}^{\infty}
			\gamma^{-\ell}\|\mathbf M_\ell\|_2
			\right)
			U_K
			\nonumber\\
			&=
			\gamma^{-1}
			\mathcal G_1(1)\mathcal G_1(\gamma)U_K,
			\label{gt_nonsym_D_reindexed}
		\end{align}
		where
		$U_K\define
		\sum_{t=0}^{K}\gamma^{-t}N^{-1}\Ex\|\mathbf u^t\|^2$.
		The second inequality follows by bounding the inner $t$-sum by $U_K$
		and extending the nonnegative $\ell$-sum to infinity.	Since $\gamma\geq\omega$,
		$\mathcal G_1(\gamma)\leq\mathcal G_1(\omega)$, while
		$\mathcal G_1(1)\leq\mathcal G_1(\omega)$. Moreover,
		\eqref{gt_nonsym_gamma_inverse} gives $\gamma^{-1}\leq2$.
		Consequently,
		\begin{equation}\label{gt_nonsym_D_weighted}
			\sum_{k=0}^{K}
			\gamma^{-k}\frac1N\Ex\|\D^k\|^2
			\overset{\eqref{gt_nonsym_D_reindexed}}{\leq}
			2\mathcal G_1(\omega)^2U_K.
		\end{equation}
		By \eqref{eq:appendix_u_split} and
		\eqref{eq:appendix_cocoercive},
		\begin{align}
			\frac1N\Ex\|\mathbf u^t\|^2
			&\leq
			\frac{2L^2}{N}\Ex\|\widehat{\x}^t\|^2
			+4L\Ex\!\left[f(\bar x^t)-f(x^\star)\right]
			\nonumber\\
			&=
			2L^2E_t+4L\tilde f_t.
			\label{gt_nonsym_u_bound_gap}
		\end{align}
		Consequently,
		$U_K
		\leq
		2L^2\mathcal D_K
		+4L\sum_{t=0}^{K}\gamma^{-t}\tilde f_t$.
		Combining this with \eqref{gt_nonsym_D_weighted} yields
		\begin{equation}\label{gt_nonsym_D_weighted_final}
			3\alpha^2
			\sum_{k=0}^{K}
			\gamma^{-k}\frac1N\Ex\|\D^k\|^2
			\leq
			6\alpha^2\mathcal G_1(\omega)^2
			\left(
			2L^2\mathcal D_K
			+4L\sum_{t=0}^{K}\gamma^{-t}\tilde f_t
			\right).
		\end{equation}
		
		\emph{(iii) Stochastic term.}
		For any $s\neq t$, assume without loss of generality that $s<t$.
		Since $\v^s$ is $\cF^t$-measurable and
		$\Ex[\v^t\mid\cF^t]=0$, for every deterministic matrix $A$,
		\[
		\Ex\!\left[(\v^s)\tran A\v^t\right]
		=
		\Ex\!\left[
		(\v^s)\tran A
		\Ex[\v^t\mid\cF^t]
		\right]
		=0.
		\]
		Hence noise terms from distinct iterations are orthogonal in
		expectation, so all cross terms in
		$\Ex\|\mathbf S^k\|^2$ vanish and
		\[
		\frac1N\Ex\|\mathbf S^k\|^2
		=
		\frac1N
		\sum_{\ell=0}^{k-1}
		\Ex\left\|
		\mathbf M_\ell(\I-\mathbf\Pi)
		\v^{k-1-\ell}
		\right\|^2.
		\]
		Moreover, write $\mathbf M_\ell=B_\ell\otimes I_d$ with
		$B_\ell\define m_\ell(W-\Pi)$, and let $v_j^t$ be the $j$th block of
		$\v^t$. 		Since the $i$th block of $\mathbf M_\ell\v^t$ is
		$\sum_j(B_\ell)_{ij}v_j^t$, and, given $\cF^t$, the blocks $v_j^t$ are
		independent with zero mean, the cross terms
		$\Ex[(v_j^t)\tran v_{j'}^t\mid\cF^t]$, $j\neq j'$, vanish. Together with
		\eqref{gt_semi_variance}, this gives
		\[
		\Ex[\|\mathbf M_\ell\v^t\|^2\mid\cF^t]
		=\sum_{i,j}(B_\ell)_{ij}^2\Ex[\|v_j^t\|^2\mid\cF^t]
		\leq\sigma^2\|B_\ell\|_F^2=\frac{\sigma^2}{d}\|\mathbf M_\ell\|_F^2.
		\]
		Therefore, since $\mathbf M_\ell(\I-\mathbf\Pi)=\mathbf M_\ell$,
		\[
		\begin{aligned}
			\frac1N\Ex\|\mathbf S^k\|^2
			&\leq\frac{\sigma^2}{Nd}\sum_{\ell=0}^{k-1}\|\mathbf M_\ell\|_F^2
			\leq\sigma^2\mathcal G_{2,\mathrm{av}}\leq\sigma^2\mathcal G_2,
		\end{aligned}
		\]
		and, by \eqref{gt_nonsym_geometric_sum},
		\begin{equation}\label{gt_nonsym_noise_weighted}
			3\alpha^2\sum_{k=0}^{K}\gamma^{-k}\frac1N\Ex\|\mathbf S^k\|^2
			\leq\frac{6\alpha\sigma^2\mathcal G_{2,\mathrm{av}}}{\mu}\gamma^{-K}.
		\end{equation}
		Applying \eqref{gt_nonsym_initialization_sum},
		\eqref{gt_nonsym_D_weighted_final}, and
		\eqref{gt_nonsym_noise_weighted} to
		\eqref{gt_nonsym_three_way_weighted} proves
		\eqref{gt_nonsym_weighted_disagreement}.
	\end{proof}

	\subsection{Proof of Theorem \ref{thm:gt_nonsym_refined}}
	\label{app:gt_nonsym_refined}
	\begin{proof}[Proof of Theorem \ref{thm:gt_nonsym_refined}]
		Let
		\[
		\phi\define6\alpha^2\mathcal G_1(\omega)^2.
		\]
		The last bound in \eqref{gt_nonsym_stepsize_conditions} gives
		$\alpha^2L^2\mathcal G_1(\omega)^2\leq1/168$, and hence
		\begin{equation}\label{gt_nonsym_phi_bound}
			\phi L^2\leq\frac1{28}.
		\end{equation}
		Substituting \eqref{gt_nonsym_weighted_gap} into
		\eqref{gt_nonsym_weighted_disagreement} gives
		\[
		\bigl(1-14\phi L^2\bigr)\mathcal D_K
		\leq
		3\mathcal H_0
		+\frac{4\phi L\bar e_0}{\alpha}
		+
		\left(
		\frac{16\phi L\sigma^2}{\mu N}
		+
		\frac{6\alpha\sigma^2\mathcal G_{2,\mathrm{av}}}{\mu}
		\right)\gamma^{-K}.
		\]
		Hence, by \eqref{gt_nonsym_phi_bound},
		\begin{equation}\label{gt_nonsym_C_solved}
			\mathcal D_K
			\leq
			6\mathcal H_0
			+\frac{8\phi L\bar e_0}{\alpha}
			+
			\left(
			\frac{32\phi L\sigma^2}{\mu N}
			+
			\frac{12\alpha\sigma^2\mathcal G_{2,\mathrm{av}}}{\mu}
			\right)\gamma^{-K}.
		\end{equation}
		Next, unrolling \eqref{gt_nonsym_average_bound_simple} and dropping
		$-\alpha\tilde f_k$ yields
		\begin{align}
			\bar e_K
			\leq{}&
			(1-\mu\alpha)^K\bar e_0
			+\frac{3\alpha L}{2}
			\sum_{k=0}^{K-1}
			(1-\mu\alpha)^{K-1-k}E_k
			\nonumber\\
			&+
			\frac{\alpha^2\sigma^2}{N}
			\sum_{j=0}^{K-1}(1-\mu\alpha)^j.
		\end{align}
		Since $1-\mu\alpha\leq\gamma$,
		$\gamma^{-1}\leq2$, and
		$\sum_{j\geq0}(1-\mu\alpha)^j=(\mu\alpha)^{-1}$,
		\begin{equation}\label{gt_nonsym_aK_intermediate}
			\bar e_K
			\leq
			\gamma^K\bar e_0
			+3\alpha L\gamma^K\mathcal D_K
			+\frac{\alpha\sigma^2}{\mu N}.
		\end{equation}
		Multiplying \eqref{gt_nonsym_C_solved} by
		$3\alpha L\gamma^K$ gives
		\[
		3\alpha L\gamma^K\mathcal D_K
		\leq
		\gamma^K
		\left(
		18\alpha L\mathcal H_0
		+24\phi L^2\bar e_0
		\right)
		+\frac{96\phi L^2\alpha\sigma^2}{\mu N}
		+
		\frac{36\alpha^2L\sigma^2\mathcal G_{2,\mathrm{av}}}{\mu}.
		\]
		By \eqref{gt_nonsym_phi_bound},
		\[
		24\phi L^2\leq\frac67,
		\qquad
		\frac{96\phi L^2\alpha\sigma^2}{\mu N}
		\leq
		\frac{4\alpha\sigma^2}{\mu N}.
		\]
		Moreover, using \eqref{gt_nonsym_H0} and
		$\alpha\leq1/(4L)$,
		\[
		18\alpha L\mathcal H_0
		\leq
		\frac{9\kappa_W^2}{\Delta_\lambda}E_0
		+\widetilde c_2\alpha^3.
		\]
		Substituting these bounds into
		\eqref{gt_nonsym_aK_intermediate}, using
		$(1+6/7)\bar e_0\leq2\bar e_0$, and
		$\gamma^K\leq\exp(-\mu\alpha K/2)$ proves
		\eqref{gt_nonsym_main_neighborhood}.
	\end{proof}
	
	\subsection{Proof of Corollary \ref{cor:gt_nonsym_linear_speedup}}
	\label{app:gt_nonsym_linear_speedup}
	\begin{proof}[Proof of Corollary \ref{cor:gt_nonsym_linear_speedup}]
		By \eqref{gt_nonsym_alpha_bar},
		\[
		\bar\alpha\leq\frac{1}{4L},
		\qquad
		\mu\bar\alpha\leq\frac{\Delta_\lambda}{2},
		\qquad
		168L^2\mathcal G_1(\omega)^2\bar\alpha^2\leq1.
		\]
				Hence every $\alpha\in(0,\bar\alpha]$, in particular the choice in
		\eqref{eq:tuning_alpha}, satisfies \eqref{gt_nonsym_stepsize_conditions},
		so \eqref{gt_nonsym_main_neighborhood} holds on $(0,\bar\alpha]$. Applying
		Lemma~\ref{lem:log_stepsize} to
		\eqref{gt_nonsym_main_neighborhood} gives
		\[
		\bar e_K
		=
		\widetilde\cO\!\left(
		\frac{c_1}{\mu K}
		+
		\frac{c_2}{\mu^2K^2}
		\right)
		+
		\left(
		\mathcal I_0+\widetilde c_2\bar\alpha^3
		\right)
		\exp\!\left(
		-\frac{\mu\bar\alpha K}{2}
		\right).
		\]
		By the definitions of $c_1$ and $c_2$ in
		Theorem~\ref{thm:gt_nonsym_refined},
		\[
		\frac{c_1}{\mu K}
		=
		\frac{5\sigma^2}{\mu^2NK},
		\qquad
		\frac{c_2}{\mu^2K^2}
		=
		\frac{36L\sigma^2\mathcal G_{2,\mathrm{av}}}{\mu^3K^2}.
		\]
		Moreover,
		\[
		\mathcal I_0+\widetilde c_2\bar\alpha^3
		=
		2\bar e_0
		+\frac{9\kappa_W^2}{\Delta_\lambda}E_0
		+\frac{
			72L\kappa_W^2\varsigma_\star^2
		}{
			\Delta_\lambda^3
		}
		\bar\alpha^3
		=
		C_{\mathrm{init}}.
		\]
		Absorbing the numerical constants into
		$\widetilde\cO(\cdot)$ proves
		\eqref{gt_nonsym_big_o_neighborhood}.
	\end{proof}

	\section{Proof of Theorem~\ref{thm:jlb_joint}}
	\label{app:finite_lower_bound}
	\begin{proof}[Proof of Theorem~\ref{thm:jlb_joint}]
		The stepsize condition gives
		$0<2\alpha\mu\leq\nu^2/6\leq1/96$, and hence
		$0<\rho_{\mathrm{lb}}<1$. Since the gradients are globally
		Lipschitz and the noises have finite second moments, all iterates
		have finite second moments at every finite horizon.
		
		The proof has three steps. First, we reduce the averaged-iterate error to
		a lower bound involving the weighted amplitude $\mathcal B_K$. Second, we
		use the positive and negative eigenmodes of $W$ to obtain a forcing
		inequality for $a_k$. Finally, we sum this inequality with geometric
		weights, absorb the feedback terms, and lower-bound the accumulated
		noise and heterogeneity contributions.
		
		\noindent
		\emph{Step 1: Reduction to $\mathcal B_K$.}
		By \eqref{eq:glb_objectives} and $\one\tran s^{\mathrm{lb}}=0$,
		\[
		\overline{\grad\f}(\x^k)
		=
		\frac1N\sum_i f_i'(x_i^k)
		=
		2\mu\bar x^k+\frac{\mu}{N}\|\x^k\|_1.
		\]
		Thus the centroid recursion \eqref{gt_semi_centroid} becomes
		$
		\bar x^{k+1}
		=
		\rho_{\mathrm{lb}}\bar x^k
		-\frac{\alpha\mu}{N}\|\x^k\|_1
		-\alpha\bar v^k$.
		Iterating from $\bar x^0=-\sqrt{\bar e_0}$ gives
		\begin{equation}\label{eq:glb_scalar_identity}
			\begin{aligned}
				\bar x^K
				&=
				-\sqrt{\bar e_0}\,\rho_{\mathrm{lb}}^K
				-U_K^{\mathrm{lb}}
				+G_K^{\mathrm{lb}},\\
				U_K^{\mathrm{lb}}
				&\define
				\frac{\alpha\mu}{N}
				\sum_{k<K}
				\rho_{\mathrm{lb}}^{K-1-k}\|\x^k\|_1
				\geq0,\\
				G_K^{\mathrm{lb}}
				&\define
				-\alpha\sum_{k<K}
				\rho_{\mathrm{lb}}^{K-1-k}\bar v^k.
			\end{aligned}
		\end{equation}
		Define
		\begin{equation}\label{eq:jlb_B_definition}
			a_k
			\define
			\frac1N\Ex\|\x^k\|_1,
			\qquad
			\mathcal B_K
			\define
			\sum_{k<K}
			\rho_{\mathrm{lb}}^{K-1-k}a_k.
		\end{equation}
		Then, by \eqref{eq:glb_scalar_identity} and
		\eqref{eq:jlb_B_definition},
		$\Ex U_K^{\mathrm{lb}}=\alpha\mu\mathcal B_K$.
		
		By \eqref{eq:glb_noise_initialization},
		$\bar v^k=N^{-1}\sum_i\xi_i^k$ are independent across iterations and
		$\bar v^k\sim\mathcal N(0,\sigma^2/N)$. Therefore
		$G_K^{\mathrm{lb}}$ is centered Gaussian and
		\begin{equation}\label{eq:jlb_G_variance}
			\begin{aligned}
				\Ex[(G_K^{\mathrm{lb}})^2]
				&=
				\frac{\alpha^2\sigma^2}{N}
				\sum_{k=0}^{K-1}
				\rho_{\mathrm{lb}}^{2(K-1-k)}\\
				&=
				\frac{\alpha^2\sigma^2}{N}
				\frac{1-\rho_{\mathrm{lb}}^{2K}}
				{1-\rho_{\mathrm{lb}}^2}\\
				&=
				\frac{\alpha\sigma^2}
				{4\mu N(1-\alpha\mu)}
				\bigl(1-\rho_{\mathrm{lb}}^{2K}\bigr),
			\end{aligned}
		\end{equation}
		where we used
		$1-\rho_{\mathrm{lb}}^2
		=
		4\alpha\mu(1-\alpha\mu)$.
		No independence between $G_K^{\mathrm{lb}}$ and
		$U_K^{\mathrm{lb}}$ is required.
		
		Since the first two terms on the right of
		\eqref{eq:glb_scalar_identity} are nonpositive,
		$\bar x^K\leq G_K^{\mathrm{lb}}$ pathwise. Therefore, on
		$\{G_K^{\mathrm{lb}}\leq0\}$,
		$
		\bar x^K
		\leq
		G_K^{\mathrm{lb}}
		\leq0$,
		and hence $|\bar x^K|\geq|G_K^{\mathrm{lb}}|$. Thus
		\begin{equation}\label{eq:jlb_G_event}
			\bar e_K
			\geq
			\Ex\!\left[
			(G_K^{\mathrm{lb}})^2
			\mathds 1\{G_K^{\mathrm{lb}}\leq0\}
			\right].
		\end{equation}
		Since $G_K^{\mathrm{lb}}$ is centered Gaussian, its distribution is
		symmetric about zero and $z^2$ is even. Hence
		\begin{equation}\label{eq:jlb_G_symmetry}
			\Ex\!\left[
			(G_K^{\mathrm{lb}})^2
			\mathds 1\{G_K^{\mathrm{lb}}\leq0\}
			\right]
			=
			\frac12\Ex[(G_K^{\mathrm{lb}})^2].
		\end{equation}
		Combining \eqref{eq:jlb_G_event} and
		\eqref{eq:jlb_G_symmetry},
		\begin{equation}\label{eq:jlb_noise_lower}
			\bar e_K
			\geq
			\frac12\Ex[(G_K^{\mathrm{lb}})^2].
		\end{equation}
		On the other hand, taking expectations in
		\eqref{eq:glb_scalar_identity} and using
		$\Ex U_K^{\mathrm{lb}}
		=
		\alpha\mu\mathcal B_K$ and $\Ex G_K^{\mathrm{lb}}=0$ gives
		\begin{equation}\label{eq:jlb_centroid_mean}
			\Ex\bar x^K
			=
			-\sqrt{\bar e_0}\rho_{\mathrm{lb}}^K
			-\alpha\mu\mathcal B_K.
		\end{equation}
		Therefore, Jensen's inequality and
		\eqref{eq:jlb_centroid_mean} imply
		\begin{equation}\label{eq:jlb_Jensen_lower}
			\begin{aligned}
				\bar e_K
				&\geq
				(\Ex\bar x^K)^2\\
				&=
				\left(
				\sqrt{\bar e_0}\rho_{\mathrm{lb}}^K
				+\alpha\mu\mathcal B_K
				\right)^2\\
				&\geq
				\bar e_0\rho_{\mathrm{lb}}^{2K}
				+(\alpha\mu\mathcal B_K)^2.
			\end{aligned}
		\end{equation}
		Combining \eqref{eq:jlb_noise_lower} with weight $2/3$ and
		\eqref{eq:jlb_Jensen_lower} with weight $1/3$, and then using
		\eqref{eq:jlb_G_variance} together with
		$(1-\alpha\mu)^{-1}\geq1$, gives
		\begin{equation}\label{eq:glb_scalar_master}
			\bar e_K
			\geq
			\frac13\,\bar e_0\rho_{\mathrm{lb}}^{2K}
			+\frac{\alpha\sigma^2}{12\mu N}
			\bigl(1-\rho_{\mathrm{lb}}^{2K}\bigr)
			+\frac13(\alpha\mu\mathcal B_K)^2.
		\end{equation}
		Thus the initialization and averaged-noise terms are already obtained.
		It remains to lower-bound $\mathcal B_K$ by the positive-mode noise,
		negative-mode noise, and heterogeneity forcing.
		
		\noindent
		\emph{Step 2: Modal forcing of $\mathcal B_K$.}
		We use two families of eigenvectors. The $+$ family extracts the mode
		$1-p$ and is aligned with $s^{\mathrm{lb}}$, so it captures both
		positive-mode noise and heterogeneity. The $-$ family extracts the mode
		$-1+\nu$ and captures the negative-mode noise.
		
		Let $e_j$ denote the $j$th standard basis vector of $\real^N$.
		We use the real vectors
		\begin{equation}\label{eq:jlb_vectors}
			v_{j,-}
			\define
			\frac{e_{2j-1}-e_{2j}}2,
			~ j=1,\ldots,N/2,
			~~
			v_{j,+}
			\define
			\frac{e_{4j-3}+e_{4j-2}-e_{4j-1}-e_{4j}}4,
			~ j=1,\ldots,N/4
		\end{equation}
		satisfying
		\begin{equation}\label{eq:jlb_vector_norms}
			\|v_{j,-}\|_1=\|v_{j,+}\|_1=1,
			\qquad
			\|v_{j,-}\|^2=\frac12,
			\qquad
			\|v_{j,+}\|^2=\frac14,
		\end{equation}
		and
		\begin{equation}\label{eq:jlb_vector_alignment}
			v_{j,-}\tran s^{\mathrm{lb}}=0,
			\qquad
			v_{j,+}\tran s^{\mathrm{lb}}=1.
		\end{equation}
		From the properties of $W$ in \eqref{eq:jlb_mixing}, we have for
		$v=v_{j,\pm}$,
		\begin{equation}\label{eq:jlb_modal_properties}
			v\tran\one=0,
			\qquad
			v\tran W=z_\pm v\tran,
			\qquad
			v\tran\widetilde W^\ell=z_\pm^\ell v\tran,
			\qquad
			v\tran\mathbf M_\ell=m_\ell^\pm v\tran,
		\end{equation}
		where $z_+=1-p$ for $v=v_{j,+}$, $z_-=-1+\nu$ for $v=v_{j,-}$, and
		$
		m_\ell^\pm
		\define
		\ell z_\pm^\ell-(\ell+1)z_\pm^{\ell+1}$
		as in \eqref{ml_poly}.
		Define
		\begin{equation}\label{eq:jlb_S_definition}
			S_k^\pm
			\define
			\left(
			\sum_{\ell=0}^{k-1}|m_\ell^\pm|^2
			\right)^{1/2}.
		\end{equation}
		In our construction,
		$\widehat{\x}^0=\zero$,
		$\widetilde{\y}^0=\grad\f(\x^\star)
		=\varsigma_\star s^{\mathrm{lb}}$, and
		\[
		\widehat{\mathbf r}^t
		=
		(I_N-\Pi)(\mathbf u^t+\v^t),
		\qquad
		\v^t=\bxi^t,
		\qquad
		u_i^t
		=
		f_i'(x_i^t)-f_i'(0)
		=
		2\mu x_i^t+\mu|x_i^t|.
		\]
		Since $v\tran\one=0$, we have
		$v\tran\widehat{\x}^k=v\tran\x^k$.	Multiplying \eqref{gt_nonsym_convolution_identity} by $v\tran$,
		and using \eqref{eq:jlb_modal_properties} and
		$v\tran(I_N-\Pi)=v\tran$, gives
		\begin{equation}\label{eq:glb_projected}
			\begin{aligned}
				v\tran\x^k
				&=
				-\alpha\varsigma_\star
				(v\tran s^{\mathrm{lb}})kz_\pm^k
				+Z_{j,k}^\pm
				+\alpha\sum_{t<k}
				m_{k-1-t}^\pm\,v\tran\mathbf u^t,\\
				Z_{j,k}^\pm
				&\define
				\alpha\sum_{t<k}
				m_{k-1-t}^\pm\,v\tran\bxi^t.
			\end{aligned}
		\end{equation}
		The variables $v\tran\bxi^t$ are independent across $t$ and
		distributed as $\mathcal N(0,\sigma^2\|v\|^2)$. Therefore, using
		\eqref{eq:jlb_vector_norms} and \eqref{eq:jlb_S_definition},
		\begin{equation}\label{eq:jlb_projected_noise}
			\text{Var}(Z_{j,k}^\pm)
			=
			\alpha^2\sigma^2\|v_{j,\pm}\|^2(S_k^\pm)^2,
		\end{equation}
		and consequently \cite{papoulis2002probability}
		\begin{equation}\label{eq:jlb_projected_noise_absolute}
			\Ex|Z_{j,k}^-|
			=
			\frac{\alpha\sigma S_k^-}{\sqrt\pi},
			\qquad
			\Ex|Z_{j,k}^+|
			=
			\frac{\alpha\sigma S_k^+}{\sqrt{2\pi}}.
		\end{equation}
		For deterministic $b$ and centered real Gaussian $Z$, Jensen's
		inequality gives
		\[
		\Ex|b+Z|
		\geq
		|\Ex(b+Z)|
		=
		|b|.
		\]
		Also, the map $b\mapsto\Ex|b+Z|$ is convex and even, hence minimized
		at zero, so $\Ex|b+Z|\geq\Ex|Z|$. Combining these two observations,
		\begin{equation}\label{eq:jlb_shifted_gaussian}
			\Ex|b+Z|
			\geq
			\max\{|b|,\Ex|Z|\}
			\geq
			\frac{|b|+\Ex|Z|}{2}.
		\end{equation}
		Write \eqref{eq:glb_projected} as
		\[
		v\tran\x^k
		=
		b_{j,k}^\pm+Z_{j,k}^\pm+F_{j,k}^\pm,
		\]
		where
		\begin{equation}\label{eq:jlb_bF_definition}
			b_{j,k}^\pm
			\define
			-\alpha\varsigma_\star
			(v\tran s^{\mathrm{lb}})kz_\pm^k,
			\qquad
			F_{j,k}^\pm
			\define
			\alpha\sum_{t<k}
			m_{k-1-t}^\pm v\tran\mathbf u^t.
		\end{equation}
		By \eqref{eq:jlb_shifted_gaussian},
		\[
		\frac12
		\left(
		|b_{j,k}^\pm|
		+\Ex|Z_{j,k}^\pm|
		\right)
		\leq
		\Ex|b_{j,k}^\pm+Z_{j,k}^\pm|.
		\]
		On the other hand, since
		$v\tran\x^k=b_{j,k}^\pm+Z_{j,k}^\pm+F_{j,k}^\pm$,
		the triangle inequality gives
		\[
		\Ex|b_{j,k}^\pm+Z_{j,k}^\pm|
		\leq
		\Ex|v\tran\x^k|
		+\Ex|F_{j,k}^\pm|.
		\]
		Moreover, by the definition of $F_{j,k}^\pm$,
		\[
		\Ex|F_{j,k}^\pm|
		\leq
		\alpha\sum_{t<k}
		|m_{k-1-t}^\pm|
		\Ex|v\tran\mathbf u^t|.
		\]
		Combining these inequalities gives
		\begin{equation}\label{eq:jlb_projected_forcing}
			\frac12
			\left(
			|b_{j,k}^\pm|
			+\Ex|Z_{j,k}^\pm|
			\right)
			\leq
			\Ex|v\tran\x^k|
			+
			\alpha\sum_{t<k}
			|m_{k-1-t}^\pm|
			\Ex|v\tran\mathbf u^t|.
		\end{equation}
		Within either family, the supports of $v_{j,\pm}$ partition the agents.
		Using the coefficients in \eqref{eq:jlb_vectors},
		\begin{subequations}\label{eq:jlb_family_averages}
			\begin{alignat}{2}
				\frac2N\sum_{j=1}^{N/2}
				\Ex|v_{j,-}\tran\x^k|
				&\leq a_k,
				&\qquad
				\frac4N\sum_{j=1}^{N/4}
				\Ex|v_{j,+}\tran\x^k|
				&\leq a_k,
				\label{eq:jlb_family_x}
				\\
				\frac2N\sum_{j=1}^{N/2}
				\Ex|v_{j,-}\tran\mathbf u^t|
				&\leq La_t,
				&\qquad
				\frac4N\sum_{j=1}^{N/4}
				\Ex|v_{j,+}\tran\mathbf u^t|
				&\leq La_t,
				\label{eq:jlb_family_u}
			\end{alignat}
		\end{subequations}
		where \eqref{eq:jlb_family_u} uses
		$|u_i^t|\leq L|x_i^t|$.
		
		For the negative family,
		\eqref{eq:jlb_vector_alignment} gives
		$b_{j,k}^-=0$. For the positive family,
		\eqref{eq:jlb_vector_alignment} and $z_+=1-p$ give
		\begin{equation}\label{eq:jlb_b_values}
			|b_{j,k}^-|=0,
			\qquad
			|b_{j,k}^+|
			=
			\alpha\varsigma_\star k(1-p)^k.
		\end{equation}
		Averaging \eqref{eq:jlb_projected_forcing} over the negative family with
		weight $2/N$, and using \eqref{eq:jlb_family_averages},
		\eqref{eq:jlb_projected_noise_absolute}, and
		$|b_{j,k}^-|=0$, gives
		\[
		a_k
		+\alpha L\sum_{t<k}|m_{k-1-t}^-|a_t
		\geq
		\frac{\alpha\sigma S_k^-}{2\sqrt{\pi}}.
		\]
		Likewise, averaging over the positive family with weight $4/N$ gives
		\[
		a_k
		+\alpha L\sum_{t<k}|m_{k-1-t}^+|a_t
		\geq
		\frac{\alpha}{2}
		\left(
		\frac{\sigma S_k^+}{\sqrt{2\pi}}
		+\varsigma_\star k(1-p)^k
		\right).
		\]
		Adding these two inequalities yields
		\begin{equation}\label{eq:glb_forcing}
			\begin{aligned}
				2a_k
				+\alpha L\sum_{t<k}
				\bigl(
				|m_{k-1-t}^+|
				+
				|m_{k-1-t}^-|
				\bigr)a_t
				&\geq
				\frac{\alpha}{2}
				\left[
				\frac{\sigma S_k^-}{\sqrt\pi}
				+\frac{\sigma S_k^+}{\sqrt{2\pi}}
				+\varsigma_\star k(1-p)^k
				\right]
				\\
				&\gtrsim
				\alpha
				\left[
				\sigma(S_k^++S_k^-)
				+\varsigma_\star k(1-p)^k
				\right].
			\end{aligned}
		\end{equation}
		
		\noindent
		\emph{Step 3: Weighted accumulation.}
		
		\noindent
		\emph{Feedback absorption.}
		Define
		\begin{equation}\label{eq:jlb_Q_definition}
			Q_\pm
			\define
			\sum_{\ell\geq0}
			\rho_{\mathrm{lb}}^{-(\ell+1)}|m_\ell^\pm|.
		\end{equation}
		From \eqref{eq:jlb_stepsize} and $L=3\mu$,
		\[
		2\alpha\mu
		\leq
		\frac{\min\{p,\nu^2\}}6.
		\]
		Hence
		\begin{equation}\label{eq:jlb_rho_gaps}
			\rho_{\mathrm{lb}}-1+p
			=
			p-2\alpha\mu
			\geq
			\frac56p,
			\qquad
			\rho_{\mathrm{lb}}-1+\nu
			=
			\nu-2\alpha\mu
			\geq
			\frac{23}{24}\nu,
		\end{equation}
		where the second inequality uses $\nu\leq1/4$.
		In particular,
		$\rho_{\mathrm{lb}}>|z_\pm|$, so the sums in
		\eqref{eq:jlb_Q_definition} are finite.
		
		Since
		\[
		|m_\ell^+|
		\leq
		(1-p)^\ell(\ell p+1-p),
		\qquad
		|m_\ell^-|
		=
		(1-\nu)^\ell
		[\ell+(\ell+1)(1-\nu)],
		\]
		geometric summation together with \eqref{eq:jlb_rho_gaps} gives
		\begin{equation}\label{eq:jlb_Q_bounds}
			\begin{aligned}
				Q_+
				&\leq
				\frac{p(1-p)}
				{(\rho_{\mathrm{lb}}-1+p)^2}
				+
				\frac{1-p}
				{\rho_{\mathrm{lb}}-1+p}
				\leq
				\frac3p,\\
				Q_-
				&=
				\frac{(1-\nu)(1+\rho_{\mathrm{lb}})}
				{(\rho_{\mathrm{lb}}-1+\nu)^2}
				\leq
				\frac3{\nu^2}.
			\end{aligned}
		\end{equation}
		Using \eqref{eq:jlb_stepsize} in
		\eqref{eq:jlb_Q_bounds},
		\begin{equation}\label{eq:jlb_feedback_gain}
			\alpha LQ_+\leq\frac34,
			\qquad
			\alpha LQ_-\leq\frac34,
			\qquad
			2+\alpha L(Q_++Q_-)\leq\frac72.
		\end{equation}
		
		For later use, for either sign, the weighted feedback convolution
		satisfies
		\begin{equation}\label{eq:jlb_convolution_bound}
			\begin{aligned}
				&\sum_{k<K}
				\rho_{\mathrm{lb}}^{K-1-k}
				\sum_{t<k}
				|m_{k-1-t}^\pm|a_t\\
				&\qquad=
				\sum_{t<K-1}
				\rho_{\mathrm{lb}}^{K-1-t}a_t
				\sum_{\ell\leq K-2-t}
				\rho_{\mathrm{lb}}^{-(\ell+1)}
				|m_\ell^\pm|\\
				&\qquad\leq
				Q_\pm\mathcal B_K,
			\end{aligned}
		\end{equation}
		where we used $\ell=k-1-t$ and exchanged the order of summation.
		
		\noindent
		\emph{Finite-window forcing.}
		We next lower-bound the three weighted forcing terms in
		\eqref{eq:glb_forcing}.
		For $a\in\{p,\nu\}$,
		\begin{equation}\label{eq:jlb_T_bounds}
			T_a
			=
			\left\lceil\frac1{4a}\right\rceil,
			\qquad
			\frac1{4a}\leq T_a\leq\frac1{2a},
			\qquad
			(1-a)^{T_a}\geq\frac12,
		\end{equation}
		where the last inequality follows from Bernoulli's inequality.
		
		For $\ell<T_p$,
		\eqref{eq:jlb_T_bounds} gives
		\[
		(1-p)^\ell\geq\frac12,
		\qquad
		1-p-\ell p
		\geq
		1-T_pp
		\geq
		\frac12,
		\]
		and hence
		\begin{equation}\label{eq:jlb_positive_kernel}
			|m_\ell^+|
			\geq
			\frac14,
			\qquad
			\ell<T_p.
		\end{equation}
		Likewise, for $\ell<T_\nu$,
		\begin{equation}\label{eq:jlb_negative_kernel}
			|m_\ell^-|
			\geq
			(1-\nu)^{\ell+1}(\ell+1)
			\geq
			\frac{\ell+1}{2}.
		\end{equation}
		Using \eqref{eq:jlb_S_definition},
		\eqref{eq:jlb_positive_kernel}, and
		\eqref{eq:jlb_negative_kernel},
		\[
		(S_{T_p}^+)^2
		\geq
		\frac{T_p}{16},
		\qquad
		(S_{T_\nu}^-)^2
		\geq
		\frac14\sum_{a=1}^{T_\nu}a^2.
		\]
		Since $T_p\geq(4p)^{-1}$,
		$T_\nu\geq(4\nu)^{-1}$, and
		$\sum_{a=1}^{T}a^2\geq T^3/3$, we obtain
		\begin{equation}\label{eq:jlb_energy_bounds}
			(S_{T_p}^+)^2
			\geq
			\frac1{64p},
			\qquad
			(S_{T_\nu}^-)^2
			\geq
			\frac1{768\nu^3}.
		\end{equation}
		For the heterogeneity term, retain only $1\leq k\leq T_p$.
		By \eqref{eq:jlb_stepsize} and \eqref{eq:jlb_T_bounds},
		\[
		\rho_{\mathrm{lb}}^{T_p}
		\geq
		1-2\alpha\mu T_p
		\geq
		\frac12,
		\qquad
		(1-p)^{T_p}\geq\frac12.
		\]
		Therefore, for $1\leq k\leq T_p$,
		\begin{equation}\label{eq:jlb_heterogeneity_window}
			\rho_{\mathrm{lb}}^{K-1-k}(1-p)^k
			\geq
			\frac14
			\rho_{\mathrm{lb}}^{K-T_p-1}.
		\end{equation}
		Using \eqref{eq:jlb_heterogeneity_window} and
		$\sum_{k=1}^{T_p}k\geq T_p^2/2$, together with
		\eqref{eq:jlb_T_bounds}, gives
		\begin{equation}\label{eq:jlb_heterogeneity_sum}
			\sum_{k=1}^{K-1}
			\rho_{\mathrm{lb}}^{K-1-k}k(1-p)^k
			\geq
			\frac{
				\rho_{\mathrm{lb}}^{K-T_p-1}
			}{128p^2}.
		\end{equation}
		For the two noise terms, $S_k^\pm$ is nondecreasing in $k$. Hence, for
		either sign and any $T<K$,
		\[
		\begin{aligned}
			\sum_{k<K}
			\rho_{\mathrm{lb}}^{K-1-k}S_k^\pm
			&\geq
			S_T^\pm
			\sum_{k=T}^{K-1}
			\rho_{\mathrm{lb}}^{K-1-k}
			=
			S_T^\pm
			\frac{1-\rho_{\mathrm{lb}}^{K-T}}
			{1-\rho_{\mathrm{lb}}}.
		\end{aligned}
		\]
		Since $1-\rho_{\mathrm{lb}}=2\alpha\mu$, this gives
		\begin{equation}\label{eq:jlb_noise_sum}
			\sum_{k<K}
			\rho_{\mathrm{lb}}^{K-1-k}S_k^\pm
			\geq
			\frac{S_T^\pm}{2\alpha\mu}
			\left(
			1-\rho_{\mathrm{lb}}^{K-T}
			\right),
		\end{equation}
		where $T=T_p$ for the $+$ mode and $T=T_\nu$ for the $-$ mode.
		
		Now multiply \eqref{eq:glb_forcing} by
		$\rho_{\mathrm{lb}}^{K-1-k}$ and sum over $k<K$. By
		\eqref{eq:jlb_feedback_gain} and
		\eqref{eq:jlb_convolution_bound}, its weighted left-hand side is at most
		$\tfrac72\mathcal B_K$. Therefore,
		\[
		\begin{aligned}
			\frac72\mathcal B_K
			\geq \frac{\alpha}{2}\Bigg[
			&\frac{\sigma}{\sqrt\pi}
			\sum_{k<K}\rho_{\mathrm{lb}}^{K-1-k}S_k^-
			+
			\frac{\sigma}{\sqrt{2\pi}}
			\sum_{k<K}\rho_{\mathrm{lb}}^{K-1-k}S_k^+
			+
			\varsigma_\star
			\sum_{k=1}^{K-1}
			\rho_{\mathrm{lb}}^{K-1-k}k(1-p)^k
			\Bigg].
		\end{aligned}
		\]
		Using \eqref{eq:jlb_energy_bounds},
		\eqref{eq:jlb_heterogeneity_sum}, and
		\eqref{eq:jlb_noise_sum}, and then multiplying by $\alpha\mu$, gives
		\begin{equation}\label{eq:jlb_amplitudes}
			\begin{aligned}
				\alpha\mu\mathcal B_K
				\geq{}&
				\frac{\alpha\sigma}
				{112\sqrt{2\pi p}}
				\left(
				1-\rho_{\mathrm{lb}}^{K-T_p}
				\right)
				+
				\frac{\alpha\sigma}
				{224\sqrt{3\pi}\,\nu^{3/2}}
				\left(
				1-\rho_{\mathrm{lb}}^{K-T_\nu}
				\right)
				\\
				&+
				\frac{\alpha^2\mu\varsigma_\star}
				{896p^2}
				\rho_{\mathrm{lb}}^{K-T_p-1}.
			\end{aligned}
		\end{equation}
		All three terms are nonnegative, so
		$(a+b+c)^2\geq a^2+b^2+c^2$. Moreover, since $p\leq\nu$, we have $T_\nu\leq T_p$ and hence $1-\rho_{\mathrm{lb}}^{K-T_\nu}\geq1-\rho_{\mathrm{lb}}^{K-T_p}$. Substituting
		\eqref{eq:jlb_amplitudes} into the last term of
		\eqref{eq:glb_scalar_master} gives the coefficients $1/3$ for the
		initialization term, $1/12$ for the averaged-noise term,
		$1/(75264\pi)$ for the positive-mode noise term,
		$1/(451584\pi)$ for the negative-mode noise term, and
		$1/2408448$ for the heterogeneity term. Taking the smallest of these
		universal constants proves \eqref{eq:jlb_main}.
	\end{proof}

	\section{Proof of the baseline convergence result}
	\label{sec:gt_two_recursion_appendix}
	
	This appendix proves Theorem~\ref{thm:gt_baseline_convergence}; it uses
	Lemma~\ref{lemma:gt_shared_estimates} and Lemma~\ref{lem:log_stepsize}
	of Appendix~\ref{app:shared_estimates}. The proof of
	Theorem~\ref{thm:cubic_network_lower_bound} is given in
	Appendix~\ref{app:cubic_network_lower_bound}.
	
	\subsection{Network contraction and input-gain bounds}

	We begin with the two norm bounds needed for the joint error recursion.
	
	\begin{lemma}[\bfseries\small Contraction and input-gain bounds]
		\label{lemma:gt_semi_matrix_bounds}
		Suppose Assumption~\ref{assump:gt_semi_network} holds and
		$0<\lambda<1$. Then the matrices $\C$ and $\B$ in
		\eqref{gt_semi_C_B} satisfy
		\begin{align}
			\|\C\|_2
			&\leq
			1-\delta,\quad
			\|\B\|_2^2
			\leq
			Q_\lambda .
			\label{gt_semi_matrix_bounds}
		\end{align}
	\end{lemma}
	
	\begin{proof}
		From \eqref{gt_semi_C_B},
		\[
		\C
		=
		\begin{bmatrix}
			\widetilde{\W} & \zero\\
			\zero & \widetilde{\W}
		\end{bmatrix}
		\begin{bmatrix}
			\I & -\eta\I\\
			\zero & \I
		\end{bmatrix}.
		\]
		The first factor has norm $\lambda$, while the triangle inequality gives
		norm at most $1+\eta$ for the second. Hence
		\[
		\|\C\|_2
		\leq
		\lambda(1+\eta)
		=
		\frac{1+\lambda}{2}
		=
		1-\delta.
		\]
		For $\B$, since $W$ is nonnegative and doubly stochastic,
		$\|W\|_2\leq\sqrt{\|W\|_1\|W\|_\infty}=1$, and hence
		$\|I-W\|_2\leq2$. Therefore,
		\[
		\|\B\|_2^2
		\leq
		\lambda^2+\frac4{\eta^2}
		=
		\lambda^2+\frac{16\lambda^2}{\Delta_\lambda^2}
		\leq
		\frac{17\lambda^2}{\Delta_\lambda^2}
		=
		Q_\lambda.
		\]
	\end{proof}
	
	\subsection{Lyapunov inequality}
	
	We next derive the single recursion from which the baseline theorem follows.

	\begin{lemma}[\bfseries\small One-step Lyapunov recursion]
		\label{lemma:gt_semi_one_step}
		Suppose Assumptions~\ref{assump:gt_semi_functions}--
		\ref{assump:gt_semi_noise} hold, $0<\lambda<1$, and
		\begin{equation}\label{eq:appendix_stepsize_conditions}
			0<\alpha\leq
			\min\left\{
			\frac{1}{4L},
			\frac{\delta}{L\sqrt{24Q_\lambda}},
			\frac{\delta}{2\mu}
			\right\}.
		\end{equation}
		Then, for every $k\geq0$,
		\begin{equation}
			\mathcal V_{k+1}
			\leq
			\left(1-\frac{\mu\alpha}{2}\right)\mathcal V_k
			+\frac{\alpha^2\sigma^2}{N}
			+\frac{3\alpha^3Q_\lambda L\sigma^2}{\delta},
			\label{eq:lyapunov_one_step}
		\end{equation}
		where $\mathcal V_k
		\define
		\bar e_k+\frac{3\alpha L}{\delta}b_k$.
	\end{lemma}
	
	\begin{proof}
		We first bound the average error and the joint disagreement error, and then
		combine the two bounds so that the intermediate objective term cancels.
		
		\paragraph{Average error.}
		Let $\widetilde x^k\define\bar x^k-x^\star$. By
		Lemma~\ref{lemma:gt_shared_estimates}, \eqref{gt_semi_average_before_step} holds.
		Since $\alpha\leq1/(4L)$,
		$2(1-2\alpha L)\geq1$ and $1+2\alpha L\leq3/2$.
		Moreover, by \eqref{gt_semi_eta_delta_joint},
		$\|\widehat{\x}^k\|^2\leq\|\widehat{\e}^k\|^2$.
		Taking total expectations and using
		\eqref{eq:baseline_errors_main},
		\begin{align}
			\bar e_{k+1}
			&\leq
			(1-\mu\alpha)\bar e_k
			-\alpha\Ex\!\left[
			f(\bar x^k)-f(x^\star)
			\right]
			+\frac{3\alpha L}{2}b_k
			+\frac{\alpha^2\sigma^2}{N}.
			\label{eq:appendix_average_bound}
		\end{align}
		
		\paragraph{Joint disagreement error.}
		From the exact recursion \eqref{gt_semi_joint_exact},
		conditional unbiasedness of $\v^k$ gives
		\begin{align}
			\Ex\!\left[
			\|\widehat{\e}^{k+1}\|^2
			\,\middle|\,
			\cF^k
			\right]
			&=
			\|\C\widehat{\e}^k-\alpha\B\mathbf u^k\|^2
			+\alpha^2
			\Ex\!\left[
			\|\B\v^k\|^2
			\,\middle|\,
			\cF^k
			\right].
			\label{eq:appendix_joint_noise_cancel}
		\end{align}
		Applying Young's inequality
		$\|p+q\|^2\leq\theta^{-1}\|p\|^2
		+(1-\theta)^{-1}\|q\|^2$
		with $\theta=1-\delta$, and using
		Lemma~\ref{lemma:gt_semi_matrix_bounds}, gives
		\begin{equation}
			\|\C\widehat{\e}^k-\alpha\B\mathbf u^k\|^2
			\leq
			(1-\delta)\|\widehat{\e}^k\|^2
			+\frac{\alpha^2Q_\lambda}{\delta}\|\mathbf u^k\|^2.
			\label{eq:appendix_predictable_bound}
		\end{equation}
		Similarly,
		\begin{equation}
			\Ex\!\left[
			\|\B\v^k\|^2
			\,\middle|\,
			\cF^k
			\right]
			\leq
			Q_\lambda N\sigma^2.
			\label{eq:appendix_B_noise_bound}
		\end{equation}
		It remains to bound $\mathbf u^k$, defined in \eqref{gt_semi_u_definition}; by
		Lemma~\ref{lemma:gt_shared_estimates}, \eqref{eq:appendix_u_split} and
		\eqref{eq:appendix_cocoercive} hold.
		Combining
		\eqref{eq:appendix_u_split} and
		\eqref{eq:appendix_cocoercive}, and using
		$\|\widehat{\x}^k\|^2\leq\|\widehat{\e}^k\|^2$,
		\begin{equation}
			\|\mathbf u^k\|^2
			\leq
			2L^2\|\widehat{\e}^k\|^2
			+
			4LN\bigl(f(\bar x^k)-f(x^\star)\bigr).
			\label{eq:appendix_u_bound}
		\end{equation}
		Substituting
		\eqref{eq:appendix_predictable_bound}--
		\eqref{eq:appendix_u_bound}
		into \eqref{eq:appendix_joint_noise_cancel},
		taking total expectations, and dividing by $N$ gives
		\begin{align}
			b_{k+1}
			&\leq
			\left(
			1-\delta
			+\frac{2\alpha^2Q_\lambda L^2}{\delta}
			\right)b_k
			+\frac{4\alpha^2Q_\lambda L}{\delta}
			\Ex\!\left[
			f(\bar x^k)-f(x^\star)
			\right]
			+\alpha^2Q_\lambda\sigma^2.
			\label{eq:appendix_disagreement_bound}
		\end{align}
		
		\paragraph{Combination.}
		For compactness, let
		\begin{equation}
			\chi
			\define
			\frac{2\alpha^2Q_\lambda L^2}{\delta},
			\qquad
			\tau
			\define
			\frac{3\alpha L}{\delta}.
			\label{eq:appendix_chi_gamma}
		\end{equation}
		By \eqref{eq:baseline_lyapunov_main},
		$\mathcal V_k=\bar e_k+\tau b_k$.
		Multiply \eqref{eq:appendix_disagreement_bound} by $\tau$
		and add the result to \eqref{eq:appendix_average_bound}. We obtain
		\begin{align}
			\mathcal V_{k+1}
			&\leq
			(1-\mu\alpha)\bar e_k
			+
			\left[
			\frac{3\alpha L}{2}
			+\tau(1-\delta+\chi)
			\right]b_k
			\nonumber\\
			&\quad
			-
			\left[
			\alpha
			-\frac{2\tau\chi}{L}
			\right]
			\Ex\!\left[
			f(\bar x^k)-f(x^\star)
			\right]
			+
			\frac{\alpha^2\sigma^2}{N}
			+
			\tau\alpha^2Q_\lambda\sigma^2.
			\label{eq:appendix_lyapunov_expansion}
		\end{align}
		The second stepsize condition in
		\eqref{eq:appendix_stepsize_conditions}, namely $24\alpha^2Q_\lambda L^2\leq\delta^2$, gives
		\begin{align}
			\frac{2\tau\chi}{L}
			&=
			\alpha
			\frac{12\alpha^2Q_\lambda L^2}{\delta^2}
			\leq
			\frac{\alpha}{2},
			\qquad
			\chi\leq\frac{\delta}{12}.
			\label{eq:appendix_gap_absorption}
		\end{align}
		Hence the objective-gap term in
		\eqref{eq:appendix_lyapunov_expansion} is nonpositive and may be discarded. Moreover,
		\[
		\frac{3\alpha L}{2}
		+\tau(1-\delta+\chi)
		=
		\tau\left(1-\frac{\delta}{2}+\chi\right)
		\leq
		\tau\left(1-\frac{5\delta}{12}\right)
		\leq
		\tau\left(1-\frac{\mu\alpha}{2}\right),
		\]
		where the last inequality uses $2\mu\alpha\leq\delta$.
		Also,
		$1-\mu\alpha\leq1-\mu\alpha/2$.
		Substituting these estimates into
		\eqref{eq:appendix_lyapunov_expansion} proves
		\eqref{eq:lyapunov_one_step}.
	\end{proof}
	
	\subsection{Proof of Theorem~\ref{thm:gt_baseline_convergence}}
	
	Under the assumptions of Theorem~\ref{thm:gt_baseline_convergence},
	Lemma~\ref{lemma:gt_semi_one_step} applies. Set
	\begin{equation}
		r
		\define
		1-\frac{\mu\alpha}{2},
		\qquad
		D
		\define
		\frac{\alpha^2\sigma^2}{N}
		+
		\frac{3\alpha^3Q_\lambda L\sigma^2}{\delta}.
		\label{eq:appendix_r_D}
	\end{equation}
	The stepsize conditions imply $0<r<1$, and
	\eqref{eq:lyapunov_one_step} becomes
	\begin{equation}
		\mathcal V_{k+1}
		\leq
		r\mathcal V_k+D.
	\end{equation}
	Iterating,
	\begin{align}
		\mathcal V_k
		&\leq
		r^k\mathcal V_0
		+
		D\sum_{t=0}^{k-1}r^t
		\leq
		r^k\mathcal V_0
		+
		\frac{D}{1-r}.
	\end{align}
	Since $1-r=\mu\alpha/2$,
	\begin{align}
		\frac{D}{1-r}
		&=
		\frac{2\alpha\sigma^2}{\mu N}
		+
		\frac{6\alpha^2Q_\lambda L\sigma^2}{\mu\delta},
	\end{align}
	which proves \eqref{eq:baseline_linear_bound_main}.
	
	It remains to verify the initialization term. By the consensual initialization,
	$\widehat{\x}^0=\zero$. Moreover,
	$\y^0=\zero$ by \eqref{gt_semi_y_definition}, and therefore $\widetilde{\y}^0
	=
	\grad\f(\x^\star)$
	by \eqref{gt_semi_xstar_stack}. Hence, from
	\eqref{eq:baseline_errors_main},
	\begin{equation}
		b_0
		=
		\frac{\alpha^2}{\eta^2N}
		\|\grad\f(\x^\star)\|^2
		=
		\frac{\alpha^2}{\eta^2}\varsigma_\star^2.
	\end{equation}
	Using \eqref{eq:baseline_lyapunov_main}, we have $\mathcal V_0
	=
	\bar e_0
	+
	\frac{3\alpha^3L}{\eta^2\delta}
	\varsigma_\star^2$. The horizon-dependent bound \eqref{eq:baseline_tuned_rate_main} is established in Corollary~\ref{cor:tuned_rate_baseline} below.
	
	\subsection{Horizon-dependent tuning}
	We next prove the second part of
	Theorem~\ref{thm:gt_baseline_convergence}, using
	Lemma~\ref{lem:log_stepsize} of Appendix~\ref{app:shared_estimates}.
	
	\begin{corollary}[\bfseries\small Horizon-dependent baseline rate]
		\label{cor:tuned_rate_baseline}
		Under the assumptions and initialization of
		Theorem~\ref{thm:gt_baseline_convergence}, let $K\geq1$ and let
		$\bar\alpha$ be defined in \eqref{eq:baseline_stepsize_main}.
		For $\sigma^2>0$, choose $\alpha$ according to
		\eqref{eq:tuning_alpha} with the constants specified in the proof below.
		Then
		\begin{equation}
			\mathcal V_K
			=
			\widetilde\cO\!\left(
			\frac{\sigma^2}{\mu^2NK}
			+
			\frac{L\sigma^2}
			{\mu^3\Delta_\lambda^3K^2}
			\right)
			+
			\left(
			\bar e_0
			+
			\frac{24L\lambda^2\bar\alpha^3}
			{\Delta_\lambda^3}\varsigma_\star^2
			\right)
			\exp\!\left(
			-\frac{\mu\bar\alpha K}{2}
			\right).
			\label{eq:tuned_rate_bigO}
		\end{equation}
		If $\sigma^2=0$, the same statement holds with the two stochastic terms
		equal to zero by taking $\alpha=\bar\alpha$.
	\end{corollary}
	
	\begin{proof}
		For every $\alpha\in(0,\bar\alpha]$,
		\eqref{eq:baseline_linear_bound_main} and
		$1-x\leq e^{-x}$ give
		\begin{align}
			\mathcal V_K
			&\leq
			\exp\!\left(
			-\frac{\mu\alpha K}{2}
			\right)
			\left(
			\bar e_0
			+
			\frac{3L\varsigma_\star^2}{\eta^2\delta}\alpha^3
			\right)
			\nonumber\\
			&\quad
			+\frac{2\sigma^2}{\mu N}\alpha
			+\frac{6Q_\lambda L\sigma^2}{\mu\delta}\alpha^2 .
			\label{eq:appendix_pre_tuning}
		\end{align}
		Thus Lemma~\ref{lem:log_stepsize} applies with
		$c_0=\bar e_0$,
		$\widetilde c_2=3L\varsigma_\star^2/(\eta^2\delta)$,
		$c_1=2\sigma^2/(\mu N)$, and
		$c_2=6Q_\lambda L\sigma^2/(\mu\delta)$.
		Moreover,
		\[
		\widetilde c_2\bar\alpha^3
		=
		\frac{24L\lambda^2\bar\alpha^3}
		{\Delta_\lambda^3}\varsigma_\star^2,
		\qquad
		\frac{Q_\lambda}{\delta}
		=
		\frac{34\lambda^2}{\Delta_\lambda^3}
		\leq
		\frac{34}{\Delta_\lambda^3}.
		\]
		Substituting these quantities into
		\eqref{eq:tuning_bound} proves
		\eqref{eq:tuned_rate_bigO}.
		
		When $\sigma^2=0$, the two noise terms in
		\eqref{eq:baseline_linear_bound_main} vanish, and taking
		$\alpha=\bar\alpha$ gives \eqref{eq:tuned_rate_bigO} directly.
	\end{proof}

	\bibliographystyle{IEEEtran}
	\bibliography{ref_gt_spectral}

\end{document}